\documentclass[onehalfspacing,gpscopy,11pt]{ubcdiss}

\usepackage{times,courier}
\usepackage[scaled=.92]{helvet}

\usepackage{pifont}

\usepackage{checkend}	
\usepackage{graphicx}	

\usepackage{booktabs}

\usepackage{listings}
\usepackage{color}
\definecolor{greytext}{gray}{0.5}

\usepackage{comment}

\usepackage[numbers,sort&compress]{natbib}

\usepackage[compact]{titlesec}
\titleformat*{\section}{\singlespacing\raggedright\bfseries\Large}
\titleformat*{\subsection}{\singlespacing\raggedright\bfseries\large}
\titleformat*{\subsubsection}{\singlespacing\raggedright\bfseries}
\titleformat*{\paragraph}{\singlespacing\raggedright\itshape}

\usepackage[format=hang,indention=-1cm,labelfont={bf},margin=1em]{caption}

\usepackage{url}
\ifgpscopy
  \usepackage[bookmarks,bookmarksnumbered,%
     pagebackref,linktocpage,%
     colorlinks=true,%
     linkcolor=black,%
     urlcolor=blue,%
     citecolor=black%
     ]{hyperref}
\else
  \usepackage[bookmarks,bookmarksnumbered,%
    pagebackref,linktocpage,%
    allbordercolors={0.8 0.8 0.8},%
    ]{hyperref}
\fi
\newcommand{\nocitations}{\relax}
\renewcommand*{\backrefalt}[4]{%
\textcolor{greytext}{\ifcase #1%
\nocitations%
\or
\(\rightarrow\) page #2%
\else
\(\rightarrow\) pages #2%
\fi}}

\usepackage{amsthm}
\usepackage{amscd}
\usepackage{amsfonts}         
\usepackage{amsmath}
\usepackage{amssymb}
\usepackage{tikz}
\usepackage{tikz-cd}
\usepackage{mathrsfs}

\DeclareUrlCommand\DOI{}
\newcommand{\doi}[1]{\href{http://dx.doi.org/#1}{\DOI{doi:#1}}}

\makeatletter
\newenvironment{epigraph}{%
	\begin{flushright}
	\begin{minipage}{\columnwidth-0.75in}
	\begin{flushright}
	\@ifundefined{singlespacing}{}{\singlespacing}%
    }{
	\end{flushright}
	\end{minipage}
	\end{flushright}}
\makeatother

\newcommand{\Ab}{\mathbb{A}}
\newcommand{\Bb}{\mathbb{B}}

\newcommand{\Hb}{\mathbb{H}}

\newcommand{\Lb}{\mathbb{L}}
\newcommand{\Mb}{\mathbb{M}}

\newcommand{\Pb}{\mathbb{P}}
\newcommand{\Qb}{\mathbb{Q}}

\newcommand{\Zb}{\mathbb{Z}}

\newcommand{\Cc}{\mathcal{C}}

\newcommand{\Ec}{\mathcal{E}}
\newcommand{\Fc}{\mathcal{F}}
\newcommand{\Gc}{\mathcal{G}}

\newcommand{\Ic}{\mathcal{I}}

\newcommand{\Kc}{\mathcal{K}}
\newcommand{\Lc}{\mathcal{L}}

\newcommand{\Nc}{\mathcal{N}}
\newcommand{\Oc}{\mathcal{O}}

\newcommand{\Qc}{\mathcal{Q}}
\newcommand{\Rc}{\mathcal{R}}
\newcommand{\Sc}{\mathcal{S}}

\newcommand{\Cs}{\mathscr{C}}

\newcommand{\Fs}{\mathscr{F}}

\newcommand{\Is}{\mathscr{I}}

\newcommand{\Ls}{\mathscr{L}}
\newcommand{\Ms}{\mathscr{M}}

\newcommand{\Ss}{\mathscr{S}}

\newcommand{\Pbf}{\mathbf{P}}

\newcommand{\bbf}{\mathbf{b}}

\newcommand{\xbf}{\mathbf{x}}

\newcommand{\PCob}{{\underline\Omega}}

\newcommand{\op}{\mathrm{op}}
\newcommand{\ch}{\mathrm{ch}}
\newcommand{\CH}{\mathrm{CH}}

\newcommand{\pr}{\mathrm{pr}}

\newcommand{\Spec}{\mathrm{Spec}}
\newcommand{\Gr}{\mathrm{Gr}}
\newcommand{\Fl}{\mathrm{Fl}}
\newcommand{\bl}{\mathrm{Bl}}
\newcommand{\Td}{\mathrm{Td}}
\newcommand{\wtil}{\widetilde}
\newcommand{\QCoh}{\mathrm{QCoh}}
\newcommand{\colim}{\mathrm{colim}}
\newcommand{\Perf}{\mathrm{Perf}}

\newcommand{\rank}{\mathrm{rank}}

\newcommand{\Fun}{\mathrm{Fun}}
\newcommand{\hook}{\hookrightarrow}
\newcommand{\cl}{\mathrm{cl}}
\newcommand{\modmod}{/\!\!/}

\newcommand{\Id}{\mathrm{Id}}

\newcommand{\xto}{\xrightarrow}
\newcommand{\dash}{{\text -}}
\newcommand{\Map}{\mathrm{Map}}
\newcommand{\Poly}{\mathrm{Poly}}
\newcommand{\dRing}{\mathrm{dRing}}
\newcommand{\dSch}{\mathrm{dSch}}
\newcommand{\inv}{\mathrm{inv}}
\newcommand{\naive}{\mathrm{naive}}
\newcommand{\univ}{\mathrm{univ}}
\newcommand{\pre}{\mathrm{pre}}
\newcommand{\snc}{\mathrm{snc}}

\newtheorem{theo}{Tplottin ubuntuheorem}
\theoremstyle{plain}
\newtheorem{thm}[theo]{Theorem}
\newtheorem{lem}[theo]{Lemma}
\newtheorem{prop}[theo]{Proposition}
\newtheorem{cor}[theo]{Corollary}

\newtheorem*{thm*}{Theorem}
\newtheorem*{lem*}{Lemma}
\newtheorem*{prop*}{Proposition}
\newtheorem*{cor*}{Corollary}

\theoremstyle{definition}
\newtheorem{defn}[theo]{Definition}
\newtheorem{ex}[theo]{Example}
\newtheorem{cons}[theo]{Construction}
\newtheorem{rem}[theo]{Remark}
\newtheorem{war}[theo]{Warning}

\title{Derived Algebraic Cobordism}

\author{Toni Annala}
\previousdegree{BSc., University of Helsinki, 2015}
\previousdegree{MSc., University of Helsinki, 2017}

\degreetitle{Doctor of Philosophy}

\institution{The University of British Columbia}
\campus{Vancouver}

\faculty{The Faculty of Graduate and Postdoctoral Studies}
\department{Mathematics}
\submissionmonth{March}
\submissionyear{2022}

\examiningcommittee{Kalle Karu, Mathematics, UBC}{Supervisor}
\examiningcommittee{Kai Behrend, Mathematics, UBC}%
    {Supervisory Committee Member}
\examiningcommittee{Ben Williams, Mathematics, UBC}{Supervisory Committee Member}
\examiningcommittee{Zinovy Reichstein, Mathematics, UBC}{University Examiner}
\examiningcommittee{Dragos Ghioca, Mathematics, UBC}{University Examiner}
\hypersetup{
  pdftitle={Derived Algebraic Cobordism},
  pdfauthor={Toni Annala},
  pdfkeywords={derived algebraic geometry, cobordism, intersection theory, $K$-theory}
}

\begin{document}


\maketitle

\makecommitteepage


\chapter{Abstract}

We construct and study a theory of bivariant cobordism of derived schemes. Our theory provides a vast generalization of the algebraic bordism theory of characteristic 0 algebraic schemes, constructed earlier by Levine and Morel, and a (partial) non-$\Ab^1$-invariant refinement of the motivic cohomology theory $MGL$ in Morel--Voevodsky's stable motivic homotopy theory. Our main result is that bivariant cobordism satisfies the projective bundle formula. As applications of this, we construct cobordism Chern classes of vector bundles, and establish a strong connection between the cobordism cohomology rings and the Grothendieck ring of vector bundles. We also provide several universal properties for our theory. Additionally, our algebraic cobordism is also used to construct a candidate for the elusive theory of Chow cohomology of schemes.

\cleardoublepage



\chapter{Lay Summary}

Algebraic geometry is the study of geometric objects---schemes---that admit a description as solution sets of systems of polynomial equations. Schemes are often studied via their invariants, such as their dimension and degree. In order to construct more refined invariants, one often considers the scheme itself as its own invariant, up to a suitably chosen notion of equivalence. A variation of this idea naturally leads one to consider algebraic cobordism, where, roughly speaking, two schemes are considered equivalent if they can be interpolated by a continuous one-dimensional family. We provide a vast generalization of the previous theories of algebraic cobordism using techniques originating from homotopy theory. Moreover, by computing the algebraic cobordism group of a ``projective bundle'' in terms of its ``base'', we are able to establish strong relationship of our invariant with other fundamental invariants in algebraic geometry.

\cleardoublepage


\chapter{Preface}

This thesis provides a streamlined account of part of the material from the articles \cite{annala-cob, annala-yokura, annala-chern, annala-pre-and-cob, annala-base-ind-cob}. The research of \cite{annala-yokura} was conducted under an equal collaboration with Professor Shoji Yokura from Kagoshima University, and the other articles were a product of independent research. Moreover, we have also included unpublished results from an equal collaboration with Ryomei Iwasa.

\begin{itemize}
\item Chapter \ref{ch:biv} introduces the reader to the bivariant formalism. The chapter is mostly expository, but contains some original contributions from \cite{annala-cob, annala-pre-and-cob}, such as the definition and the basic properties of bivariant ideals.

\item Chapter \ref{ch:dac} is an introduction to derived algebraic geometry. Section \ref{sect:dacbasics} is expository, whereas Sections \ref{sect:qproj} and \ref{sect:dci} contain some original material from \cite{annala-base-ind-cob} generalizing well-known classical facts to derived algebraic geometry. Section \ref{sect:blowpup} contains some original results from \cite{annala-cob,annala-chern,annala-yokura}.

\item Chapter \ref{ch:cob} presents a streamlined construction of bivariant algebraic cobordism, based on the articles \cite{annala-cob, annala-yokura, annala-chern, annala-base-ind-cob}.

\item Chapter \ref{ch:pbf} provides a streamlined proof for the projective bundle formula, and various applications, based on the articles \cite{annala-cob, annala-yokura, annala-chern, annala-base-ind-cob}.

\item Chapter \ref{ch:univprop} provides a universal property for bivariant algebraic cobordism, which is based on unpublished joint work with Ryomei Iwasa, related to the article \cite{annala-iwasa}.

\item In Chapter \ref{ch:LMcomp}, we prove that bivariant algebraic cobordism extends Levine--Morel's algebraic bordism, based on the articles \cite{annala-cob, annala-pre-and-cob}.

\item In Chapter \ref{ch:bivtwist}, we study the operation of twisting oriented bivariant theories, and its applications, based on the articles \cite{annala-cob, annala-chern, annala-base-ind-cob}.
\end{itemize}
\cleardoublepage

\tableofcontents
\cleardoublepage	

\textspacing		


\chapter{Acknowledgments}

First and foremost, I would like to thank my advisor Kalle Karu. Of the more than 20 graduate school applications I sent out, University of British Columbia was the only one that accepted me, so, without his choice to take me as a graduate student, I might not be doing mathematics. Moreover, the research problem he proposed---extending Levine--Morel's algebraic bordism to a bivariant theory---led to an interesting research program. Moreover, I have benefited from several discussions with several faculty members at the department, including Ben Williams and Kai Behrend.

Secondly, I would like to thank my colleagues, who have made me feel welcome to the community. I cannot overstate the psychological impact of meeting other people who are interested in your research. Special thanks go to Marc Hoyois, Ryomei Iwasa, Adeel Khan, Marc Levine, Denis Nardin, Gabriele Vezzosi and Shoji Yokura, with whom I have had several pleasant and productive conversations. 

When doing a doctorate in mathematics, it is easy to forget the other aspects of life. I would therefore like to thank Professor Mikko Möttönen who, during my national service in 2021, reminded me that, in real life, one has to consider physics as well. Joking aside, I would like to thank my family, a small circle of friends from high school, and the Finnish algorithm community from keeping me (relatively) sane during these years. I would also like to thank all the friends I have met in Vancouver, with special thanks to Abhishek for being an almost perfect roommate for two and a half years. 

The support from Vilho, Yrjö and Kalle Väisälä Foundation of the Finnish Academy of Science and letters has been instrumental in avoiding all teaching responsibilities for the past couple of years.


\mainmatter



\chapter{Introduction}
\label{ch:Introduction}

\begin{epigraph}
\emph{Ask, and it will be given to you; seek, and you will find; knock, and the door will be opened to you.} ---~Matthew 7:7
\end{epigraph}
%
%
%
In topology, complex cobordism is the universal cohomology theory with a complex orientation. In algebraic geometry, similar role should be played by \emph{algebraic cobordism}, but unfortunately the situation is less clear. First models of algebraic cobordism appeared as a piece of Voevodsky's original approach to the proof of the Milnor conjecture \cite{voevodsky:1996}, but later its role in the proof was suppressed \cite{voevodsky:2003a, voevodsky:2003b}. In this approach, algebraic cobordism is treated, rather abstractly, as the cohomology theory represented by the motivic spectrum $MGL$ in the motivic homotopy theory of Morel and Voevodsky \cite{morel:1999}. This approach has two major downsides. First of all, the theory produced this way is automatically $\Ab^1$-homotopy invariant, which should not be true for a genuine theory of algebraic cobordism. Secondly, it does not easily yield geometric models similar to those in topology \cite{quillen:1971} (however, there exists recent work in this direction \cite{EHKSY1, EHKSY2, EHKSY3}).

A geometric model for algebraic cobordism was constructed by Levine and Morel in their seminal treatise \cite{levine-morel}, and later simplified by Levine and Pandharipande \cite{levine-pandharipande}. What they actually construct is the corresponding homology theory, algebraic bordism, for finite type schemes over a field of characteristic 0. However, as in the case of Chow groups \cite{fulton:1998}, this theory restricts to a multiplicative cohomology theory on smooth schemes. The biggest drawback of this approach is its degree of generality: since Levine and Morel made liberal use of resolution of singularities \cite{hironaka:1964a, hironaka:1964b} and weak factorization \cite{abramovich:2002}, their methods are restricted to algebraic schemes over a field of characteristic 0. Still, this approach has led to several interesting applications \cite{levine-pandharipande,hudson:2017,hudson:2019,sechin:2018a, sechin:2018b}.

Derived algebraic geometry was introduced to the subject by Lowrey and Schürg in their innovative work \cite{lowrey--schurg}, providing an alternative construction of Levine--Morel's algebraic bordism, transforming the most intricate part of the original construction---the construction of \emph{Gysin pullbacks}---completely transparent. However, Lowrey and Schürg did not exploit the full potential of derived geometry in studying algebraic cobordism and related cohomology theories in algebraic geometry. In this thesis, I explain some of the results of my research program, the aim of which is to construct a workable and concrete theory of algebraic cobordism using derived geometry. A lot of material had to be left out: most importantly, we completely overlook the results of \cite{annala-spivak}, by which algebraic bordism groups over fields and nice enough discrete valuation rings are, after inverting the residual characteristic exponent in the coefficients, $\Ab^1$-homotopy invariant and generated by bordism classes of regular schemes; and the results of the joint article with Ryomei Iwasa \cite{annala-iwasa}, in which the first steps towards a geometric construction of non-$\Ab^1$-invariant higher algebraic cobordism are taken. Most of our attention will be focused on studying the universal precobordism theory, which is an intermediate step in the construction of algebraic cobordism. The associated cohomology theory of universal precobordism already satisfies many of the properties expected from algebraic cobordism, and this, together with the simplicity of its construction, justifies the attention we give to it.

\section*{Contents of the thesis} 

Chapters \ref{ch:biv} and \ref{ch:dac} are essentially background sections that introduce the reader to bivariant theories and derived algebraic geometry, respectively. In Chapter \ref{ch:cob}, we construct the universal precobordism as well as the algebraic cobordism of $S$-schemes, where $S$ is a nice enough fixed base scheme. The main result of this Chapter is that the Euler class of the tensor product of line bundles may be computed by the means of a formal group law. In Chapter \ref{ch:pbf}, we prove the main result of this thesis: universal precobordism satisfies the projective bundle formula. As an application, we construct Chern classes of vector bundles in universal precobordism, and prove that the Grothendieck ring of vector bundles can be recovered from the universal precobordism rings by ``enforcing the multiplicative formal group law''. In Chapter \ref{ch:univprop}, following unpublished joint work with Ryomei Iwasa, we prove that universal precobordism is the universal bivariant theory that satisfies a weak version of the projective bundle formula. We also provide similar universal properties for the associated homology and cohomology theories. In Chapter \ref{ch:LMcomp} we show that our bivariant algebraic cobordism extends Levine--Morel's algebraic bordism to a bivariant theory of quasi-projective derived algebraic $k$-schemes, where $k$ is a field of characteristic 0. Finally, in Chapter \ref{ch:bivtwist} we study universal precobordism and algebraic cobordism with rational coefficients, using the operation of twisting bivariant theories. We also establish a cohomological counterpart of Baum--Fulton--MacPherson's Grothendieck--Riemann--Roch theorem \cite{baum:1975}.

Throughout the thesis, we will freely use the language of $\infty$-categories (see e.g. \cite{HTT}). Moreover, a basic familiarity with the theory schemes (see e.g. \cite{hartshorne:1977}) and, to a much lesser extent, formal group laws (see e.g. \cite{hazewinkel:1978}), is assumed.

\chapter{Bivariant theories}\label{ch:biv}

In this preliminary chapter, we recall the bivariant formalism. Bivariant theories were originally introduced by Fulton and MacPherson to aid in the study of characteristic classes of singular spaces, and to unify several types of Riemann--Roch theorems \cite{fulton-macpherson}. Such a theory, assigning an Abelian group to each morphism in a category, contains both a homological and a cohomological theory, and the bivariant axioms ensure that these theories interact with each other in a well-behaved manner. The operational bivariant theories have attracted attention \cite{kimura:1992, fulton:1998, krishna:2012, anderson:2015, gonzalez:2015, annala:OBMToric}, especially in the context of toric varieties, as the operational cohomology rings of toric varieties tend to admit pleasant combinatorial descriptions. More relevant for us, however, is Yokura's work on universal bivariant theories \cite{yokura:2009, yokura:2019}, which is one of the main ingredients of our construction of bivariant algebraic cobordism, the other one being derived algebraic geometry. Unlike the operational theories, the universal bivariant theories admit a natural geometric description; moreover, they contain more refined information.

Our own contributions in the content of this chapter are rather modest. Mostly it consists on the choice of subject material, as well as some simple definitions (such as that of a bivariant ideal), which make using the bivariant formalism slightly more pleasant. In Section \ref{sect:bivform}, we recall the basics of bivariant theories. Section \ref{sect:univbiv} is dedicated to the study of universal bivariant theories, following the work of Yokura \cite{yokura:2009}. In Section \ref{sect:bivmanip}, we define two useful tools to manipulate bivariant theories, namely quotients by bivariant ideals and extensions of scalars. 

\section{Bivariant formalism}\label{sect:bivform}

In this section, we carefully lay out the bivariant formalism of Fulton--MacPherson \cite{fulton-macpherson}. We start by defining bivariant theories in Section \ref{ssect:bivdef}, after which we recall the construction of bivariant cross product in Section \ref{ssect:bivcrossprod}. A bivariant theory contains a homological and a cohomological theory whose basic properties are recalled in Section \ref{ssect:indthy}. One of the most interesting aspect of bivariant theories are the orientations they may admit. The theory of orientations of bivariant theories is recalled in Section \ref{ssect:gysin}, where it is also noted that the associated (co)homology theory of oriented bivariant theories satisfies many of the axioms of oriented (co)homology theories. In Section \ref{ssect:poincare} we recall the definition of strong orientations, which implies an analogue Poincaré duality for a subclass of objects. We also record categorical conditions, under which such statements automatically hold. Finally, in Section \ref{ssect:bivprop}, we study additive bivariant theories, as well as rings of coefficients for bivariant theories.

\subsection{Basic definitions}\label{ssect:bivdef}

Before defining bivariant theories, we need to fix the desired functoriality.

\begin{defn}\label{def:bivfunct}
A \emph{(bivariant) functoriality} is a tuple $\Fs = (\Cc, \Cs, \Is, \Ss)$, consisting of 
\begin{enumerate}
\item $\Cc$, an $\infty$-category with a final object $pt$ and all fiber products;

\item $\Cs$, a class of morphisms in $\Cc$ regarded as \emph{confined morphisms}, which contains all equivalences and is closed under compositions and pullbacks; moreover, if $f, f'$ are equivalent morphisms $X \to Y$ in $\Cc$, and if $f$ in $\Cs$, then so is $f'$: in short, $\Cs$ is \emph{closed under equivalences};

\item $\Is$, a class of Cartesian squares in $\Cc$ regarded as \emph{independent squares}, which contains all squares of form
\begin{center}
\begin{tikzcd}
X \arrow[]{r}{f} \arrow[]{d}{\mathrm{Id}_X} & Y \arrow[]{d}{\mathrm{Id}_Y} \\
X \arrow[]{r}{f} & Y
\end{tikzcd}
\ \ and \ \ \
\begin{tikzcd}
X \arrow[]{r}{\mathrm{Id}_X} \arrow[]{d}{f} & X \arrow[]{d}{f} \\
Y \arrow[]{r}{\mathrm{Id}_Y} & Y,
\end{tikzcd}
\end{center}
is closed under vertical and horizontal compositions, and is closed under equivalences of Cartesian squares;

\item $\Ss$, a class of morphisms in $\Cc$ regarded as \emph{specialized morphisms}, which contains all equivalences and is closed under compositions and equivalences.
\end{enumerate}
\end{defn}

Many of the classical bivariant functorialities (e.g. that of bivariant $K$-theory \cite{fulton-macpherson} or of operational bivariant theories \cite{fulton:1998,anderson:2015,gonzalez:2015}) are straightforward modifications of the following functoriality.

\begin{ex}\label{ex:bivfunctclass}
Let $\Cc$ be the category of schemes. Then the tuple $(\Cc, \Cs, \Is, \Ss)$ is a bivariant functoriality, where $\Cs$ consists of proper morphisms, $\Is$ consists either of all Cartesian squares or of all Tor-independent Cartesian squares, and $\Ss$ consists of formally lci morphisms. 
\end{ex}

On the other hand, the derived bivariant theories we are going to consider in this work are fairly straightforward modifications of the following functoriality. 

\begin{ex}\label{ex:bivfunct}
Let $\Cc$ be the $\infty$-category of derived schemes. Then the tuple $(\Cc, \Cs, \Is, \Ss)$ is a bivariant functoriality, where $\Cs$ consists of proper morphisms, $\Is$ consists of all Cartesian squares and $\Ss$ consists of formally quasi-smooth morphisms\footnote{Quasi-smooth morphism is a derived geometric analogue of an lci morphism.}. 
\end{ex}

We are now able to define bivariant theories with a fixed functoriality.

\begin{defn}\label{def:bivthy}
Let $\Fc$ be a bivariant functoriality. A \emph{bivariant theory $\Bb$ (with functoriality $\Fc$)} is the assignment of
\begin{enumerate}
\item an Abelian group $\Bb(X \xto{f} Y)$ to each morphisms $X \to Y$ in $\Cc$, depending only on the equivalence class of $f$ in $\Cc$; the morphism $f$ is often omitted from the notation for simplicity;

\item a \emph{bivariant pushforward} homomorphism $g_*: \Bb(X \xto{h \circ g} Y) \to \Bb(X' \xto{h} Y)$ to each factorization $f \simeq h \circ g$ such that $g$ is confined; the morphism $g_*$ depends only of the equivalence class of $g \in \Cc$;

\item a \emph{bivariant pullback} homomorphism $g^*: \Bb(X \stackrel{f}{\to} Y) \to \Bb(X' \xto{f'} Y')$ to each independent square
$$
\begin{tikzcd}
X' \arrow[]{r}{f'} \arrow[]{d}{g'} & Y' \arrow[]{d}{g} \\
X \arrow[]{r}{f} & Y,
\end{tikzcd}
$$
depending only on the equivalence class of the square, i.e., the equivalence class of $g \in \Cc$; 

\item a bi-additive \emph{bivariant product} $\bullet: \Bb(X \xto{f} Y) \times \Bb(Y \xto{g}  Z) \to \Bb(X \xto{g \circ f} Z)$ to each composition, depending only on the equivalence classes of $f$, $g$ and $g \circ f$.
\end{enumerate}
The above structure is required to satisfy the \emph{bivariant axioms}:

\begin{enumerate}
\item[($A_1$)] \emph{associativity of $\bullet$}: given morphisms $X \to Y \to Z \to W$ in $\Cc$ and bivariant elements $\alpha \in \Bb(X \to Y)$, $\beta \in \Bb(Y \to Z)$ and $\gamma \in \Bb(Z \to W)$, then
$$(\alpha \bullet \beta) \bullet \gamma = \alpha \bullet (\beta \bullet \gamma) \in \Bb(X \to W);$$

\item[($A_2$)] \emph{covariant functoriality of bivariant pushforward}: if $X \to Y$ factors through a composition $X \xrightarrow{f} X' \xrightarrow{g} X''$, where $f$ and $g$ are confined, then
$$g_* \circ f_*  = (g \circ f)_*;$$
moreover, $\Id_* = \Id$;

\item[($A_3$)] \emph{contravariant functoriality of bivariant pullback}: if  the Cartesian squares
$$
\begin{tikzcd}
X' \arrow[]{d} \arrow[]{r} & Y' \arrow[]{d}{g} \\
X \arrow[]{r} & Y
\end{tikzcd}
\text{ \ and \ }
\begin{tikzcd}
X'' \arrow[]{d} \arrow[]{r} & Y'' \arrow[]{d}{f} \\
X' \arrow[]{r} & Y'
\end{tikzcd}
$$
are independent, then
$$f^* \circ g^* = (g \circ f)^*;$$
moreover, $\Id^* = \Id$;

\item[($A_{12}$)] \emph{product and pushforward commute}: given $\alpha \in \Bb(X \to Y)$, $\beta \in \Bb(Y \to Z)$, and a confined morphism $f: X \to X'$ through which $X \to Y$ factors, then
$$f_*(\alpha \bullet \beta) = f_*(\alpha) \bullet \beta \in \Bb(X' \to Y);$$

\item[($A_{13}$)] \emph{bivariant pullback is multiplicative}: given $\alpha \in \Bb(X \to Y)$ and $\beta \in \Bb(Y \to Z)$ and suppose that the Cartesian squares
$$
\begin{tikzcd}
X' \arrow[]{d}{h''} \arrow[]{r} & Y' \arrow[]{d}{h'} \\
X \arrow[]{r} & Y
\end{tikzcd}
\text{ \ and \ }
\begin{tikzcd}
Y' \arrow[]{d}{h'} \arrow[]{r} & Z' \arrow[]{d}{h} \\
Y \arrow[]{r} & Z
\end{tikzcd}
$$ 
are independent\footnote{This implies that their horizontal composition is independent as well.}, then 
$$h^*(\alpha \bullet \beta) = h'^*(\alpha) \bullet h^*(\beta) \in \Bb(X' \to Z');$$

\item[($A_{23}$)] \emph{bivariant push-pull formula}: given $\alpha \in \Bb(X \to Z)$ and a Cartesian diagram
$$
\begin{tikzcd}
X' \arrow[]{d}{h''} \arrow[]{r}{f'} & Y' \arrow[]{d}{h'} \arrow[]{r} & Z' \arrow[]{d}{h} \\
X \arrow[]{r}{f} & Y \arrow[]{r} & Z
\end{tikzcd}
$$
such that $f$ is confined and 
$$
\begin{tikzcd}
X' \arrow[]{d}{h''} \arrow[]{r} & Z' \arrow[]{d}{h} \\
X \arrow[]{r} & Z
\end{tikzcd}
\text{ \ and \ }
\begin{tikzcd}
Y' \arrow[]{d}{h'} \arrow[]{r} & Z' \arrow[]{d}{h} \\
Y \arrow[]{r} & Z
\end{tikzcd}
$$ 
are independent, then
$$h^*(f_*(\alpha)) = f'_*(h^*(\alpha)) \in \Bb(Y' \to Z');$$

\item[($A_{123}$)] \emph{bivariant projection formula}: given a commutative diagram
$$
\begin{tikzcd}
X \arrow[]{d}{g'} \arrow[]{r} & Y \arrow[]{d}{g} \arrow[]{rd} & \\
X' \arrow[]{r} & Y' \arrow[]{r} & Z
\end{tikzcd}
$$
such that $g$ is confined and the small square on the left is independent, and elements $\alpha \in \Bb(Y \to Z)$ and $\beta \in \Bb(X' \to Y')$, then
$$g'_*(g^*(\beta) \bullet \alpha) = \beta \bullet g_*(\alpha) \in \Bb(X' \to Z);$$

\item[($U$)] \emph{existence of units}: for each $X \in \Cc$, there exists an unit\footnote{The unit is necessarily unique.} $1_X \in \Bb(X \to X)$ satisfying that for all $\alpha \in \Bb(X \to Y)$
$$1_X \bullet \alpha = \alpha \in \Bb(X \to Y),$$
and for all $\beta \in \Bb(Z \to X)$
$$\beta \bullet 1_X = \beta \in \Bb(Z \to X).$$
\end{enumerate}
\end{defn}

The majority of bivariant theories defined in geometric terms are commutative in the following sense.

\begin{defn}\label{def:bivcomm}
The class of independent squares is \emph{symmetric} if 
$$
\begin{tikzcd}
X' \arrow[]{d} \arrow[]{r} & Y' \arrow[]{d} \\
X \arrow[]{r} & Y
\end{tikzcd}
$$
being independent implies the independence of
$$
\begin{tikzcd}
X' \arrow[]{d} \arrow[]{r} & X \arrow[]{d} \\
Y' \arrow[]{r} & Y.
\end{tikzcd}
$$
A bivariant theory with a symmetric class of independent squares is \emph{commutative} if for all independent squares
$$
\begin{tikzcd}
X' \arrow[]{d} \arrow[]{r} & Y' \arrow[]{d}{g} \\
X \arrow[]{r}{f} & Y,
\end{tikzcd}
$$
and elements $\alpha \in \Bb(X \to Y)$ and $\beta \in \Bb(Y' \to Y)$, the equality
$$g^*(\alpha) \bullet \beta = f^*(\beta) \bullet \alpha$$
holds.
\end{defn}

Often, bivariant theories admit a useful (partial) grading.

\begin{defn}\label{def:bivgrad}
A \emph{grading} of a bivariant theory consists of isomorphisms 
$$\Bb(X \to Y) \cong \bigoplus_{i \in \Zb} \Bb^i(X \to Y),$$
depending only on the equivalence class of $X \to Y$ in $\Cc$, such that
\begin{enumerate}
\item bivariant pushforwards and pullbacks preserve the grading;
\item the bivariant product is graded, i.e., induces bi-additive pairings
$$\bullet: \Bb^i(X \to Y) \times \Bb^j(Y \to Z) \to \Bb^{i+j}(X \to Z).$$
\end{enumerate}
A \emph{graded bivariant theory} $\Bb^*$ is a bivariant theory $\Bb$ together with a specified grading.
\end{defn}

Unfortunately, bivariant algebraic cobordism seems not to admit a natural grading in general, but only for the bivariant groups associated to morphisms of finite type. This motivates the following definition.

\begin{defn}\label{def:bivpartgrad}
Let $\Bb$ be a bivariant theory with functoriality $\Fc = (\Cc, \Cs, \Is, \Ss)$, and let $\Cc'$ be a sub-$\infty$-category of $\Cc$ such that $\Cc'$ contains a final object of $\Cc$ and is closed under fiber products. Then, we define $\Fc' := (\Cc', \Cs', \Is', \Ss')$, where $\Cs', \Is'$ and $\Ss'$ consists of those morphisms and squares in $\Cs$, $\Is$ and $\Ss$ respectively, which are contained in $\Cc'$. Clearly the tuple $\Fc'$ is a bivariant functoriality and the bivariant theory $\Bb$ restricts to a bivariant theory $\Bb\vert_{\Fc'}$ on $\Fc'$.

A \emph{($\Fc'$)-partial grading} on $\Bb$ is a grading of $\Bb\vert_{\Fc'}$. A bivariant theory with a specified partial grading is referred to as being ($\Fc'$)-partially graded. Such a theory is denoted by $\Bb$; however, whenever a morphism $X \to Y$ belongs to $\Cc'$, the associated bivariant group is denoted by $\Bb^*(X \to Y)$, and it is considered as a graded Abelian group.
\end{defn}

Of course, in addition to bivariant theories, we will be interested in bivariant transformations between them.  

\begin{defn}\label{def:grotrans}
Let $\Bb$ and $\Bb'$ be bivariant theories with functoriality $\Fc$. Then a \emph{Grothendieck transformation} $\eta: \Bb \to \Bb'$ is a collection of homomorphisms
$$\eta_{X \to Y}: \Bb(X \to Y) \to \Bb'(X \to Y),$$
depending only on the equivalence class of $X \to Y$ in $\Cc$, such that 
\begin{enumerate}
\item if $X \to Y$ factors through a confined morphism $f: X \to X'$ and $\alpha \in \Bb(X \to Y)$, then
$$\eta_{X' \to Y}(f_*(\alpha)) = f_*(\eta_{X \to Y}(\alpha));$$

\item if the Cartesian square
$$
\begin{tikzcd}
X' \arrow[]{d}{g'} \arrow[]{r} & Y' \arrow[]{d}{g} \\
X \arrow[]{r} & Y
\end{tikzcd}
$$
is independent and $\alpha \in \Bb(X \to Y)$, then
$$\eta_{X' \to Y'}(g^*(\alpha)) = g^*(\eta_{X \to Y}(\alpha));$$

\item if $\alpha \in \Bb(X \to Y)$ and $\beta \in \Bb(Y \to Z)$, then
$$\eta_{X \to Z}(\alpha \bullet \beta) = \eta_{X \to Y}(\alpha) \bullet \eta_{Y \to Z}(\beta).$$
\end{enumerate}
For simplicity, the morphism $X \to Y$ is often omitted from the notation. A \emph{(partially) graded Grothendieck transformation} $\eta$ between (partially) graded bivariant theories is a Grothendieck transformation that preserves the (partial) grading.
\end{defn}

\subsection{Bivariant cross product}\label{ssect:bivcrossprod}

Before considering bivariant cross products, it is useful to restrict our attention to functorialities, in which independent squares satisfy the usual cancellation property of Cartesian squares.

\begin{defn}\label{def:sqcancelprop}
The class of independent squares has the \emph{cancellation property} if the independence of the large square and the bottom small square in
$$
\begin{tikzcd}
X' \arrow[]{r} \arrow[]{d} & X \arrow[]{d} \\
Y' \arrow[]{r} \arrow[]{d} & Y \arrow[]{d} \\
Z' \arrow[]{r} & Z
\end{tikzcd}
$$
imply the independence of the top small square, and if the independence of the large square and the rightmost small square in
$$
\begin{tikzcd}
X' \arrow[]{r} \arrow[]{d} & Y' \arrow[]{r} \arrow[]{d} & Z' \arrow[]{d} \\
X \arrow[]{r} & Y \arrow[]{r} & Z
\end{tikzcd}
$$
imply the independence of the leftmost small square.
\end{defn}

\begin{ex}\label{ex:torindcancel}
The class of Tor-independent squares in schemes satisfies the cancellation property.
\end{ex}

In order to define the cross product pairing $\times: \Bb(X_1 \to Y_1) \times \Bb(X_2 \to Y_2) \to \Bb(X_1 \times X_2 \to Y_1 \times Y_2)$, we form the following Cartesian diagram
$$
\begin{tikzcd}
X_1 \times X_2 \arrow[]{r} \arrow[]{d} & X_1 \times Y_2 \arrow[]{r} \arrow[]{d} & X_1  \arrow[]{d} \\
Y_1 \times X_2 \arrow[]{r}{g} \arrow[]{d} & Y_1 \times Y_2 \arrow[]{r}{f} \arrow[]{d}{h} & Y_1 \arrow[]{d} \\
X_2 \arrow[]{r} & Y_2 \arrow[]{r} & pt.
\end{tikzcd}
$$


\begin{defn}\label{def:bivcrossprod}
Let $\Bb$ be a bivariant theory. Moreover, suppose that each absolute product square
$$
\begin{tikzcd}
X \times Y \arrow[]{d} \arrow[]{r} & Y \arrow[]{d} \\
X \arrow[]{r} & pt
\end{tikzcd}
$$
is independent and that the independent squares satisfy the cancellation property. Then each square of the above diagram is independent, and the \emph{cross product} of $\alpha_1 \in \Bb(X_1 \to Y_1)$ and $\alpha_2 \in \Bb(X_2 \to Y_2)$ is be defined as 
$$\alpha_1 \times \alpha_2 := g^*(f^*(\alpha_1)) \bullet h^*(\alpha_2) \in \Bb(X_1 \times X_2 \to Y_1 \times Y_2).$$
\end{defn}

The cross product is compatible with pushforwards.

\begin{prop}\label{prop:crossprodpush}
Let $\alpha_1 \in \Bb(X_1 \to Y_1)$ and $\alpha_2 \in \Bb(X_2 \to Y_2)$, and suppose that $X_1 \to Y_1$ and $X_2 \to Y_2$ factor through confined morphisms $f_1: X_1 \to X
'_1$ and $f_2: X_2 \to X'_2$, respectively. Then $f_1 \times f_2: X_1 \times X_2 \to X'_1 \times X'_2$ is confined and the equality
$$(f_1 \times f_2)_*(\alpha_1 \times \alpha_2) = f_{1*}(\alpha_1) \times f_{2*}(\alpha_2)$$
holds.
\end{prop}
\begin{proof}
That $f_1 \times f_2$ is confined follows easily from the stability of confined morphisms in pullbacks and compositions. In order to prove the rest, we consider the Cartesian diagram
$$
\begin{tikzcd}
X_1 \times X_2 \arrow[]{r}{\Id_{X_1} \times f_2} \arrow[]{d}{f_1 \times \Id_{X_2}} & X_1 \times X'_2 \arrow[]{d} \arrow[]{r} & X_1 \times Y_2 \arrow[]{r} \arrow[]{d} & X_1  \arrow[]{d}{f_1} \\
X'_1 \times X_2 \arrow[]{r}{\Id_{X'_1} \times f_2} \arrow[]{d} & X'_1 \times X'_2 \arrow[]{d} \arrow[]{r} & X'_1 \times Y_2 \arrow[]{r} \arrow[]{d} & X'_1  \arrow[]{d} \\
Y_1 \times X_2 \arrow[]{r}{\Id_{Y_1} \times f_2} \arrow[]{d} & Y_1 \times X'_2 \arrow[]{d} \arrow[]{r}{g'} & Y_1 \times Y_2 \arrow[]{r}{f} \arrow[]{d}{h} & Y_1 \arrow[]{d} \\
X_2 \arrow[]{r}{f_2} & X'_2 \arrow[]{r} & Y_2 \arrow[]{r} & pt.
\end{tikzcd}
$$
In the notation of the above diagram, we compute that
\begin{align*}
&(f_1 \times f_2)_*(\alpha_1 \times \alpha_2) \\
&= (\Id_{X'_1} \times f_2)_* \bigg( (f_1 \times \Id_{X_2})_* \Big((\Id_{Y_1} \times f_2)^*\big(g'^*(f^*(\alpha_1))\big) \bullet h^*(\alpha_2)\Big)\bigg) \\
&= (\Id_{X'_1} \times f_2)_* \bigg( (f_1 \times \Id_{X_2})_* \Big((\Id_{Y_1} \times f_2)^*\big(g'^*(f^*(\alpha_1))\big)\Big) \bullet h^*(\alpha_2)\bigg) & (A_{12}) \\
&= (\Id_{X'_1} \times f_2)_* \bigg((\Id_{Y_1} \times f_2)^* \Big(g'^*\big(f^*(f_{1*}(\alpha_1))\big)\Big) \bullet h^*(\alpha_2)\bigg) & (A_{23}) \\
&= g'^*\big(f^*(f_{1*}(\alpha_1))\big) \bullet  (\Id_{Y_1} \times f_2)_*(h^*(\alpha_2)) & (A_{123}) \\
&= g'^*\big(f^*(f_{1*}(\alpha_1))\big) \bullet h^*(f_{2*}(\alpha_2)) & (A_{23}) \\
&= f_{1*}(\alpha_1) \times f_{2*}(\alpha_2),
\end{align*}
as desired.
\end{proof}

\subsection{The induced homology and cohomology theory}\label{ssect:indthy}

A bivariant theory contains both a homology theory and a cohomology theory as a part of it. Moreover, the bivariant axioms ensure that these theories interact in a well-behaved manner.

\begin{defn}\label{def:bivcoh}
Let $\Bb$ be a bivariant theory with functoriality $\Fc = (\Cc, \Cs, \Is, \Ss)$. Then the \emph{induced cohomology theory} $\Bb^\bullet$ assigns for every $X \in \Cc$ the ring 
$$\Bb^\bullet(X) := \Bb(X \xto{\Id} X),$$
the ring structure of which is given by the bivariant product. Bivariant pullbacks induce a contravariant functoriality along all morphisms in $\Cc$ for the rings $\Bb^\bullet(X)$. The induced cohomology theory of a graded bivariant theory inherits a grading, and is denoted by $\Bb^*$ instead of $\Bb^\bullet$. For partially graded bivariant theory, the induced cohomology rings are denoted by $\Bb^*(X)$ rather than $\Bb^\bullet(X)$ only in the case that they admit a grading.
\end{defn}

\begin{defn}\label{def:bivhom}
Let $\Bb$ be a bivariant theory with functoriality $\Fc = (\Cc, \Cs, \Is, \Ss)$. Then the \emph{induced homology theory} $\Bb_\bullet$ assigns, for every $X \in \Cc$, the group 
$$\Bb_\bullet(X) := \Bb(X \to pt).$$
The bivariant pushforwards induce a covariant functoriality along the confined morphisms of $\Cc$ for the groups $\Bb_\bullet(X)$. The induced homology theory of a graded bivariant theory inherits a grading with the convention that
$$\Bb_{i}(X) := \Bb^{-i}(X \to pt).$$ 
In this situation, the induced homology theory is denoted by $\Bb_*$ rather than $\Bb_\bullet$. For partially graded bivariant theory, the induced homology groups are denoted by $\Bb_*(X)$ rather than $\Bb_\bullet(X)$ only the case that they admit a grading.
\end{defn}

Let us then record some immediate results. First, the bivariant cross product induces a cross product for the associated cohomology and homology theories. For the cohomology theory, this product admits a pleasant alternative description.

\begin{prop}\label{prop:bihcohcrossprod}
Let $\Bb$ be a bivariant theory. Then the induced cross product on the associated cohomology rings satisfies the formula
$$\alpha_1 \times \alpha_2 = \pr_1^*(\alpha_1) \bullet \pr_2^*(\alpha_2) \in \Bb^\bullet(X \times Y),$$
where $\pr_1$ and $\pr_2$ are the projections $X \times Y \to X$ and $X \times Y \to Y$, respectively.
\end{prop} 
\begin{proof}
Follows immediately from the definition.
\end{proof}

The bivariant product makes the associated homology groups $\Bb_\bullet(X)$ left $\Bb^\bullet(X)$-modules. The bivariant axioms imply that a projection formula holds for this structure.

\begin{prop}\label{prop:cohhomprojform}
Suppose that $f: X \to Y$ is confined. Then, for all $\alpha \in \Bb^\bullet(Y)$ and $\beta \in \Bb_\bullet(X)$, the equation
$$f_*(f^*(\alpha) \bullet \beta) = \alpha \bullet f_*(\beta)$$
holds.
\end{prop}
\begin{proof}
Apply the bivariant axiom $(A_{123})$ to the diagram
$$
\begin{tikzcd}
X \arrow[]{d}{g} \arrow[]{r}{\Id} & X \arrow[]{d}{g} \arrow[]{rd} \\
Y \arrow[]{r}{\Id} & Y \arrow[]{r} & pt. 
\end{tikzcd}
$$
\end{proof}

\subsection{Orientations and Gysin morphisms}\label{ssect:gysin}

Orientations are one of the main reasons that make the study of bivariant theories worthwhile. They induce functoriality in the ``wrong direction'' for the associated homology and cohomology theory, generalizing, among other things, the lci pullback of Chow groups \cite{fulton:1998}. Moreover, Grothendieck transformations that do not conserve the orientation often lead to several Grothendieck--Riemann--Roch type formulas. 

\begin{defn}\label{def:bivor}
Let $\Bb$ be a bivariant theory. An \emph{orientation} of $\Bb$ is a choice of an element $\theta(f) \in \Bb(X \xto{f} Y)$ for each specialized map $f: X \to Y$, satisfying
\begin{enumerate}
\item $\theta(\Id) = 1_X \in \Bb(X \xto{\Id} X)$;
\item $\theta(g \circ f) = \theta(f) \bullet \theta(g) \in \Bb(X \to Z)$ whenever $f: X \to Y$ and $g: Y \to Z$ are specialized.
\end{enumerate}
An \emph{oriented bivariant theory}\footnote{Here, our terminology differs from that used by Yokura \cite{yokura:2009}. Yokura defines an oriented bivariant theory to be a bivariant theory together with special characteristic operators for objects in a category $\Lc$ fibered over $\Cc$ (e.g. $\Cc$ is a category of schemes and $\Lc$ is the category of line bundles over the schemes in $\Cc$). However, the notion of orientation we consider here seems more fundamental, because the special operators are essentially canonical Euler classes of line bundles, and these can be defined using the bivariant orientation.} is a bivariant theory together with a specified orientation. A Grothendieck transformation between oriented bivariant theories is \emph{orientation preserving} if it conserves the orientation.

If the class $\Ss$ of specialized morphisms is stable under independent pullbacks, then an orientation $\theta$ that is stable under independent pullbacks is referred to as being \emph{stable under pullbacks}. A bivariant theory equipped with such an orientation is called \emph{stably oriented}. 
\end{defn}

\begin{rem}\label{rem:stableor}
The derived bivariant theories we are going to consider in this work are stable under pullbacks. However, many classical bivariant theories, such as the operational bivariant theories \cite{fulton:1998,anderson:2015,gonzalez:2015}, do not have this property: there, the orientation is stable only under Tor-independent pullbacks. 
\end{rem}

Another interesting class of oriented bivariant theories is given by the centrally oriented theories, which generalize oriented commutative bivariant theories.

\begin{defn}\label{def:centralor}
An oriented bivariant theory $\Bb$ with a symmetric class of independent squares is \emph{centrally oriented} if, whenever
$$
\begin{tikzcd}
X' \arrow[d] \arrow[r] & Y' \arrow[d]{}{g} \\
X \arrow[r]{}{f} & Y
\end{tikzcd}
$$
is independent and $f$ is specialized, then, for all $\alpha\in \Bb(Y' \to Y)$, the equality
$$g^*(\theta(f)) \bullet \alpha = f^*(\alpha) \bullet \theta(f)$$
holds.
\end{defn}

\begin{ex}\label{ex:centralor}
The operational bivariant theories \cite{fulton:1998, anderson:2015, gonzalez:2015} are centrally oriented. However, it is not clear in which generality they are commutative and in any case, commutativity is not clear from the definition. 
\end{ex}

The orientation of a bivariant theory induces an interesting structure on the associated cohomology and homology theories. 

\begin{defn}\label{def:gysin}
Let $\Bb$ be an oriented bivariant theory. Then for each confined and specialized morphism $f: X \to Y$ and $\alpha \in \Bb^\bullet(X)$, one defines the \emph{Gysin pushforward} by
$$f_!(\alpha) := f_*(\alpha \bullet \theta(f)) \in \Bb^\bullet(Y),$$
and for each specialized morphism $g: X' \to Y'$ and $\beta \in \Bb_\bullet(Y')$, one defines the \emph{Gysin pullback} by
$$g^!(\beta) := \theta(g) \bullet \beta.$$
\end{defn}

The associated cohomology theory of an oriented bivariant theory has pleasant formal properties. These properties are analogous to a subset of the axioms of an oriented cohomology theory, used in Levine--Morel's treatment of algebraic cobordism (\cite{levine-morel} Definition 1.1.2).

\begin{prop}\label{prop:orcohformal}
Let $\Bb$ be an oriented bivariant theory. Then the induced cohomology theory $\Bb^\bullet$ satisfies the following properties:
\begin{enumerate}
\item \emph{covariant functoriality:} the Gysin pushforwards are functorial;
\item \emph{projection formula\footnote{The projection formula is equivalent to requiring $f_!$ to be a map of $\Bb^\bullet(Y)$-modules.}:} if $f: X \to Y$ is both specialized and confined, then, for all $\alpha \in \Bb^\bullet(Y)$ and $\beta \in \Bb^\bullet(X)$, the equality
$$f_!(f^*(\alpha) \bullet \beta) = \alpha \bullet f_!(\beta)$$
holds.
\end{enumerate}
Moreover, if $\Bb$ is stably oriented, then the associated cohomology theory satisfies the following further property:
\begin{enumerate}
\item[3.] \emph{push-pull formula:} if the square
$$
\begin{tikzcd}
X' \arrow[]{r}{g'}  \arrow[]{d}{f'} & Y' \arrow[]{d}{f} \\
X \arrow[]{r}{g} & Y
\end{tikzcd}
$$
is independent and if $f$ is both specialized and confined, then the equality
$$g^* \circ f_! = f'_! \circ g'^*: \Bb^\bullet(Y') \to \Bb^\bullet(X)$$
holds.
\end{enumerate}
\end{prop}
\begin{proof}
\begin{enumerate}
\item Clearly $\Id_! = \Id$. It is therefore enough to check that, for specialized and confined morphisms $f: X \to Y$ and $g: Y \to Z$, the equality $g_! \circ f_! = (g \circ f)_!$ holds. But this is straightforward: for any $\alpha \in \Bb^\bullet(X)$, we compute that
\begin{align*}
g_!(f_!(\alpha)) &= g_*\Big(f_*\big(\alpha \bullet \theta(f)\big) \bullet \theta(g)\Big) \\
&= g_*\Big(f_*\big(\alpha \bullet \theta(f) \bullet \theta(g) \big)\Big) & (A_{12})\\
&= (g \circ f)_*\big(\alpha \bullet \theta(g \circ f) \big) \\
&= (g \circ f)_!(\alpha), 
\end{align*}
as desired.

\item By a straightforward computation, the equality
\begin{align*}
f_!(f^*(\alpha) \bullet \beta) &= f_*\big(f^*(\alpha) \bullet \beta \bullet \theta(f)\big) \\
&= \alpha \bullet f_*\big(\beta \bullet \theta(f)\big) & (A_{123}) \\
&= \alpha \bullet f_!(\beta)
\end{align*}
holds, as desired. 

\item Let $\alpha \in \Bb^\bullet(Y')$. Then the equality
\begin{align*}
g^*(f_!(\alpha)) &= g^*\big(f_*(\alpha \bullet \theta(f))\big) \\
&= f'_*\big(g^*(\alpha \bullet \theta(f))\big) & (A_{23}) \\
&= f'_*\big(g'^*(\alpha) \bullet g^*(\theta(f))\big) & (A_{13}) \\
&= f'_*\big(g'^*(\alpha) \bullet \theta(f')\big) & (\text{stability}) \\
&= f'_!(g'^*(\alpha))
\end{align*}
holds, as desired.
\qedhere
\end{enumerate}
\end{proof}

Similar result holds for the associated homology theory. The properties are analogous to a subset of the axioms of an oriented Borel--Moore homology theory, used in Levine--Morel's treatment of algebraic cobordism (\cite{levine-morel} Definition 5.1.3).

\begin{prop}\label{prop:orhomformal}
Let $\Bb$ be an oriented bivariant theory. Then, the induced homology theory $\Bb_\bullet$ satisfies the following properties:
\begin{enumerate}
\item \emph{contravariant functoriality:} the Gysin pullbacks are functorial;

\item \emph{cross product is compatible with pushforward:} if all absolute product squares are independent and if independent squares satisfy the cancellation property, then, for all $\alpha \in \Bb_\bullet(X)$ and $\beta \in \Bb_\bullet(Y)$, and confined morphisms $f: X \to X'$ and $g: Y \to Y'$, 
$$(f \times g)_*(\alpha \times \beta) = f_*(\alpha) \times g_*(\beta);$$

\item \emph{cross product is compatible with Gysin pullback:} if $\Bb$ is centrally and stably oriented, if all absolute product squares are independent, and if the independent squares satisfy the cancellation property, then, for all $\alpha \in \Bb_\bullet(X)$ and $\beta \in \Bb_\bullet(Y)$, and for all specialized morphisms $f: X' \to X$ and $g: Y' \to Y$,
$$(f \times g)^!(\alpha \times \beta) = f^!(\alpha) \times g^!(\beta);$$

\item \emph{push-pull formula:} if $\Bb$ is stably oriented and the square
$$
\begin{tikzcd}
X' \arrow[]{r}{g'}  \arrow[]{d}{f'} & Y' \arrow[]{d}{f} \\
X \arrow[]{r}{g} & Y
\end{tikzcd}
$$
is independent, if $f$ is confined and if $g$ is specialized, then the equality
$$g^! \circ f_* = f'_* \circ g'^!: \Bb_\bullet(Y') \to \Bb_\bullet(X)$$
holds.
\end{enumerate}
\end{prop}
\begin{proof}
\begin{enumerate}
\item Follows immediately from the definition.

\item This is a special case of Proposition \ref{prop:crossprodpush}.

\item Consider the Cartesian diagram 
$$
\begin{tikzcd}
X' \times Y' \arrow[]{r}{\Id_{X'} \times g} \arrow[]{d} & X' \times Y \arrow[]{r} \arrow[]{d}{f \times \Id_Y} & X' \arrow[]{d}{f} \\
X \times Y' \arrow[]{r} \arrow[]{d} & X \times Y \arrow[]{r}{\pr_1} \arrow[]{d} & X \arrow[]{d}\\
Y' \arrow[]{r}{g} & Y \arrow[]{r}{\pi_Y} & pt.
\end{tikzcd}
$$
By definition, $\alpha \times \beta = \pi_Y^*(\alpha) \bullet \beta$, and
\begin{align*}
&(f \times g)^!(\alpha \times \beta) \\
&= (\Id_{X'} \times g)^!\Big((f \times \Id_Y)^! \big(\pi_Y^*(\alpha) \bullet \beta\big) \Big) \\
&= \theta(\Id_{X'} \times g) \bullet \theta(f \times \Id_Y) \bullet \pi^*_Y(\alpha) \bullet \beta \\
&= \theta(\Id_{X'} \times g) \bullet  \pi^*_Y(\theta(f) \bullet \alpha) \bullet \beta & (\text{stability, $A_{13}$}) \\
&= g^*(\pi^*_Y(\theta(f) \bullet \alpha)) \bullet \theta(g) \bullet \beta & (\text{stability, centrality}) \\
&= \pi^*_{Y'}(f^!(\alpha)) \bullet g^!(\beta) \\
&= f^!(\alpha) \times g^!(\beta)
\end{align*}
holds, as desired.

\item Suppose $\alpha \in \Bb_\bullet(Y')$. Then, 
\begin{align*}
g^!(f_*(\alpha)) &= \theta(g) \bullet f_*(\alpha) \\
&= f'_* \big(f^*(\theta(g)) \bullet \alpha\big) & (A_{123}) \\
&= f'_* \big(\theta(g') \bullet \alpha\big)  & (\text{stability}) \\
&= f'_*(g'^!(\alpha)), 
\end{align*}
 as desired. \qedhere
\end{enumerate}
\end{proof}

\subsection{Strong orientations and Poincaré duality}\label{ssect:poincare}

Often, in oriented bivariant theories that naturally appear in algebraic geometry, the   orientations along smooth morphisms have the following property, which gives rise to phenomena analogous to Poincaré duality.

\begin{defn}\label{def:strongor}
Let $\Bb$ be an oriented bivariant theory and $f: Y \to Z$ a specialized morphism. The orientation of $f$ is \emph{strong} if for all $X \to Y$, the map
$$- \bullet \theta(f): \Bb(X \to Y) \to \Bb(X \to Z)$$
is an isomorphism.
\end{defn}

In particular, if $\pi_X: X \to pt$ is specialized and has strong orientation, then the above morphism, ``capping by the fundamental class'', gives rise to an isomorphism $\Bb^\bullet(X) \to \Bb_\bullet(X)$. In order to obtain stronger results, it is useful to consider the following definition.

\begin{defn}\label{def:veryspecial}
A subclass $\Ss_v \subset \Ss$ is called a class of \emph{very special morphisms}\footnote{In \cite{annala-cob}, the terminology \emph{specialized projection} was used instead.} if
\begin{enumerate}
\item $\Ss_v$ is stable under compositions and pullbacks;
\item if $f: X \to Y$ is very specialized, then all Cartesian squares of form
$$
\begin{tikzcd}
X' \arrow[]{d} \arrow[]{r} & X \arrow[]{d}{f} \\
Y' \arrow[]{r} & Y
\end{tikzcd}
$$
are independent;
\item if $X \to Y$ is very specialized, and if $s$ is a morphism $Y \to X$ such that $f \circ s = \Id$, then $s$ is specialized.
\end{enumerate}
\end{defn}

The following result is just the formalization of \cite{fulton:1998} Proposition 17.4.2.

\begin{prop}\label{prop:veryspecialstrong}
Let $\Bb$ be a centrally and stably oriented bivariant theory such that the independent squares satisfy the cancellation property (Definition \ref{def:sqcancelprop}). Let $\Ss_v \subset \Ss$ be a class of very specialized morphisms. Then, the orientation of $\Bb$ is strong along morphisms in $\Ss_v$.
\end{prop} 
\begin{proof}
Let $g: Y \to Z$ be a very specialized morphism, and let $f: X \to Y$ be an arbitrary morphisms. Then all squares in the Cartesian diagram
$$
\begin{tikzcd}
X \arrow[]{r}{f} \arrow[]{d}{\Gamma_f} & Y \arrow[]{d}{\Delta} \\
X \times_Z Y \arrow[]{r} \arrow[]{d}{g''} & Y \times_Z Y \arrow[]{d}{g'} \arrow[]{r}{g'''} & Y \arrow[]{d}{g} \\
X \arrow[]{r}{f} & Y \arrow[]{r}{g} & Z
\end{tikzcd}
$$
are independent. We claim that the inverse to right multiplication by $\theta(g)$ is given by $\theta(\Gamma_f) \bullet g^*(-)$. Indeed, for any $\alpha \in \Bb(X \to Y)$, 
\begin{align*}
\theta(\Gamma_f) \bullet g^*(\alpha \bullet \theta(g)) &= \theta(\Gamma_f) \bullet g'^*(\alpha) \bullet \theta(g''') & (\text{$A_{13}$, stability}) \\
&= \alpha \bullet \theta(\Delta) \bullet \theta(g''') & (\text{centrality, stability}) \\
&= \alpha
\end{align*}
and, for any $\beta \in \Bb(X \to Z)$,
\begin{align*}
\theta(\Gamma_f) \bullet g^*(\beta) \bullet \theta(g)  &= \theta(\Gamma_f) \bullet \theta(g'') \bullet \beta & (\text{centrality, stability}) \\
&= \beta,
\end{align*}
as desired.
\end{proof}

If a class of very specialized morphisms is chosen, then $X \in \Cc$ is a \emph{very specialized object} if the structure morphism $\pi_X: X \to pt$ is very specialized. The homology group of such an object admits an intersection product.

\begin{defn}\label{def:intersectionprod}
Let $X$ be a very specialized object. Then $\Delta: X \to X \times X$ is specialized, and one defines the \emph{intersection product}
$$\frown: \Bb_\bullet(X) \times \Bb_\bullet(X) \to \Bb_\bullet(X)$$
as 
$$\alpha \frown \beta := \Delta^!(\alpha \times \beta).$$
\end{defn}

The intersection product equips the homology group of a very specialized object by a ring structure. The following analogue of Poincaré duality identifies it with the associated cohomology ring of $X$.

\begin{prop}\label{prop:poincare}
Let $\Bb$ be a centrally and stably oriented bivariant theory such that all absolute product squares are independent and the independent squares satisfy the cancellation property (Definition \ref{def:sqcancelprop}), and let $X$ be a very specialized object. Then, $- \bullet \theta(\pi_X): \Bb^\bullet(X) \to \Bb_\bullet(X)$ is an isomorphism of Abelian groups, and for all $\alpha, \beta \in \Bb^\bullet(X)$, 
$$(\alpha \bullet \beta) \bullet \theta(\pi_X) = (\alpha \bullet \theta(\pi_X)) \frown (\beta \bullet \theta(\pi_X)).$$
In particular, the intersection product $\frown$ is associative.
\end{prop}
\begin{proof}
The first claim follows immediately from Proposition \ref{prop:veryspecialstrong}. To prove the second claim, we consider the Cartesian diagram
$$
\begin{tikzcd}
X \times X \arrow[]{r}{\pr_1} \arrow[]{d}{\pr_2} & X \arrow[]{d}{\pi_X}\\
X \arrow[]{r}{\pi_X}& pt.
\end{tikzcd}
$$
By the proof of Proposition \ref{prop:veryspecialstrong}, the equality
$$\theta(\Delta)\bullet \pi_X^*(\alpha \bullet \theta(\pi_X)) = \alpha$$
holds, and therefore
\begin{align*}
(\alpha \bullet \theta(\pi_X)) \frown (\beta \bullet \theta(\pi_X)) &= \Delta^!\big((\alpha \bullet \theta(\pi_X)) \times (\beta \bullet \theta(\pi_X))\big) \\
&= \theta(\Delta)\bullet \pi_X^*(\alpha \bullet \theta(\pi_X)) \bullet \beta \bullet \theta(\pi_X) \\
&=\alpha \bullet \beta \bullet \theta(\pi_X),
\end{align*}
as desired.
\end{proof}

\subsection{Further properties of bivariant theories}\label{ssect:bivprop}

In this section, we recall miscellaneous properties of bivariant theories that will be relevant to our work. 

\subsubsection{Additivity}\label{sssect:add}

Here, we study additive bivariant theories. Note that our treatment differs slightly from that of Yokura \cite{yokura:2009}. We start by fixing kind of bivariant theories, where the notion of additivity is sensible.

\begin{defn}\label{def:goodcoprod}
Let $\Fc = (\Cc, \Cs, \Is, \Ss)$ be a bivariant functoriality. Then $\Fc$ has \emph{good coproducts} if 
\begin{enumerate}
\item the $\infty$-category $\Cc$ has all finite coproducts;
\item for all non-initial objects $X \in \Cc$, the mapping space $\Map_\Cc(X, \emptyset)$ is empty, where $\emptyset$ is an initial object;
\item a morphism $f: X \sqcup Y \to Z$ is confined (specialized) if and only if the morphisms $f \vert_X: X \to Z$ and $f \vert_Y: Y \to Z$ are confined (specialized);
\item all squares of form
$$
\begin{tikzcd}
\emptyset \arrow[]{r} \arrow[]{d} & X \arrow[]{d}{\iota_1} \\
Y \arrow[]{r}{\iota_2} & X \sqcup Y
\end{tikzcd}
\text{ \ and \ }
\begin{tikzcd}
X \times_Z Z' \sqcup Y \times_Z Z' \arrow[]{r} \arrow[]{d} & Z' \arrow[]{d}\\
X \sqcup Y \arrow[]{r} & Z
\end{tikzcd}
$$
are Cartesian, where $\iota_i$ are canonical inclusions;
\item all Cartesian squares of form
$$
\begin{tikzcd}
V' \arrow[r] \arrow[d] &  X \arrow[d]{}{\iota} \\
V \arrow[r] & X \sqcup Y 
\end{tikzcd}
\text{ \ and \ }
\begin{tikzcd}
V' \arrow[r] \arrow[d] &  V \arrow[d]\\
X \arrow[r]{}{\iota}  & X \sqcup Y 
\end{tikzcd}
$$
are independent.
\end{enumerate}
\end{defn}

Several other properties follow from the above assumptions.

\begin{lem}\label{lem:goodcoprodprop}
Let $\Fc$ be a bivariant functoriality with good coproducts. Then 
\begin{enumerate}
\item canonical inclusions are stable under pullbacks;
\item a square of form 
$$
\begin{tikzcd}
X \arrow[]{r}{\Id} \arrow[]{d}{\Id} & X \arrow[]{d}{\iota_1} \\
X \arrow[]{r}{\iota_1} & X \sqcup Y
\end{tikzcd}
$$
is Cartesian.
\end{enumerate}
\end{lem}
\begin{proof}
\begin{enumerate}
\item As $V \simeq \big(V \times_{X \sqcup Y} X \big) \sqcup \big(V \times_{X \sqcup Y} Y \big)$, the map $V \times_{X \sqcup Y} X \to V$ is a canonical inclusion.

\item As $X \simeq \big(X \times_{X \sqcup Y} X \big) \sqcup \big(Y \times_{X \sqcup Y} X \big)$, and as $Y \times_{X \sqcup Y} X \simeq \emptyset$, it follows that the natural map $X \times_{X \sqcup Y} X \to X$ is an equivalence. \qedhere
\end{enumerate}
\end{proof}

The following further ingredients are needed for a sensible notion of additivity.

\begin{defn}\label{def:wellor}
Let $\Bb$ be an oriented bivariant theory with good coproducts. Then $\Bb$ is \emph{well oriented} if the orientations along canonical inclusions are stable under (independent) pullbacks.
\end{defn}

We are now ready to define additive bivariant theories.

\begin{defn}\label{def:add}
Let $\Bb$ be a well oriented bivariant theory. Then $\Bb$ is \emph{additive}, if, for all $X \simeq \coprod_{i=1}^r X_i$, the equation
$$\sum_{i = 1}^r \iota_{i_*}(\theta(\iota_i)) = 1_X$$
holds, where $\iota_i$ are the canonical inclusions $X_i \to X$.
\end{defn}

Additive theories have many useful properties. In particular, the associated homology and cohomology groups of an additive bivariant theory are additive functors in the usual sense. 

\begin{prop}\label{prop:add}
Let $\Bb$ be an additive bivariant theory. Then
\begin{enumerate}
\item for all $Y \in \Cc$, the group $\Bb(\emptyset \to Y)$ is empty;

\item for all $X_1 \sqcup X_2 \to Y$, the map
$$
\iota_* := 
\begin{bmatrix}
\iota_{1*} & \iota_{2*}
\end{bmatrix}
:\Bb(X_1 \to Y) \oplus \Bb(X_2 \to Y) \to \Bb(X_1 \sqcup X_2 \to Y)
$$
is an isomorphism, the inverse of which is given by
\begin{align*}
\theta(\iota) \bullet - := 
\begin{bmatrix}
\theta(\iota_1) \bullet -  \\ 
\theta(\iota_2) \bullet -
\end{bmatrix} &: \Bb(X_1 \sqcup X_2 \to Y) \\
&\to \Bb(X_1 \to Y) \oplus \Bb(X_2 \to Y);
\end{align*}

\item for all $X \to Y_1 \sqcup Y_2$, the map
$$
\iota^* :=
\begin{bmatrix}
\iota^*_1 \\ 
\iota^*_2 
\end{bmatrix}
: \Bb(X \to Y_1 \sqcup Y_2) \to \Bb(X_1 \to Y_1) \oplus \Bb(X_2 \to Y_2)
$$
is an isomorphism, the inverse of which is given by
\begin{align*}
\begin{bmatrix}
\iota'_{1*}(- \bullet \theta(\iota_1)) & \iota'_{2*}(- \bullet \theta(\iota_2))
\end{bmatrix}
& : \Bb(X_1 \to Y_1) \oplus \Bb(X_2 \to Y_2) \\
& \to \Bb(X \to Y_1 \sqcup Y_2),
\end{align*}
where $X_i$ and $\iota'_i$ are as in the Cartesian diagram
$$
\begin{tikzcd}
X_i \arrow[d]{}{\iota'_i} \arrow[r] & Y_i \arrow[d]{}{\iota_i} \\
X \arrow[r] & Y_1 \sqcup Y_2.
\end{tikzcd}
$$
\end{enumerate}
\end{prop}
\begin{proof}
\begin{enumerate}
\item By assumption $1_\emptyset$ equals an empty sum, i.e., it is zero. Hence, for any $\alpha \in \Bb(\emptyset \to Y)$, we have that $\alpha = 1_\emptyset \bullet \alpha = 0$, proving that $\Bb(\emptyset \to Y) \cong 0$.

\item Considering the commutative diagram
$$
\begin{tikzcd}
X_{ij} \arrow[d] \arrow[r] & X_i \arrow[d]{} \arrow[rd] \\
X_j \arrow[r] & X_1 \sqcup X_2 \arrow[r] & Y,
\end{tikzcd}
$$
where the left square is an independent Cartesian square, one deduces from the bivariant axiom $(A_{123})$ and the well orientedness of $\Bb$, that for all $\alpha \in \Bb(X_i \to Y)$,
$$
\theta(\iota_j) \bullet \iota_{i*} (\alpha) = 
\begin{cases}
\alpha & \text{if $i = j$;} \\
0 & \text{if $i \not = j$.}
\end{cases}
$$
Hence $\theta(\iota) \bullet \iota_* (-)$ is the identity morphism. To conclude that $ \iota_*(\theta(\iota) \bullet -)$ is the identity, we compute that, for every $\beta \in \Bb(X_1 \sqcup X_2 \to Y)$,
\begin{align*}
\sum_{i=1}^2 \iota_{i*} (\theta(\iota_i) \bullet \beta) &=  \sum_{i=1}^2 \iota_{i*} (\theta(\iota_i)) \bullet \beta & (A_{12}) \\
&= \beta.
\end{align*}

\item Considering the independent Cartesian diagram
$$
\begin{tikzcd}
X_{ij} \arrow[d] \arrow[r] & X_j \arrow[d]{}{\iota'_j} \arrow[r] & Y_j \arrow[d]{}{\iota_j} \\
X_i \arrow[r]{}{\iota'_i} & X \arrow[r] & Y_1 \sqcup Y_2 
\end{tikzcd}
$$
we deduce from the bivariant axiom $(A_{23})$, that, for all $\alpha \in \Bb(X_i \to Y_i)$,
$$
\iota^*_i\big(\iota'_j(\alpha \bullet \iota(\iota_j))\big) =
\begin{cases}
\iota^*_i(\alpha \bullet \theta(\iota_i)) & \text{if $i = j$;} \\
0 & \text{if $i \not = j$.}
\end{cases} 
$$
Moreover, considering the independent Cartesian diagram
$$
\begin{tikzcd}
X_i \arrow[d]{}{\Id} \arrow[r] & Y_i \arrow[d]{}{\Id} \arrow[r]{}{\Id} & Y_i \arrow[d]{}{\iota_i} \\
X_i \arrow[r] & Y_i \arrow[r]{}{\iota_i} & Y_1 \sqcup Y_2,
\end{tikzcd}
$$
it follows from the bivariant axiom $(A_{13})$ and the well orientedness of $\Bb$ that $\iota^*_i(\alpha \bullet \theta(\iota_i)) = \alpha$. Hence the composition 
$$\Bb(X_1 \to Y_1) \oplus \Bb(X_2 \to Y_2) \to \Bb(X_1 \to Y_1) \oplus \Bb(X_2 \to Y_2)$$
is the identity. To verify that the other composition is the identity as well, we compute that, for every $\beta \in \Bb(X \to Y_1 \sqcup Y_2)$,
\begin{align*}
\sum_{i=1}^2 \iota'_{i*}\big( \iota_i^*(\beta) \bullet \theta(\iota_i) \big) &= \sum_{i=1}^2 \beta \bullet \iota_{i*}(\theta(\iota_i)) & (A_{123}) \\
&= \beta,
\end{align*}
as desired.
\qedhere
\end{enumerate}
\end{proof}

%

\subsubsection{Rings of coefficients}\label{sssect:bivcoeff}

It is often useful to consider bivariant theories taking values in modules over a ring, rather than just Abelian groups.

\begin{defn}\label{def:bivcoeff}
A \emph{ring of coefficients} for a bivariant theory $\Bb$ is a ring $R$ together with a map of rings $R \to \Bb^\bullet(pt)$. Such a bivariant theory $\Bb$ is referred to as being \emph{$R$-linear}.
\end{defn}

The terminology is explained by the following properties of $R$-linear theories.

\begin{prop}\label{prop:rlinearthyprop}
Let $\Bb$ be an $R$-linear bivariant theory. Then 
\begin{enumerate}
\item $\Bb(X \to Y)$ is an $R \dash R$-bimodule with the actions defined as
$$r \alpha := \pi_X^*(r) \bullet \alpha$$
and  
$$\alpha r := \alpha \bullet \pi^*_{Y}(r);$$

\item bivariant pullbacks and pushforwards are maps of $R \dash R$-bimodules;

\item bivariant products are $R$-balanced and they induce maps of $R \dash R$-bimodules
$$\Bb(X \to Y) \otimes_R \Bb(Y \to Z) \to \Bb(X \to Z).$$
\end{enumerate}
Moreover, if $\Bb$ is oriented, then
\begin{enumerate}
\item[4.] the Gysin-morphisms are maps of $R \dash R$-bimodules.
\end{enumerate}
\end{prop}
\begin{proof}
The fourth claim is an immediate consequence of the first three, so we will only prove them.
\begin{enumerate}
\item Since the bivariant product $\bullet$ is $\Zb$-bilinear, it follows that $\Bb(X \to Y)$ is both a left and a right $R$-module. The equality
$$r(\alpha s) = (r \alpha) s$$
follows from the associativity of $\bullet$.

\item Suppose that $X \to Y$ factors through a confined morphism $f: X \to X'$, and let $\alpha \in \Bb(X \to Y)$. Then,
\begin{align*}
f_*(s \alpha) &= f_* (\pi_X^*(s) \bullet \alpha) \\
&= f_*(f^*(\pi_{X'}^*(s)) \bullet \alpha) \\
&= \pi_{X'}^*(s) \bullet f_*(\alpha) & (A_{123}) \\
&= s f_*(\alpha),
\end{align*} 
and 
\begin{align*}
f_*(\alpha s) &= f_*(\alpha \bullet \pi_Y^*(s)) \\
&= f_*(\alpha) \bullet \pi^*_Y(s) & (A_{12}) \\
&= f_*(\alpha)s,
\end{align*}
so bivariant pushforwards are maps of $R \dash R$-bimodules.

Consider then a morphism $g: Y' \to Y$ and form the Cartesian square
$$
\begin{tikzcd}
X' \arrow[]{r} \arrow[]{d}{g'} & Y' \arrow[]{d}{g} \\
X \arrow[]{r} & Y.
\end{tikzcd}
$$ 
Then, given $\alpha \in \Bb(X \to Y)$, we compute that 
\begin{align*}
g^*(r \alpha) &= g^*(\pi_X^*(r) \bullet \alpha) \\
&= g'^*(\pi_X^*(r)) \bullet g^*(\alpha) & (A_{13}) \\
&= \pi^*_{X'}(r) \bullet g^*(\alpha) \\
&= r g^*(\alpha)
\end{align*}
and similarly one shows that $g^*(\alpha r) = g^*(\alpha) r$, so bivariant pullbacks are maps of $R \dash R$-bimodules as well.

\item Let us have $\alpha \in \Bb(X \to Y)$ and $\beta \in \Bb(Y \to Z)$. Then
\begin{align*}
(\alpha r) \bullet \beta &= (\alpha \bullet \pi^*_Y(r)) \bullet \beta \\
&= \alpha \bullet (\pi^*_Y(r) \bullet \beta)  \\
&= \alpha \bullet (r \beta),
\end{align*}
i.e., $\bullet$ is $R$-balanced. Moreover, the induced morphism
$$\Bb(X \to Y) \otimes_R \Bb(Y \to Z) \to \Bb(X \to Z)$$
is $R \dash R$-linear since
\begin{align*}
(r \alpha) \otimes \beta &\mapsto (r \alpha) \bullet \beta \\
&= (\pi^*_X(r) \bullet \alpha) \bullet \beta \\
&= \pi^*_X(r) \bullet (\alpha \bullet \beta)  \\
&= r(\alpha \bullet \beta)
\end{align*}
and similarly one shows that $\alpha \bullet (\beta r) \mapsto (\alpha \bullet \beta) r$. \qedhere
\end{enumerate}
\end{proof}

Of course, if $\Bb$ is a commutative bivariant theory and $R$ is a commutative ring, then the structure of a bimodule is extraneous: in such a case, it is simply the case that $\Bb$ takes values in $R$-modules and all the relevant structure morphisms are $R$-linear. 

There exists a natural notion of transformation of $R$-linear bivariant theories.

\begin{defn}\label{def:rlineartrans}
A Grothendieck transformation $\eta$ from an $R$-linear bivariant theory $\Bb_1$ to an $R$-linear bivariant theory $\Bb_2$ is \emph{$R$-linear} if the induced map
$$\eta: \Bb^\bullet(pt) \to \Bb^\bullet(pt)$$
is a map of $R$-algebras. 
\end{defn}

\begin{prop}\label{prop:rlineartrans}
A Grothendieck transformation $\eta: \Bb_1 \to \Bb_2$ between $R$-linear theories is $R$-linear if and only if all the induced morphisms 
$$\eta: \Bb_1(X \to Y) \to \Bb_2(X \to Y)$$
are maps of $R \dash R$-bimodules.
\end{prop}
\begin{proof}
Clearly the condition above is at least as strong as $R$-linearity. To prove that the two conditions coincide, we consider an $R$-linear Grothendieck transformation $\eta$, and $\alpha \in \Bb_1(X \to Y)$, and compute that
\begin{align*}
\eta (r \alpha) &= \eta(\pi^*_X(r) \bullet \alpha) \\
&= \eta(\pi^*_X(r)) \bullet \eta(\alpha) \\
&= \pi^*_X(\eta(r)) \bullet \eta(\alpha) \\
&= \pi^*_X(r) \bullet \eta(\alpha) \\
&= r \eta(\alpha).
\end{align*} 
One proves that $\eta(\alpha r) = \eta(\alpha)r$ in a similar fashion.
\end{proof}

\section{Universal bivariant theories}\label{sect:univbiv}

In this section, we recall Yokura's construction of the universal bivariant theory \cite{yokura:2009}. In addition, we consider a variant of this theory, namely the universal additive bivariant theory, which provides a convenient basis for our construction of bivariant algebraic cobordism.

\subsection{Universal stably oriented bivariant theory}\label{ssect:univbivthy}

Here, we recall the construction and the universal property of Yokura's universal bivariant theory. Throughout this subsection, $\Fc = (\Cc, \Cs, \Is, \Ss)$ is a bivariant functoriality in which all Cartesian squares of form
$$
\begin{tikzcd}
X' \arrow[]{d} \arrow[]{r} & Y' \arrow[]{d}{f} \\
X \arrow[]{r} & Y
\end{tikzcd}
\text{ \ and \ }
\begin{tikzcd}
X' \arrow[]{d} \arrow[]{r} & X \arrow[]{d} \\
Y' \arrow[]{r}{f} & Y,
\end{tikzcd}
$$
such that $f$ is confined, are independent ($\Fc$ satisfies \emph{$\Cs$-independence}), the class of independent squares is symmetric, and in which specialized morphisms are stable under independent pullbacks.

\begin{defn}\label{def:univbivthy}
Let $\Fc$ be as above. Define $\Mb_\Fc(X \to Y)$ as the free Abelian group on equivalence classes $[V \to X]$ of confined morphisms $V \to X$ such that the composition $V \to Y$ is specialized. One defines the bivariant operations as follows:
\begin{enumerate}
\item \emph{bivariant pushforward}: if $X \to  Y$ factors through a confined map $g: X \to X'$, then the bivariant pushforward is defined by linearly extending the formula
$$g_*([V \xto{f} X]) := [V \xto{g \circ f} X'];$$

\item \emph{bivariant pullback}: if
$$
\begin{tikzcd}
X' \arrow[]{d}{g'} \arrow[]{r} & Y' \arrow[]{d}{g} \\
X \arrow[]{r} & Y
\end{tikzcd}
$$
is an independent square, then the bivariant pullback is defined by linearly extending the formula
$$g^*([V \xto f X]) := [V' \xto{f'} X'],$$
where $f'$ is the pullback of $f$ along $g'$; the element $[V' \xto{f'} X']$ is well defined because the outer square in
$$
\begin{tikzcd}
V' \arrow[]{d} \arrow[]{r} & X' \arrow[]{d} \arrow[]{r} & Y' \arrow[]{d}  \\
V \arrow[]{r} & X \arrow[]{r} & Y
\end{tikzcd}
$$
is independent (by $\Cs$-independence), and because specialized morphisms are stable under independent pullbacks;

\item \emph{bivariant product}: the bivariant product is defined by bilinearly extending the formula
$$[V \to X] \bullet [W \to X] := [V' \to X],$$
where the morphism $V' \to X$ is as in the diagram
$$
\begin{tikzcd}
V' \arrow[]{r} \arrow[]{d} & X' \arrow[]{r} \arrow[]{d} & W \arrow[]{d} \\
V \arrow[]{r} & X \arrow[]{r} & Y \arrow[]{r} & Z;
\end{tikzcd}
$$
the element $[V' \to X]$ is well defined because $V' \to Z$ is the composition of the specialized maps $V' \to W$ and $W \to Z$.
\end{enumerate}
\end{defn}

Our first task is to verify that the above data satisfies the axioms of a bivariant theory.

\begin{lem}\label{lem:univthyisbiv}
$\Mb_\Fc$, equipped with the orientation defined by the formula
$$\theta(f) := [X \xto{\Id} X] \in \Mb_\Fc(X \xto{f} Y),$$ 
is a stably oriented commutative bivariant theory.
\end{lem}
\begin{proof}
We begin by verifying the bivariant axioms.
\begin{enumerate}
\item[$(A_1)$] Given elements $[V \to X] \in \Mb_\Fc(X \to Y)$, $[W \to Y] \in \Mb_\Fc(Y \to Z)$ and $[U \to Z] \in \Mb_\Fc(Z \to P)$, then both $([V \to X] \bullet [W \to Y]) \bullet [U \to Z]$ and $[V \to X] \bullet ([W \to Y] \bullet [U \to Z])$ coincide with $[V'_2 \to X] \in \Mb_\Fc(X \to P)$, where the map $V'_2 \to X$ is as in the three-dimensional Cartesian diagram
$$
\begin{tikzcd}[row sep=10, column sep=10]
 & V'_2 \arrow[dd] \arrow[rr] \arrow[ld] && X'_2 \arrow[dd] \arrow[ld] \arrow[rr] &&  W_2 \arrow[dd] \arrow[ld] \\
V' \arrow[crossing over]{rr} \arrow[]{dd} &&  X' \arrow[crossing over]{rr} && W  \\
 & V_2 \arrow[ld] \arrow[rr] && X_2 \arrow[ld] \arrow[rr] && Y_2 \arrow[ld] \arrow[rr] && U \arrow[ld]  \\
V \arrow[]{rr} && X \arrow[]{rr} \arrow[uu, <-, crossing over] & & Y \arrow[]{rr} \arrow[uu, <-, crossing over] && Z \arrow[rr] && P.
\end{tikzcd}
$$

\item[$(A_2)$] Follows directly from the definition.

\item[$(A_3)$] Follows directly from the definition.

\item[$(A_{12})$] If $[V \to X] \in \Mb_\Fc(X \to Y)$, $[W \to Y] \in \Mb_\Fc(X \to Y)$, and if $X \to Y$ factors through a confined map $f: X \to U$, then, by investigating the diagram
$$
\begin{tikzcd}
V' \arrow[]{r} \arrow[]{d} & X' \arrow[]{r} \arrow[]{d} & U' \arrow[]{r} \arrow[]{d} & W \arrow[]{d} \\
V \arrow[]{r} & X \arrow[]{r}{f} & U \arrow[]{r} & Y \arrow[]{r} & Z,
\end{tikzcd}
$$
we conclude that $f_*([V \to X] \bullet [W \to Y]) = f_*([V \to X]) \bullet [W \to Y]$;

\item[$(A_{13})$] If all the squares of the Cartesian diagram
$$
\begin{tikzcd}
X_2 \arrow[d] \arrow[r] & Y_2 \arrow[d]{}{g'} \arrow[r] & Z_2 \arrow[d]{}{g}  \\
X \arrow[r] & Y \arrow[r] & Z 
\end{tikzcd}
$$
are independent, and if $[V \to X] \in \Mb_\Fc(X \to Y)$ and $[W \to Y] \in \Mb_\Fc(Y \to Z)$, then both $g^*([V \to X] \bullet [W \to Y])$ and $g'^*([V \to X]) \bullet g^*([W \to Y])$ coincide with $[V'_2 \to X_2] \in \Mb_\Fc(X_2 \to Z_2)$, where $V'_2 \to X_2$ is as in the three-dimensional Cartesian diagram
$$
\begin{tikzcd}[row sep=10, column sep=10]
& V'_2 \arrow[dd] \arrow[rr] \arrow[ld] && X'_2 \arrow[dd] \arrow[rr] \arrow[ld] && W_2 \arrow[dd] \arrow[ld]\\
V' \arrow[rr, crossing over] && X' \arrow[rr, crossing over] && W  \\
& V_2 \arrow[rr] \arrow[ld] && X_2 \arrow[rr] \arrow[ld] && Y_2 \arrow[rr] \arrow[ld] && Z_2 \arrow[ld] \\
V \arrow[rr] \arrow[uu, <-] && X \arrow[rr] \arrow[uu, <-, crossing over] && Y \arrow[rr] \arrow[uu, <-, crossing over] && Z.
\end{tikzcd}
$$

\item[$(A_{23})$] If the large square and the rightmost little square of the Cartesian diagram
$$
\begin{tikzcd}
X' \arrow[d]{}{h''} \arrow[r]{}{f'} & Y' \arrow[d]{}{h'} \arrow[r] & Z' \arrow[d]{}{h}  \\
X \arrow[r]{}{f} & Y \arrow[r] & Z
\end{tikzcd}
$$
are independent, if $f$ is confined, and if $[V \to X] \in \Mb_\Fc(X \to Z)$, then both $h^*(f_*([V \to X]))$ and $f'_*(h^*([V \to X]))$ coincide with $[V' \to Y'] \in \Mb_\Fc(Y' \to Z')$, where $V' := V \times_Y Y'$.

\item[$(A_{123})$] Suppose
$$
\begin{tikzcd}
X' \arrow[d]{}{g'} \arrow[r] & Y' \arrow[d]{}{g} \\
X \arrow[r] & Y 
\end{tikzcd}
$$
is an independent square such that $g$ is confined. If $Y \to Z$ is a morphism, $[V \to Y']$ belongs to the bivariant group $\Mb_\Fc(Y' \to Z)$ of the composition, and $[W \to X] \in \Mb_\Fc(X \to Y)$, then both $g'_*(g^*([W \to X]) \bullet [V \to Y'])$ and $[W \to X] \bullet g_*([V \to Y'])$ coincide with $[V' \to X] \in \Mb_\Fc(X \to Z)$, where $V' := W \times_Y V$.

\item[$(U)$] The element $1_X := [X \xto{\Id} X] \in \Mb_\Fc(X \xto{\Id} X)$ satisfies the requirements of an unit.
\end{enumerate}
It is straightforward to check that $\theta$ is multiplicative, that $\theta(\Id) = 1$, and that $\theta$ is stable under pullbacks. The commutativity of $\Mb_\Fc$ follows from the commutativity of fiber products.
\end{proof}

Now that we have verified that $\Mb_\Fc$ is indeed a bivariant theory, we prove a universal property for it.

\begin{thm}\label{thm:univbivthyunivprop}
Let $\Fc$ be a bivariant functoriality and suppose that the class of independent squares is symmetric, that $\Fc$ satisfies $\Cs$-independence, and that the class of specialized morphisms is stable under independent pullbacks. Let $\Bb$ be a stably oriented bivariant theory with functoriality $\Fc$. Then, there exists a unique orientation preserving Grothendieck transformation $\eta: \Mb_\Fc \to \Bb$.
\end{thm}
\begin{proof}
Let us denote the orientation of $\Bb$ by $\theta'$. As
$$[V \xto{f} X] = f_*(\theta(g)) \in \Mb_\Fc(X \to Y),$$
where $g$ is the composition $V \to Y$, the desired Grothendieck transformation must satisfy
$$\eta([V \to X]) = f_*(\theta'(g)).$$
This formula gives rise to well-defined group homomorphisms
$$\eta_{X \to Y}: \Mb_\Fc(X \to Y) \to \Bb(X \to Y);$$
we need to verify that these homomorphisms are compatible with the bivariant structure.
\begin{enumerate}
\item \emph{The maps $\eta$ are compatible with pushforward}: this is obvious.

\item \emph{The maps $\eta$ are compatible with pullbacks}: consider a Cartesian diagram
$$
\begin{tikzcd}
V' \arrow[r]{}{f'} \arrow[d] & X' \arrow[r] \arrow[d] & Y'  \arrow[d]{}{h} \\
V \arrow[r]{}{f} & X \arrow[r] & Y
\end{tikzcd}
$$
in which the rightmost little square is independent, $f$ is confined, and the composition $V \to Y$ is specialized. From $\Cs$-independence it then follows that all squares of the diagram are independent. Let us denote the horizontal compositions by $g'$ and $g$ respectively. To verify the claim, we compute that 
\begin{align*}
\eta(h^*([V \to X])) &= \eta([V' \to X']) \\
&= f'_*(\theta'(g')) \\
&= h^*(f_*(\theta'(g))) & (\text{$A_{23}$, stability}) \\
&= h^*(\eta([V \to X])),
\end{align*}
as desired.

\item \emph{The maps $\eta$ are compatible with bivariant products}: consider a Cartesian diagram
$$
\begin{tikzcd}
V' \arrow[r]{}{f'} \arrow[d]{}{h''} & X' \arrow[r] \arrow[d]{}{h'} & W \arrow[d]{}{h} \\
V \arrow[r]{}{f} & X \arrow[r] & Y \arrow[r] & Z 
\end{tikzcd}
$$
in which $f$ and $h$ are confined, and the compositions $g: V \to Y$ and $g': W \to Z$ are specialized. By $\Cs$-independence, all the squares of the diagram are independent, and therefore the composition $g'': V' \to W$ is specialized. In order to verify the claim, we compute that
\begin{align*}
\eta([V \to X] \bullet [W \to Y]) &= \eta([V' \to X]) \\
&= h'_*\big(f'_*(\theta'(g' \circ g''))\big) \\
&= h'_*\big(f'_*(\theta'(g'') \bullet \theta'(g'))\big) \\
&= h'_*\big(f'_*(\theta'(g'')) \bullet \theta'(g')\big) & (A_{12}) \\
&= h'_*\big(h^*(f_*(\theta'(g))) \bullet \theta'(g')\big) & (\text{$A_{23}$, stability}) \\
&= f_*(\theta'(g)) \bullet h_*(\theta'(g')) & (A_{123}) \\
&= \eta([V \to X]) \bullet \eta([W \to Y]),
\end{align*}
as desired. \qedhere
\end{enumerate}
\end{proof}

\subsection{Universal additive bivariant theory}\label{ssect:univaddthy}

Here, we consider a straightforward variant of the the universal bivariant theory $\Mb_\Fc$, where the addition is compatible with disjoint unions.

\begin{defn}\label{def:univaddthy}
Suppose that the functoriality $\Fc$ has good coproducts, that it satisfies $\Cs$-independence, that the class of independent squares is symmetric, and that specialized morphisms are stable under independent pullbacks. 

Then, we define $\Ab_\Fc(X \to Y)$ as the free Abelian group on equivalence classes $[V \to X]$ of confined morphisms $V \to X$ such that the composition $V \to X$ is specialized, modulo all the \emph{decomposition relations}, i.e., the relations of form
$$[V_1 \sqcup V_2 \to X] = [V_1 \to X] + [V_2 \to X].$$
The bivariant operations are defined similarly to those of $\Mb_\Fc$ in Definition \ref{def:univbivthy}. These operations respect the decomposition relation because coproducts distribute over fiber products.
\end{defn}

\begin{lem}\label{lem:univaddthyisbiv}
$\Ab_\Fc$, equipped with the orientation defined by the formula
$$\theta(f) := [X \xto{\Id} X] \in \Ab_\Fc(X \xto{f} Y),$$ 
is an additive bivariant theory. Moreover, it is commutative and stably oriented.
\end{lem}
\begin{proof}
That $\Ab_\Fc$ is a commutative and stably oriented bivariant theory follows immediately from Lemma \ref{lem:univthyisbiv}. It is also clearly additive.
\end{proof}

\begin{thm}\label{thm:univaddthyunivprop}
Suppose that the functoriality $\Fc$ has good coproducts, that it satisfies $\Cs$-independence, that the class of independent squares is symmetric and that specialized morphisms are stable under independent pullbacks. Let $\Bb$ is a stably oriented additive bivariant theory with functoriality $\Fc$. Then, there exists a unique orientation preserving Grothendieck transformation $\eta: \Ab_\Fc \to \Bb$.
\end{thm}
\begin{proof}
Let us denote the orientation of $\Bb$ by $\theta'$. Following the proof of Theorem \ref{thm:univbivthyunivprop}, it suffices to show that the formula
$$\eta([V \xto{f} X]) := f_*(\theta'(g)),$$
where $g$ is the composition $V \to Y$, gives rise to well defined group homomorphism $\eta: \Ab_\Fc(X \to Y) \to \Bb(X \to Y)$ for each morphism $X \to Y$ in $\Cc$. In order for $\eta$ to be compatible with the decomposition relations, it suffices to show that, when $V \simeq V_1 \sqcup V_2$, then $\theta'(g) = \iota_{1*}(\theta'(g_1)) + \iota_{2*}(\theta'(g_2))$, where $\iota_i$ are the canonical inclusions and $g_i := g \circ \iota_i$. Indeed, if this was the case, then
\begin{align*}
\eta([V \xto{f} X]) &= f_*(\theta'(g)) \\
&= f_*(\iota_{1*}(\theta'(g_1))) + f_*(\iota_{2*}(\theta'(g_2))) \\
&= f_{1*}(\theta'(g_1)) + f_{2*}(\theta'(g_2)) \\
&= \eta([V_1 \xto{f_1} X]) + \eta([V_2 \xto{f_2} X]),
\end{align*}
where $f_i := f \circ \iota_i$, as desired. 

The claim follows from the fact that $1_V = \iota_{1*}(\theta(\iota_1)) + \iota_{2*}(\theta(\iota_2))$. Indeed, one computes that
\begin{align*}
\theta'(g) &= 1_V \bullet \theta'(g) \\ 
&= \big( \iota_{1*}(\theta(\iota_1)) + \iota_{2*}(\theta(\iota_2)) \big) \bullet \theta'(g) \\
&=  \iota_{1*}(\theta(\iota_1) \bullet \theta'(g)) + \iota_{2*}(\theta(\iota_2) \bullet \theta'(g)) & (A_{12}) \\
&= \iota_{1*}(\theta(g_1)) + \iota_{2*}(\theta(g_2)),
\end{align*}
which proves the claim.
\end{proof}

\section{Manipulation of bivariant theories}\label{sect:bivmanip}

In this section, we define some basic tools that have proven to be useful in the manipulation of bivariant theories. We start by studying bivariant ideals, the possible kernels of Grothendieck transformations, in Section \ref{ssect:bivid}. After this, in Section \ref{ssect:homcohid}, we investigate the kind of stability properties a bivariant ideal has after being restricted to the associated homology or cohomology theory. Finally, in Section \ref{ssect:extscal}, we recall how to extend the ring of coefficients of a bivariant ideal along a ring homomorphism.

\subsection{Bivariant ideals}\label{ssect:bivid}

In a sense, bivariant theories are like rings. Hence, in order to impose relations on bivariant theories, it is natural to consider ideals in them.

\begin{defn}\label{def:bivid}
Let $\Bb$ be a bivariant theory. A \emph{bivariant subset} $\Sc \subset \Bb$ is a collection of subsets
$$\Sc(X \to Y) \subset \Bb(X \to Y),$$
one for each morphism $X \to Y$ in $\Cc$. A bivariant subset $\Ic \subset \Bb$ is called a \emph{bivariant ideal} if
\begin{enumerate}
\item the subsets $\Ic(X \to Y) \subset \Bb(X \to Y)$ are subgroups;
\item the subsets $\Ic(X \to Y)$ are closed under bivariant pushforwards and pullbacks;
\item for any $s \in \Ic (X \to Y)$ and any $\alpha \in \Bb(V \to X)$ and $\beta \in \Bb(Y \to Z)$, we have that $\alpha \bullet s \in \Ic(V \to Y)$ and $s \bullet \beta \in \Ic(X \to Z)$.
\end{enumerate}
If $\Sc \subset \Bb$ is a bivariant subset, then the smallest bivariant ideal containing $\Sc$, the bivariant ideal \emph{generated by $\Sc$}, is denoted by $\langle \Sc \rangle_\Bb$. The subscript $\Bb$ is often omitted from the notation.$  $
\end{defn}

The following results are immediate.

\begin{prop}\label{prop:bividquot}
Let $\Bb$ be a bivariant theory. Then,
\begin{enumerate}
\item if $\Ic \subset \Bb$ is a bivariant ideal, then the quotient theory $\Bb/\Ic$, defined as
$$(\Bb / \Ic)(X \to Y) := \Bb(X \to Y) / \Ic(X \to Y),$$
and with the operations inherited from $\Bb$, is a bivariant theory;

\item if $\eta: \Bb \to \Bb'$ is a Grothendieck transformation, then the \emph{kernel}  $\Kc \subset \Bb$ of $\eta$, defined as
$$\Kc(X \to Y) := \ker(\eta_{X \to Y}),$$
is a bivariant ideal. \qed
\end{enumerate}
\end{prop}

\begin{prop}\label{prop:bivinherit}
Let $\Ic \subset \Bb$ be a bivariant ideal. Then
\begin{enumerate}
\item if $\Bb$ is an oriented bivariant theory with orientation $\theta$, then its image $\bar \theta$ provides an orientation for $\Bb/\Ic$;
\item the following properties of (oriented) bivariant theories are inherited by quotients:
\begin{enumerate}
\item commutativity;
\item being stably oriented;
\item being well oriented;
\item additivity;
\item $R$-linearity. \qed
\end{enumerate}
\end{enumerate}
\end{prop}

For certain kinds of bivariant theories, the generated bivariant ideal $\langle \Sc \rangle$ admits a concrete description.

\begin{prop}\label{prop:bividpres}
Let $\Bb$ be a bivariant theory with functoriality $\Fc$ satisfying $\Cs$-independence (Section \ref{ssect:univbivthy}). Let $\Sc \subset \Bb$ be a bivariant subset. Then $\langle \Sc \rangle(X \to Y)$ is the group generated by elements of form
$$f_*(\alpha \bullet g^*(s) \bullet \beta),$$
where $f$ and $g$ are as in the commutative diagram
$$
\begin{tikzcd}
& X \arrow[<-, ld, bend right]{}[swap]{f} \arrow[rrd, bend left] \\
A'' \arrow[r] & A' \arrow[r] \arrow[d] & B' \arrow[r] \arrow[d]{}{g} & Y \\
& A \arrow[r] & B,
\end{tikzcd}
$$
$f$ is confined, the bottom square is independent, $\alpha \in \Bb(A'' \to A')$, $s \in \Sc(A \to B)$, and $\beta \in \Bb(B' \to Y)$.
\end{prop}
\begin{proof}
Clearly such elements are contained in the generated ideal $\langle \Sc \rangle$, so it is enough to show that they form a bivariant ideal.
\begin{enumerate}
\item \emph{Stability under pushforwards}: this is obvious.

\item \emph{Stability under pullbacks}: let 
$$
\begin{tikzcd}
X' \arrow[r] \arrow[d]{}{h'} & Y' \arrow[d]{}{h} \\
X \arrow[r] & Y
\end{tikzcd}
$$
be an independent Cartesian square. Then, by $\Cs$-independence, so is
$$
\begin{tikzcd}
A''' \arrow[r]{}{f'} \arrow[d]{}{h''} & X' \arrow[d]{}{h'} \\
A'' \arrow[r]{}{f} & X,
\end{tikzcd}
$$
and we deduce from the bivariant axiom $(A_{23})$ that
$$h^*\big(f_*(\alpha \bullet g^*(s) \bullet \beta) \big) = f'_*\big(h^*(\alpha \bullet g^*(s) \bullet \beta) \big).$$
The right-hand-side element is of the correct form by bivariant axiom $(A_{13})$.

\item \emph{Stability under left multiplication}: let $\alpha \in \Bb(X' \to X)$ and consider the commutative diagram
$$
\begin{tikzcd}
A''' \arrow[r] \arrow[d]{}{f'} & A'' \arrow[d]{}{f} \arrow[rd] \\
X' \arrow[r] & X \arrow[r] & Y.
\end{tikzcd}
$$
As $f$ is confined, the square is independent. Using bivariant axiom $(A_{123})$, we compute that
$$\alpha' \bullet f_*(\alpha \bullet g^*(s) \bullet \beta) = f'_*\big((f^*(\alpha') \bullet \alpha) \bullet g^*(s) \bullet \beta \big).$$
The right-hand-side element is of the correct form.

\item \emph{Stability under right multiplication}: let $\beta' \in \Bb(Y \to Y')$. Then, by bivariant axiom $(A_{12})$,
$$f_*(\alpha \bullet g^*(s) \bullet \beta) \bullet \beta' = f_*(\alpha \bullet g^*(s) \bullet (\beta \bullet \beta')).$$
The right-hand-side element is of the correct form. \qedhere
\end{enumerate}
\end{proof}

\subsection{Homological and cohomological ideals}\label{ssect:homcohid}

Relating the universal properties of a bivariant theory and of the associated homology and cohomology theories is often done by relating quotients of bivariant theories to quotients of homology and cohomology theories. To do so efficiently, we consider the following structures.

\begin{defn}\label{def:homid}
Let $\Bb$ be a bivariant theory. Then a \emph{homological subset} $\Sc \subset \Bb$ is a bivariant subset such that, if $Y$ is not a final object, then $\Sc(X \to Y) = \emptyset$. A homological subset $\Ic$ is a \emph{homological ideal} if, for all $X \in \Cc$, $\langle \Ic \rangle(X \to pt) = \Ic (X \to pt)$. The smallest homological ideal containing the homological subset $\Sc$, the \emph{homological ideal generated} by $\Sc$, is denoted by $\langle  \Sc \rangle_{\Bb_\bullet}$, or by $\langle  \Sc \rangle_\bullet$ if the bivariant theory $\Bb$ is clear from the context.
\end{defn} 

\begin{defn}\label{def:cohid}
Let $\Bb$ be a bivariant theory. Then a \emph{cohomological subset} $\Sc \subset \Bb$ is a bivariant subset such that, if $X \to Y$ is not the identity morphism, then $\Sc(X \to Y) = \emptyset$. A cohomological subset $\Ic$ is a \emph{cohomological ideal} if, for all $X \in \Cc$, $\langle \Ic \rangle(X \to X) = \Ic (X \to X)$. The smallest cohomological ideal containing the cohomological subset $\Sc$, the \emph{cohomological ideal generated} by $\Sc$, is denoted by $\langle  \Sc \rangle_{\Bb^\bullet}$, or by $\langle  \Sc \rangle^\bullet$ if the bivariant theory $\Bb$ is clear from the context.
\end{defn} 

It turns out that, under certain assumptions, homological and cohomological ideals admit a concrete characterization. The following notion will be useful for the characterization.

\begin{defn}\label{def:defgenor}
Let $\Bb$ be an oriented bivariant theory. Then $\Bb$ is \emph{generated by orientations} if, for each map $X \to Y$ in $\Cc$, $\Bb(X \to Y)$
is generated, as a $\Bb^\bullet(pt) \dash \Bb^\bullet(pt)$-bimodule, by elements of form
$$f_*(\theta(g)),$$
where $f: V \to X$ is a confined morphism such that the composition $g: V \to Y$ is specialized.
\end{defn}

We next characterize homological ideals. Note the similarity between this characterization and that of the sets of relations used in e.g. \cite{levine-morel, lowrey--schurg}.

\begin{prop}\label{prop:homidchar}
Let $\Bb$ be a bivariant theory with functoriality $\Fc$ that satisfies $\Cs$-independence, and in which all absolute product squares are independent. Then a homological subset $\Ic$ is a homological ideal if and only if
\begin{enumerate}
\item for all $X \in \Cc$, $\Ic(X \to pt) \subset \Bb(X \to pt)$ is a subgroup;
\item $\Ic$ is closed under pushforwards in the associated homology theory;
\item given $s \in \Ic(X \to pt)$ and $\alpha \in \Bb(Y \to pt)$, we have that $s \times \alpha \in \Bb(X \times Y \to pt)$ and $\alpha \times s \in \Bb(Y \times X \to pt)$;
\item for all $\alpha \in \Bb(X \to Y)$ and $s \in \Ic(Y \to pt)$, $\alpha \bullet s \in \Ic(X \to pt)$.
\end{enumerate}
If $\Bb$ is oriented and generated by orientations, then the condition 4. is equivalent to
\begin{enumerate}
\item[4$'.$] $\Ic$ is closed under Gysin pullbacks in the associated homology theory.
\end{enumerate}
\end{prop}
\begin{proof}
Since absolute product squares are independent, homological cross products are well defined (Section \ref{ssect:bivcrossprod}). By Proposition \ref{prop:bividpres}, $\Ic(X \to pt)$ is generated by elements of form
$$f_*(\alpha \bullet \pi_{Y'}^*(s) \bullet \beta),$$
where the notation is as in the commutative diagram
$$
\begin{tikzcd}
& X \arrow[<-, ld, bend right]{}[swap]{f} \arrow[rrd, bend left] \\
V \arrow[r] & A \times B \arrow[r] \arrow[d] & B \arrow[r] \arrow[d]{}{\pi_{Y'}} & pt \\
& A \arrow[r] & pt,
\end{tikzcd}
$$
$f$ is confined, $\alpha \in \Bb(V \to A \times B)$, $s \in \Bb(A \to pt)$ and $\beta \in \Bb(B \to pt)$. As 
\begin{align*}
s' :&= \pi^*_{Y'}(s) \bullet \beta \\
&= s \times \beta \in \Ic(A \times B \to pt),
\end{align*}
such elements are contained in $\Ic$ satisfying the assumptions 1.-4. above. On the other hand, a homological ideal $\Ic$ always satisfies the assumptions 1.-4.: to see that $\alpha \times s \in \Ic(Y \times X \to pt)$, we note that $\alpha \times s = \pi^*_X(\alpha) \bullet s$, and such an element is contained in $\Ic(Y \times X \to pt)$ by stability under left multiplication.

Suppose then that $\Bb$ is an oriented bivariant theory generated by its orientations. It is clear that 4. implies 4$'.$  To see that the converse holds, it is enough to prove that $\Ic$ is stable under left multiplications by elements of form $af_*(\theta(g))b$, where $a, b \in \Bb^\bullet(pt)$. As the left $\Bb^\bullet(pt)$-module structure on $\Bb^*(X)$ coincides with that given by the cross product, it suffices to show that $\Ic$ is stable under left multiplications by elements of form $f_*(\theta(g))$. But this is obvious, as
\begin{align*}
f_*(\theta(g)) \bullet s &= f_*(\theta(g) \bullet s) \\
&=  f_*(g^!(s)),
\end{align*}
and the last element belongs to $\Ic$ by $4'.$ and $2.$
\end{proof}

There exists the following analogous characterization for cohomological ideals.

\begin{prop}
Let $\Bb$ be a stably oriented bivariant theory with functoriality $\Fc$ that satisfies $\Cs$-independence. Suppose that $\Bb$ is generated by orientations. Then, a cohomological subset $\Ic$ is a cohomological ideal if and only if
\begin{enumerate}
\item for all $X \in \Cc$, $\Ic(X \to X) \subset \Bb(X \to X)$ is an ideal;

\item $\Ic$ is closed under pullbacks in the associated cohomology theory;

\item $\Ic$ is closed under Gysin pushforwards in the associated cohomology theory.
\end{enumerate}
\end{prop}
\begin{proof}
By Proposition \ref{prop:bividpres}, $\langle \Ic \rangle (Y \to Y)$ is generated by elements of form $f_*(\alpha \bullet s \bullet \beta)$, where $f$ is a confined morphism that factors as
$$V \to X \to Y,$$
$\alpha \in \Bb(V \to X)$, $s \in \Ic(X \to X)$ and $\beta \in \Bb(X \to Y)$. As the sets $\Ic(X \to X)$ are ideals, and as $\Bb$ is generated by orientations, we may assume that $\alpha = g_*(\theta(h))$, where $g$ is a confined morphism $W \to V$ such that the composition $h: W \to X$ is confined. 

Our first claim is that for all $r \in \Bb(X \to X)$
$$g_*(\theta(h)) \bullet r = g_*\big(h^*(r) \bullet \theta(h)\big).$$
It is enough to verify that this formula holds for $r$ of form $l_*(\theta(l))$ where $l: U \to X$ is a specialized and confined morphism. Since both $g_*(\theta(h))$ and $l_*(\theta(l))$ are in the image of the map $\Mb_\Fc \to \Bb$ from the universal theory $\Mb_\Fc$ (Section \ref{ssect:univbivthy}), it is enough to verify that this formula holds in $\Mb_\Fc$. Forming the Cartesian diagram
$$
\begin{tikzcd}
W' \arrow[r]{}{g'} \arrow[d]{}{l''} & V' \arrow[r] \arrow[d]{}{l'} & U \arrow[d]{}{l} \arrow[rd] \\
W \arrow[r]{}{g} & V \arrow[r] & X \arrow[r]{}{\Id} & X,
\end{tikzcd}
$$
we observe that the product equals $g_*\big(l''_*(\theta(l'')) \bullet \theta(h) \big)$. As 
\begin{align*}
l''_*(\theta(l'')) = h^*(l_*(\theta(l))),
\end{align*}
this proves the claim.

Thus, applying bivariant axiom $(A_{12})$, we conclude that $\langle \Ic \rangle (Y \to Y)$ is generated by elements of form $f_*(s \bullet \beta)$, where $f: X \to Y$ is confined, $\beta \in \Bb(X \to Y)$, and $s \in \Ic(X \to X)$. Without loss of generality, $\beta = g_*(\theta(h))$, where $g: V \to X$ is a confined morphism such that the composition $h: V \to Y$ is specialized. Forming the commutative diagram
$$
\begin{tikzcd}
V \arrow[d]{}{g} \arrow[r] & V \arrow[d]{}{g} \arrow[rd]{}{h} \\
X \arrow[r] & X \arrow[r]{}{f} & Y,
\end{tikzcd}
$$
we compute that 
\begin{align*}
f_*(s \bullet g_*(\theta(h)) ) &= f_*\big(g_*(g^*(s) \bullet \theta(h))\big) & (A_{123}) \\
&= h_*(g^*(s) \bullet \theta(h)) \\
&= h_!(g^*(s)),
\end{align*}
which is contained in $\Ic$ by the assumptions made above. Hence, $\Ic$ satisfying the above conditions is a cohomological ideal.

On the other hand, it is clear that a cohomological ideal satisfies the conditions 1.-3. above.
\end{proof}

\subsection{Change of coefficients}\label{ssect:extscal}

We have considered $R$-linear bivariant theories in Section \ref{ssect:bivprop}. Given an $R$-linear theory $\Bb$, and a homomorphism $R \to R'$ of rings, it will be useful to construct, in a universal fashion, an $R'$-linear bivariant theory $\Bb'$. For example, universal bivariant theories with the desired formal group law will be constructed in this fashion later, by using the case where $R = \Lb$ and where the homomorphism $\Lb \to R'$ represents the desired formal group law. 

\begin{defn}\label{def:extscal}
Let $R$ be a commutative ring, and let $\Bb$ be a commutative $R$-linear bivariant theory. Let $R'$ be an $R$-algebra. Then the bivariant theory $R' \otimes_R \Bb$ is defined by the formula
$$(R' \otimes_R \Bb)(X \to Y) := R' \otimes_R \Bb(X \to Y)$$
and using the bivariant operations inherited from $\Bb$. It is clear that there exists a canonical $R$-linear Grothendieck transformation $\Bb \to R' \otimes_R \Bb$.
\end{defn}

Note that if $R \to R'$ is a ring homomorphism, then any $R'$-linear bivariant theory $\Bb'$ may be regarded as an $R$-linear theory. Change of coefficients is the left adjoint to this process.

\begin{prop}\label{prop:extscal}
Let $\Bb$ and $\Bb'$ be $R$ and $R'$-linear bivariant theories respectively, let $R \to R'$ be a ring homomorphism, and let $\eta: \Bb \to \Bb'$ be an $R$-linear Grothendieck transformation. Then there exists a unique $R'$-linear Grothendieck transformation
$$\eta': R' \otimes_R \Bb \to \Bb'$$
such that the composition 
$$\Bb \to R' \otimes_R \Bb \xto{\eta'} \Bb'$$
coincides with $\eta$.
\end{prop}
\begin{proof}
Follows immediately from the corresponding property of modules.
\end{proof}

Many properties of bivariant theories are conserved in extension of scalars.

\begin{prop}\label{prop:bivinherit2}
Let $\Bb$ be an $R$-linear bivariant theory and let $R'$ be an $R$-algebra. Then
\begin{enumerate}
\item if $\Bb$ is an oriented bivariant theory with orientation $\theta$, then its image $\theta'$ provides an orientation for $R' \otimes_R \Bb$;
\item the following properties of (oriented) bivariant theories are preserved in change of coefficients:
\begin{enumerate}
\item being stably oriented;
\item being well oriented;
\item additivity. \qed
\end{enumerate}
\end{enumerate}
\end{prop}

\chapter{Derived algebraic geometry}\label{ch:dac}

In this preliminary chapter, we will introduce the reader to the basics of derived algebraic geometry. The fundamental references of the subject are the joint work of Toën and Vezzosi \cite{HAG1, HAG2}, as well as the work of Jacob Lurie \cite{HTT, HA, SAG}. Besides these, other useful references include the book of Gaitsgory and Rozenblym \cite{GR} and the survey article by Toën \cite{toen:2014}. 

In derived algebraic geometry, one studies derived schemes, which are geometric objects that are locally modeled by the spectra of derived rings, which, roughly speaking, are the (algebraic-structure-preserving-)homotopy types of (nice enough) topological rings. The topology on these rings allows one to efficiently, and in a natural fashion, keep track of extra data, that is related to degeneracy of the equations defining the algebro-geometric object. This turns out to be rather useful when studying cohomology theories in algebraic geometry that are related to algebraic cobordism, as we will observe later in this thesis.

The structure of this chapter is as follows. In Section \ref{sect:dacbasics} we will introduce the basics of derived algebraic geometry. Due to the technical nature of the subject, our discussion will be woefully far from being self contained. In Section \ref{sect:qproj}, we will introduce the theory of ample line bundles and quasi-projective morphisms in the context of derived geometry. Section \ref{sect:dci} deals with derived complete intersections, a subject of critical importance for our construction of algebraic cobordism. Finally, in Section \ref{sect:blowpup}, we recall the derived blowup construction of Khan and Rydh, which will be an important tool for performing cobordism computations.

\section{Basics of derived geometry}\label{sect:dacbasics}

Here, we state the basic facts of derived algebraic geometry. The results of this section are well known to the experts.

\subsection{Derived rings and derived schemes}

In this section, we introduce derived rings and derived schemes. We start from the following definition from Chapter 25 of \cite{SAG}.

\begin{defn}\label{def:dring}
Let $\Poly$ be the category whose objects are the polynomial rings $\Zb[x_1,...,x_r]$, and whose morphisms are ring homomorphisms. The $\infty$-category $\dRing$ of \emph{derived rings} is defined as finite-product preserving\footnote{Products in $\Poly^\op$ are given by tensor products of commutative $\Zb$-algebras.} functors $\Poly^\op \to \Sc$, where $\Sc$ is the $\infty$-category of spaces.
\end{defn}

Let us unpack the above abstract definition. If $A$ is a derived ring, then $|A| := A(\Zb[x])$ is to be interpreted as its underlying space. Moreover, the maps $* \to |A|$, induced by the maps $\Zb[x] \to \Zb$ defined as $x \mapsto 0$ and $x \mapsto 1$, give the zero and the identity point of the derived ring. The addition and multiplication is encoded by the maps $|A| \times |A|  \to |A|$ induced by the maps $\Zb[x] \to \Zb[y, z]$ defined as $x \mapsto y + z$ and $x \mapsto yz$, respectively. Rest of the structure of a functor encodes higher coherence data for these structures.

\begin{rem}\label{rem:dringscr}
The category $\dRing$ is equivalent to the $\infty$-category obtained from simplicial commutative rings with its usual model structure \cite{quillen:1967}. The above definition is just a concise way to write this down directly in $\infty$-categorical terms.
\end{rem}

\begin{defn}\label{def:dringtrunc}
Let $A$ be a derived ring. The \emph{truncation} $\pi_0(A)$ of $A$ is the commutative ring, the elements of which are the connected components of $A$. Moreover, for each $i \geq 0$, the homotopy groups $\pi_i(|A|; 0)$ admit the structure of a $\pi_0(A)$-module. These modules are called the \emph{homotopy modules} of $A$, and they are denoted by $\pi_i(A)$.
\end{defn}

Next, we define some useful classes of derived rings and their morphisms.

\begin{defn}\label{def:dringprops}
Let $A$ be a derived ring. Then,
\begin{enumerate}
\item $A$ is \emph{Noetherian} if the truncation $\pi_0(A)$ is Noetherian, and if the homotopy modules $\pi_i(A)$ are finitely generated $\pi_0(A)$-modules; the \emph{Krull dimension} of a Noetherian derived ring is the Krull dimension of its truncation;

\item $A$ is \emph{local} if the truncation $\pi_0(A)$ is local; if $A'$ is another local derived ring, then a ring homomorphism $\psi: A \to A'$ is local if and only if the truncation $\pi_0(\psi): \pi_0(A) \to \pi_0(A)$ is local in the usual sense;

\item if $A'$ is a derived ring and $\psi: A \to A'$ is a homomorphism of derived rings, then $\psi$ is \emph{of finite type} (a \emph{surjection}) if the truncation $\pi_0(\psi): \pi_0(A) \to \pi_0(A)$ is a morphism of finite type (surjective);

\item if $A'$ is a derived ring and $\psi: A \to A'$ is a homomorphism of derived rings, then,
\begin{enumerate}
\item $\psi$ is \emph{finitely presented} if it can be obtained as a finite colimit of copies of $A[x]$; 

\item $\psi$ is \emph{locally of finite presentation} if it is a retract of a finitely presented derived $A$-algebra;

\item $\psi$ is \emph{almost of finite presentation}, if for all $n \gg 0$, there exists a finitely presented derived $A$-algebra $B$, and a derived $A$-algebra homomorphism $B \to A'$ such that $\pi_i(B) \to \pi_i(A)$ is an isomorphism for $i \leq n$.
\end{enumerate}
\end{enumerate}
\end{defn}

Concretely, finitely presented derived $A$-algebras are those that can be formed by adjoining finitely many ``free'' variables (of varying simplicial degrees) to $A$. Similarly, almost finitely presented derived $A$-algebras are those that can be formed by adjoining ``free'' variables to $A$, with only finitely many variables for each simplicial degree.

We record the following useful observation.

\begin{prop}\label{prop:noethfintyp}
A finite type morphism $\psi: A \to A'$ between Noetherian derived rings is almost of finite presentation.
\end{prop}
\begin{proof}
See \cite{lurie-thesis} Proposition 3.1.5.
\end{proof}

Next, we globalize the notion of derived rings in order to obtain derived schemes. However, before being able to do so, we introduce the following useful concept.

\begin{defn}[cf. \cite{SAG} Construction 1.1.2.2]\label{def:dringspc}
A \emph{derived ringed space} $X$ is a pair $(|X|, \Oc_X)$ consisting of a topological space $|X|$ and a sheaf $\Oc_X$ of derived rings on $|X|$. Given a sheaf of rings $\Rc$ on a topological space $X$, and a continuous map $f: X \to Y$, one defines the \emph{pushforward sheaf} $f_* \Rc$ by the formula $(f_* \Rc)(V):= \Rc(f^{-1}V)$. One defines the \emph{$\infty$-category of derived ringed spaces}, $\mathrm{Top}_\dRing$, as the source of the coCartesian fibration $\Sc_\dRing \to \mathrm{Top}$ that represents the functor assigning to each topological space $X \in \mathrm{Top}$ and to each morphism, the pushforward functor $f_*$. In other words, a (1-)morphism $f: X \to Y$ of derived ringed spaces consists of a continuous map $|f|: |X| \to |Y|$ and a map $f^\sharp: \Oc_Y \to |f|_* \Oc_X$ of sheaves of derived rings.

A derived ringed space $X$ is a \emph{derived locally ringed space} if its stalks $\Oc_{X,x}$ (defined in the obvious fashion) are local derived rings. The \emph{$\infty$-category of derived locally ringed spaces}, $\mathrm{Top}_\dRing^\mathrm{loc}$, is the subcategory of derived ringed spaces consisting of derived locally ringed spaces with those morphisms $f: X \to Y$ such that the induced morphisms $\Oc_{Y,|f|(x)} \to \Oc_{X,x}$ are local.
\end{defn}

\begin{defn}\label{def:dsch}
A \emph{derived scheme} is a derived locally ringed space $X = (|X|, \Oc_X)$, where $\Oc_X$ is a hypercomplete\footnote{Any sheaf on a finite-Krull-dimensional Noetherian topological space is automatically hypercomplete \cite{HTT}. Since we will almost exclusively focus on schemes whose underlying topological space satisfies this property, the reader may safely ignore this concept altogether.} sheaf of derived rings on $X$, satisfying
\begin{enumerate}
\item the \emph{truncation} $X_\cl := (|X|, \pi_0(\Oc_X))$ is a scheme in the classical sense;
\item each of the \emph{homotopy sheaves} $\pi_i(\Oc_X)$, defined as the sheafifications of the presheaves that assign the $\pi_0(\Oc_X(U))$-module $\pi_i(\Oc_X(U))$ to an open subset $U \subset |X|$, is quasi-coherent as a $\pi_0(\Oc_X)$-module.
\end{enumerate}
The \emph{$\infty$-category of derived schemes}, $\dSch$, is the full subcategory of derived locally ringed spaces consisting of derived schemes. 
\end{defn}

As in Chapter 1.1 of \cite{SAG}, one shows that for each derived ring $A$, one can construct a derived locally ringed space $\Spec(A)$, the \emph{spectrum of $A$}, the underlying topological space of which is the Zariski spectrum of $\pi_0(A)$, and satisfying the universal property that for each locally ringed space $X$,
$$\Map_{\mathrm{Top}_\dRing^\mathrm{loc}}(X, \Spec(A)) \simeq \Map_{\dRing}(A, \Oc_X(X)).$$
It is then possible to show that a derived locally ringed space $X$ is a derived scheme if and only it admits an open cover by spectra of derived rings. 

\begin{rem}[Functor of points of a derived scheme]\label{rem:dschfpts}
Another useful perspective to derived schemes is to consider them via their functor of points. A derived scheme $X$ defines a functor $h_X: \dRing \to \Sc$, \emph{the functor of points of $X$}, by the formula
$$h_X(A) := \Map_{\dSch}(\Spec(A), X).$$
It turns out that this is a fully faithful embedding of $\dSch$ to $\Fun(\dRing, \Sc)$ (see \cite{SAG} Chapter 1.6 for the analogous claim in spectral algebraic geometry), and the essential image of this embedding consists of functors that admit a Zariski covering\footnote{Defined in an analogous way to the Zariski coverings of set-valued functors on commutative rings.} by functors of form $\Map_\dRing(A, -)$, where $A$ is a derived ring\footnote{Note that this is just the functor of points of $\Spec(A)$.}.
\end{rem}

The $\infty$-category of derived schemes has analogous properties to the category of ordinary schemes. Importantly, it admits fiber products.

\begin{defn}\label{def:dfprod}
Let $X \to Z$ and $Y \to Z$ be maps of derived schemes. Then, the \emph{fiber product} of $X$ and $Y$ over $Z$ is a derived scheme $X \times_Z Y$, together with maps $X \times_Z Y \to X$ and $X \times_Z Y \to Y$, such that for every derived scheme $V$, the induced map
$$\Map_{\dSch}(V, X \times_Z Y) \to \Map_{\dSch}(V, X) \times_{\Map_{\dSch}(V, Z)} \Map_{\dSch}(V, Y)$$
is an equivalence of spaces. In other words, the datum of a map $V \to X \times_Z Y$ is equivalent to the data of maps $V \to X$, $V \to Y$ and a homotopy of the two compositions $V \to Z$ in $\Map_\dSch(V, Z)$.
\end{defn}

It is possible to check locally that the above functor of points description is locally modeled by tensor products of derived rings\footnote{Concretely, derived relative tensor products of simplicial commutative rings.}, and therefore is indeed a derived scheme. There is an important subtlety related to fiber products in derived algebraic geometry. Namely, even though the category of classical schemes embeds fully faithfully to $\dSch$ as the derived schemes with a discrete structure sheaf, this subcategory is not closed under fiber products.

\begin{prop}\label{prop:schfprod}
Let $X \to Z$ and $Y \to Z$ be classical schemes, and let $X \times^\cl_Z Y$ denote the classical fiber product. Then the canonical map $X \times^\cl_Z Y \to X \times_Z Y$ is an equivalence if and only if the fiber product is Tor-independent. \qed
\end{prop}

Let us then list some special classes of derived schemes and their morphisms.

\begin{defn}\label{def:dschnoeth}
A derived scheme $X$ is \emph{Noetherian} if its underlying topological space is quasi-compact and if it admits an affine open cover by spectra of Noetherian derived rings. The \emph{Krull dimension} of a Noetherian derived scheme $X$ is the Krull-dimension of its underlying topological space $|X|$.
\end{defn}

\begin{defn}\label{def:dschqcqs}
A derived scheme $X$ is \emph{qcqs} (quasi-compact and quasi-separated) if the truncation $X_\cl$ is.
\end{defn}

\begin{defn}\label{def:dschmaps}
Let $f: X \to Y$ be a map of derived schemes. Then,
\begin{enumerate}
\item $f$ is an \emph{open embedding} if $f$ identifies $X$ as $(U, \Oc_Y \vert U)$, where $U \subset |Y|$ is an open subset and $\Oc_Y \vert U$ is the restriction of the structure sheaf of $Y$ to $U$;
\item $f$ is a \emph{proper morphism} (a \emph{finite type morphism}, a \emph{qcqs morphism}, a \emph{closed embedding}) if the truncation $f_\cl: X_\cl \to Y_\cl$ is in the usual sense;
\item $f$ is a \emph{locally closed embedding} if it factors as a closed embedding followed by an open embedding.

\item $f$ is \emph{locally of finite presentation} (\emph{locally of almost finite presentation}) if, for every affine open $\Spec(A) \subset Y$ and an affine open $\Spec(B) \subset X$ mapping to $\Spec(A)$, the induced map $A \to B$ is locally of finite presentation (of almost finite presentation) as a map of derived rings.
\end{enumerate}
\end{defn}

Clearly a map $f: X \to Y$ is of finite type (a closed embedding) if and only if, for every affine open $\Spec(A) \subset Y$ and every affine open $\Spec(B) \subset X$ mapping to $\Spec(A)$, the induced map of derived rings $A \to B$ is of finite type (surjective).

\subsection{Quasi-coherent sheaves}\label{ssect:qcoh}

To every derived scheme $X$, one can associate the symmetric monoidal stable $\infty$-category $\QCoh(X)$ of \emph{quasi-coherent sheaves} on $X$. There are many equivalent characterizations of this category, but perhaps the most intuitive one is as the full sub-$\infty$-category of derived $\Oc_X$-modules\footnote{Concretely, a \emph{derived module} $M$ over a derived ring $A$ is described by a spectrum object in simplicial modules over a simplicial ring $A_\bullet$ representing $A$ in the Quillen model structure. If $M$ is connective, then it can be described by a single simplicial $A_\bullet$-module.} spanned by those modules that, on each affine open $\Spec(A) \subset \Oc_X$, are modeled by a derived $A$-module. If $X$ is a classical scheme, then $\QCoh(X)$ is closely related to the unbounded derived category of $X$ (the derived category can be recovered by its homotopy category with the induced triangulated structure, see \cite{HA} Chapter 1). For the rigorous definition, we refer to Chapter 2 of \cite{SAG}\footnote{As every derived scheme has an underlying spectral scheme, and as the (symmetric monoidal stable) $\infty$-category of quasi-coherent sheaves on a derived scheme is identified with that of the underlying spectral scheme, we may refer to \cite{SAG} for results concerning quasi-coherent sheaves in derived algebraic geometry.}.

We will denote suspensions of objects $\Fc \in \QCoh(X)$ by $\Fc[1]$. Moreover, each quasi-coherent sheaf $\Fc$ has naturally associated \emph{homotopy sheaves} $\pi_i(\Fc)$ for all $i \in \Zb$, which are naturally regarded as classical quasi-coherent sheaves on the truncation $X_\cl$. Additionally, each $\Fc \in \QCoh(X)$ takes part in a cofiber sequence
$$\tau_{\geq i + 1} \Fc \to \Fc \to \tau_{\leq i} \Fc,$$
where $\tau_{\geq n} \Fc \to \Fc$ is the universal map to $\Fc$ from a quasi-coherent sheaf with $\pi_j$ trivial for $j<n$ and where $\Fc \to \tau_{\leq n} \Fc$ is the universal map from $\Fc$ to a quasi-coherent sheaf with $\pi_j$ trivial for $j>n$. Moreover, $\pi_{j}(\tau_{\geq n} \Fc) = \pi_{j}(\Fc)$ for $j \geq n$ and  $\pi_{j}(\tau_{\leq n} \Fc) = \pi_{j}(\Fc)$ for $j \leq n$.

\begin{defn}\label{def:connective}
A quasi-coherent sheaf $\Fc$ is \emph{connective} if the natural map $\tau_{\geq 0} \Fc \to \Fc$ is an equivalence. If is \emph{eventually connective} if there exists an $i \in \Zb$ such that the natural map $\tau_{\geq i} \Fc \to \Fc$ is an equivalence.
\end{defn}

\begin{defn}\label{def:ggen}
A map $\Fc \to \Gc$ of connective quasi-coherent sheaves on $X$ is \emph{surjective} if the induced map $\pi_0(\Fc) \to \pi_0(\Gc)$ is a surjection of sheaves in the classical sense. A connective quasi-coherent sheaf $\Fc$ is \emph{globally generated} if it admits a surjection from $\bigoplus_{i \in I} \Oc_X$, where $I$ is an indexing set.
\end{defn}

Much like in classical algebraic geometry, given a qcqs morphism $f: X \to Y$ of derived schemes, there exists an adjoint pair of functors consisting of the \emph{pushforward}
$$f_*: \QCoh(X) \to \QCoh(Y)$$
and the \emph{pullback}\footnote{Pullbacks are defined along arbitrary morphisms of derived schemes.}
$$f^*: \QCoh(Y) \to \QCoh(X).$$
If the morphism $f$ is clear from the context, $f^* \Fc$ is often denoted by $\Fc \vert_X$. These operations are functorial and satisfy many good properties, as exemplified by the following result (\cite{SAG} Proposition 2.5.4.5 and Remark 3.4.2.6).

\begin{prop}\label{prop:qcohpushpullproj}
Pushforwards and pullbacks of quasi-coherent sheaves satisfy the following properties:
\begin{enumerate}
\item \emph{push-pull formula\footnote{Note that, in classical algebraic geometry, the analogous result (with derived pushforwards and pullbacks) holds only in Tor-independent Cartesian squares.}:} given a Cartesian square
$$
\begin{tikzcd}
X' \arrow[d]{}{g'} \arrow[r]{}{f'} & Y' \arrow[d]{}{g} \\
X \arrow[r]{}{f} & Y
\end{tikzcd}
$$
such that $g$ is qcqs, then there exists a natural equivalence of functors $f^* \circ g_* \simeq g'_* \circ f'^*: \QCoh(Y') \to \QCoh(X)$;
\item \emph{projection formula:} given a qcqs morphism $f: X \to Y$, $\Fc \in \QCoh(X)$, and $G \in \QCoh(Y)$, then there exists a natural equivalence
$$f_*(f^*(\Gc) \otimes \Fc) \simeq \Gc \otimes f_*(\Fc).$$ 
\end{enumerate}
\end{prop}

Another fundamental concept in the study of quasi-coherent sheaves is that of the global sections.

\begin{defn}\label{def:gsect}
Let $\Fc$ be a quasi-coherent sheaf on a qcqs derived scheme $X$. Then the \emph{spectrum of global sections} $\Gamma(X; \Fc)$ is defined as the pushforward of $\Fc$ to $\Spec(\Zb)$ with its natural spectrum structure. The \emph{space of global sections} $|\Gamma(X;\Fc)|$ is defined as the mapping space $\Map_{\QCoh(X)}(\Oc_X, \Fc)$.
\end{defn}

As $\QCoh(X)$ is a stable $\infty$-category, it has a canonical enrichment over the $\infty$-category of spectra. Moreover, it is not difficult to show that $\Gamma(X; \Fc)$ is recovered as the spectrum of maps $\Oc_X \to \Fc$. Hence, $|\Gamma(X;\Fc)|$ is the underlying space of $\Gamma(X;\Fc)$.

The following lemma is often used to reduce to a situation where only finitely many of the homotopy sheaves $\pi_i(\Fc)$ are nontrivial.

\begin{lem}[Grothendieck vanishing]\label{lem:gvanish}
Let $X$ be a Noetherian derived scheme with Krull dimension $n < \infty$, and let $\Fc \in \QCoh(X)$. Then, for all $j \in \Zb$, 
$$\pi_j\big(\Gamma(X; \Fc)\big) = \pi_j\big(\Gamma(X; \tau_{\leq i}(\Fc))\big)$$
as soon as $i \geq n+j$.
\end{lem}
\begin{proof}
Using the fundamental cofibre sequences
$$\Gamma\big(X; \pi_{i+1}(\Fc)[i+1] \big) \to \Gamma\big(X; \tau_{\leq i+1}(\Fc) \big) \to \Gamma\big(X; \tau_{\leq i}(\Fc) \big)$$
and the classical Grothendieck vanishing, we conclude that $\pi_j\big( \Gamma(X; \tau_{\leq i}(\Fc)) \big)$
stabilize for $i \gg 0$, and therefore (by Milnor sequence)
$$\pi_j \big( \Gamma(X; \Fc) \big) = \lim_i \pi_j\big(\Gamma(X; \tau_{\leq i}(\Fc))\big).$$
The bound for $i$ follows from the fact that a Noetherian topological space of Krull dimension $n$ has cohomological dimension $n$.
\end{proof}

Next, we list the most important finiteness conditions of quasi-coherent sheaves and state their basic properties.

\begin{defn}\label{def:tordim}
Let $\Fc$ be a quasi-coherent sheaf on a derived scheme $X$. Then,
\begin{enumerate}
\item $\Fc$ has \emph{Tor-amplitude} in the interval $[a, b] \subset \Zb$ if, for every discrete quasi-coherent sheaf $\Gc$, the homotopy sheaves $\pi_i(\Fc \otimes \Gc)$ are trivial whenever $i \not \in [a,b]$;

\item $\Fc$ has \emph{finite Tor-amplitude} if there exists a finite interval $[a, b] \subset \Zb$ in which $\Fc$ has its Tor-amplitude;

\item if $\Fc$ is connective, then we say that $\Fc$ has \emph{Tor-dimension} $i$ if it has Tor-amplitude in $[0,i]$;

\item if $\Fc$ is connective, then it is \emph{flat} if it has Tor-dimension 0.
\end{enumerate} 
\end{defn}

A family of concrete examples of a flat quasi-coherent sheaves is given by the locally free shaves. If $X$ is a derived scheme and $\Fc$ is a quasi-coherent sheaf on $X$, then $\Fc$ is \emph{locally free (of finite rank)} if there exists $n < \infty$, an open cover $(U_i)_{i \in I}$ of $X$, and equivalences $\Fc \vert_{U_i} \simeq \Oc_{U_i}^{\oplus n_i}$, where $n_i \leq n$. Moreover, if all the $n_i$ can be chosen to be equal to $n$, then $\Fc$ has \emph{rank} $n$.

\begin{defn}\label{def:perf}
Let $X$ be a derived scheme, and let $\Fc$ be a quasi-coherent sheaf on $X$. 
Then,
\begin{enumerate}
\item $\Fc$ is \emph{perfect} if, locally on each affine open $\Spec(A)$, $\Fc$ is a retract of a finite colimit of shifts of $A$; if $X$ is qcqs, then perfect objects are precisely the compact objects of $\QCoh(X)$;

\item $\Fc$ is \emph{almost perfect} if it is eventually connective, and for all $n \gg 0$ there exists a perfect object $\Gc$ and a map $\Gc \to \Fc$ inducing an isomorphisms $\pi_i(\Gc) \to \pi_i(\Fc)$ for $i \leq n$; if $X$ is Noetherian, then this is equivalent to $\Fc$ being eventually connective and the homotopy sheaves $\pi_i(\Fc)$ being coherent sheaves on $X_\cl$;


\item if $X$ is Noetherian, then $\Fc$ is \emph{coherent} if it is almost perfect and only finitely many of the homotopy sheaves are nontrivial.
\end{enumerate}
\end{defn}

On a qcqs scheme, an almost perfect object is perfect if and only if it has finite Tor-amplitude. Moreover, there exists a convenient characterization of those almost perfect sheaves that are locally free (clearly any locally free sheaf of finite rank is almost perfect).

\begin{prop}\label{prop:locfree}
Let $X$ be a Noetherian derived scheme. Then an almost perfect sheaf $\Fc$ is locally free if and only if it is flat.
\end{prop}

Another useful fact is that the Tor-dimension of an almost perfect sheaf can be measured using the following local procedure.

\begin{lem}\label{lem:dloctor}
Let $A$ be a Noetherian local derived ring with a residue field $\kappa$, and let $M$ be an almost perfect derived $A$-module (i.e., a quasi-coherent sheaf on $\Spec(A)$). Then $M$ has Tor-amplitude in $[a,b] \subset \Zb$ if and only if the homotopy modules $\pi_i(\kappa \otimes_A M)$ of $\kappa \otimes_A M$ vanish whenever $i \not \in [a,b]$.
\end{lem}
\begin{proof}
As any discrete $A$-module can naturally be regarded as a $\pi_0(A)$-module, it is enough to compute the Tor-amplitude of $\pi_0(A) \otimes_A M$. The claim then follows from standard facts concerning complexes of finitely generated modules over local Noetherian rings.
\end{proof}

If $X$ is a derived scheme, then we denote by $\Perf(X)$ the full subcategory of $\QCoh(X)$ spanned by the perfect objects. It is known that perfect objects coincide with the \emph{dualizable} objects in $\QCoh(X)$, i.e., to those quasi-coherent sheaves $\Fc$ that admit a dual $\Fc^\vee$ satisfying certain natural conditions. The most important properties of duals are listed below.

\begin{prop}\label{prop:duals}
Let $X$ be a derived scheme and let $\Fc \in \Perf(X)$. Then
\begin{enumerate}
\item for any $\Qc \in \QCoh(X)$, there exists a natural equivalence of spaces 
$$\Map_{\QCoh(X)}(\Fc, \Gc) \simeq |\Gamma(X; \Fc^\vee \otimes \Gc)|;$$

\item the double dual $\Fc^{\vee\vee}$ is naturally equivalent to $\Fc$;

\item if $\Fc$ is locally free (of rank $r$), then so is $\Fc^\vee$;

\item if $\Fc$ is a rank 1 locally free sheaf, then $\Fc \otimes \Fc^\vee \simeq \Oc_X$.
\end{enumerate} 
\end{prop}

We end by recalling that a morphism $\Ec \to \Fc$ of locally free sheaves is an \emph{inclusion} if the induced morphism of duals $\Fc^\vee \to \Ec^\vee$ is a surjection. As the fiber of a surjection of locally free sheaves is a vector bundle, the cofiber of an inclusion of locally free sheaves is locally free.

\subsection{The cotangent complex}\label{ssect:cotangent}

The cotangent complex is one of the most important objects in derived algebraic geometry. Much like the sheaf of Kähler-differentials, it admits an universal property as the target of the universal derived derivation, or, alternatively, as a quasi-coherent sheaf classifying infinitesimal derived extensions\footnote{Note that, in classical algebraic geometry, the cotangent complex does not admit a universal property in terms of deformation theory.}, as explained in a rather pleasant way in \cite{porta-vezzosi}. As any derived scheme is a limit of derived infinitesimal extensions that start with a classical scheme, the cotangent complex may be regarded as a key to understanding derived algebraic geometry.

As we will need only the most elementary properties of the cotangent complex, we will keep this subsection rather concise. Let us start by fixing the terminology.

\begin{defn}\label{def:cotangent}
Let $X \to Y$ be a morphism of derived schemes. Then, we will denote by $\Lb_{X/Y}$ the \emph{relative cotangent complex} of $X$ over $Y$. The cotangent complex is a connective quasi-coherent sheaf. The \emph{absolute cotangent complex} $\Lb_X$ is the relative cotangent complex associated to the morphism $X \to \Spec(\Zb)$. 
\end{defn}

Next, we record the useful properties of $\Lb_{X/Y}$.

\begin{prop}\label{prop:cotangent}
The cotangent complex satisfies the following properties:
\begin{enumerate}
\item if the square
$$
\begin{tikzcd}
X' \arrow[r] \arrow[d]{}{f} & Y' \arrow[d] \\
X \arrow[r] & Y
\end{tikzcd}
$$
is Cartesian, then there exists a natural equivalence $f^*\Lb_{X/Y} \simeq \Lb_{X'/Y'}$\footnote{Note that, in classical algebraic geometry, this is not true in general. If the classical Cartesian square is Tor-independent, then the result holds.};

\item given a sequence $X \xto{g} Y \to Z$ of derived schemes, then there exists a natural cofiber sequence
$$g^* \Lb_{Y/Z} \to \Lb_{X/Z} \to \Lb_{X/Y};$$

\item if the square
$$
\begin{tikzcd}
X' \arrow[d]{}{f} \arrow[r]{}{g} & Y' \arrow[d] \\
X \arrow[r] & Y
\end{tikzcd}
$$
is Cartesian, then the induced map $f^*\Lb_{X/Y} \oplus g^* \Lb_{Y' / Y} \to \Lb_{X' / Y}$ is an equivalence.
\end{enumerate}
\end{prop}

Finiteness properties of $f$ and of the relative cotangent complex are intimately related, as the following result (\cite{lurie-thesis} Proposition 3.2.18) shows.

\begin{prop}
Let $f: X \to Y$ be a map of derived schemes. Then the following are equivalent:
\begin{enumerate}
\item $f$ is locally of finite presentation (locally of almost finite presentation);
\item $f_\cl: X_\cl \to Y_\cl$ is finitely presented and $\Lb_{X/Y}$ is perfect (almost perfect). \qed
\end{enumerate}
\end{prop} 

\subsection{Vector bundles and related derived schemes}\label{ssect:dvbundle}

Much like in classical algebraic geometry, it is often useful to consider locally free sheaves as geometric objects.

\begin{defn}\label{def:dvbundle}
Let $X$ be a derived scheme, and let $\Ec$ on $X$ be a locally free sheaf. Then, the \emph{vector bundle} $E$ associated to $\Ec$ is the derived $X$-scheme such that, for every derived $X$-scheme $S$, there exists a natural equivalence
$$\Map_{\dSch_{/X}}(S , E) \simeq |\Gamma(S; \Ec \vert_S)|.$$
In particular, the space of sections of the structure map $E \to X$ is equivalent to the space of global sections of $\Ec$. The vector bundle of a locally free sheaf of rank $1$ is referred to as a \emph{line bundle}, and is often denoted by $\Ls$.
\end{defn}

One can check locally that the above functor of points description gives rise to a derived scheme. From now on, we will be rather careless in distinguishing vector bundles from locally free sheaves as they contain equivalent data. 

\begin{defn}\label{def:dvanish}
Let $X$ be a derived scheme, let $E$ be a vector bundle on $X$, and let $s$ be a global section of $E$. Then, the \emph{derived vanishing locus} of $s$ is defined as the fiber product
$$
\begin{tikzcd}
V_X(s) \arrow[hook,r] \arrow[hook,d] & X \arrow[d, hook]{}{0}\\
X \arrow[r, hook]{}{s}& E.
\end{tikzcd}
$$ 
Note that there exists a natural equivalence between the two inclusions $V_X(s) \hook X$ as both $0$ and $s$ are sections of $E \to X$. The derived vanishing locus of line bundle section of $\Ls$ is referred to as a \emph{virtual Cartier divisor}.
\end{defn}

Moreover, for each vector bundle, one can construct several interesting examples of derived schemes, namely the relative Grassmannians. 

\begin{defn}\label{def:grassmannian}
Let $X$ be a derived scheme and let $E$ be a vector bundle on $X$. Then the \emph{Grassmannian of $n$-planes in $E$}, $\Gr_n(E)$, is the derived $X$-scheme such that, for each derived $X$-scheme $S$, the space of $X$-morphisms $S \to \Gr_n(E)$ is equivalent to the space of inclusions 
$$F \hook  E \vert_S$$
of vector bundles on $S$, where $F$ is a vector bundle of rank $n$. The scheme $\Gr_1(E)$ is called the \emph{projective bundle} associated to $E$, and it is denoted by $\Pb(E)$\footnote{If $E$ is a vector bundle on $X$, then we will denote by $\Pbf(E)$ the universal scheme representing surjections $E \to \Ls$. However, we use this notation rarely.}.  Following the classical convention, the universal subbundle of $E$ on $\Pb(E)$ is denoted by $\Oc(-1)$. At times, it is useful to keep track of the base of the bundle, in which case we will use the notation $\Pb_X(E)$.
\end{defn}

Again, it can be verified that the above description gives rise to a derived scheme. Clearly an inclusion of vector bundles $E \hook F$ induces a closed embedding of Grassmannians $\Gr_n(E) \hook \Gr_n(F)$. In order to study this inclusion (in the special case of projective bundles) in more detail, we record the following useful result.  

\begin{lem}\label{lem:pbsect}
Let $X$ be a derived scheme and $E$ a vector bundle on $X$. Then, the natural cofiber sequence
$$\Oc(-1) \to E \to Q$$
of vector bundles on $\Pb(E)$ induces a natural equivalence $|\Gamma(X;E)| \simeq | \Gamma(\Pb(E);Q)|$ of spaces of global sections.
\end{lem}
\begin{proof}
Let $\pi$ be the structure morphism $\Pb(E) \to X$. Then one checks locally and by reducing to classical algebraic geometry that the natural map $\pi_*(\pi^*(E)) \to E$ is an equivalence. Similarly, $\pi_*(\Oc(-1)) \simeq 0$. The claim then follows from the fact that cofiber sequences push forward to cofiber sequences.
\end{proof}

We now show that projectivized linear embeddings are derived vanishing loci.

\begin{prop}\label{prop:pbundvan}
Let $E \hook F$ be an inclusion of vector bundles, and let us denote the cofiber of the inclusion by $G$. Then, the inclusion $\Pb(E) \hook \Pb(F)$ is the derived vanishing locus of the global section $s$ of $G(1)$ that corresponds to the composition
$$\Oc(-1) \hook F \to G$$
of maps of vector bundles on $\Pb(F)$.
\end{prop}
\begin{proof}
The space of $X$-morphism from a derived $X$-scheme $S$ to the derived vanishing locus $V_{\Pb(F)}(s)$ is equivalent to the space of commuting diagrams,
$$
\begin{tikzcd}
\Ls \arrow[d, hook] \arrow[rd]{}{0} \\
F \vert_S \arrow[r] & G \vert_S.
\end{tikzcd}
$$ 
As $E \hook F \to G$ is a fiber sequence in quasi-coherent sheaves, the claim follows.
\end{proof}

\section{Quasi-projectivity and ample line bundles}\label{sect:qproj}

Here, we study ample line bundles and quasi-projective morphisms in derived algebraic geometry. These results originally appeared in the background section of \cite{annala-base-ind-cob} and slightly earlier in the unpublished note \cite{annala-qpnote}. Of course, the results are rather straightforward analogues of well-known results in classical algebraic geometry.

\subsection{Ample line bundles}\label{ssect:ample}

The purpose of this subsection to define and study ample line bundles in the context of derived algebraic geometry. 

\begin{defn}\label{def:ample}
Let $X$ be a Noetherian derived scheme of finite Krull dimension and let $\Ls$ be a line bundle on $X$. Then $\Ls$ is \emph{ample} if, for any point $x \in X$, there exists $n \geq 0$ and a section $s \in \Gamma(X; \Ls^{\otimes n})$ such that $X_s := X \backslash V(s)$ is affine and contains $x$.
\end{defn}

The following result will be important when analyzing the implications of the above definition.

\begin{lem}\label{lem:ampleloc}
Let $X$ be a finite-Krull-dimensional Noetherian derived scheme, let $\Ls$ be a line bundle on $X$ and let $s$ be a global section of $\Ls$. Let $\Fc$ be a quasi-coherent sheaf on $X$, and consider the sequences
$$H^i(X; \Fc) \xrightarrow{s \cdot} H^i(X; \Ls \otimes \Fc) \xrightarrow{s \cdot} H^i(X; \Ls^{\otimes 2} \otimes \Fc) \xrightarrow{s \cdot} \cdots$$
Then there are natural isomorphisms
$$\colim_n H^i \big(X; \Ls^{\otimes n} \otimes \Fc \big) \xrightarrow{\cong} H^i(X_s; \Fc)$$
given by sending $f \in H^i(X; \Ls^{\otimes n} \otimes \Fc)$ to $f/s^n \in H^i(X_s; \Fc)$.
\end{lem}
\begin{proof}
It is enough to prove this for $i=0$. Since applying the truncation $\tau_{\geq 0}$ does not change the zeroth cohomology, we can assume that $\Fc$ is connective. Finally, by Grothendieck vanishing (Lemma \ref{lem:gvanish}) we may reduce to the situation where the nontrivial homotopy sheaves of $\Fc$ lie between degrees $0$ and $n$, where $n$ is at most the Krull dimension of $X$. Our argument proceeds by induction on $n$. Notice that the base case $n=0$ is classical (see e.g. \cite{stacks} Tag 09MR). To prove the general case, we consider the cofiber sequence
$$\Fc' \to \Fc \to \pi_0(\Fc).$$
Since the nontrivial homotopy sheaves of $\Fc'$ lie between degrees $1$ and $n$, we can apply the induction assumption to $\Fc'$ and $\pi_0(\Fc)$, and the claim follows from 5-lemma applied to the diagram comparing the cohomology long exact sequences (note that sequential colimits preserve exact sequences).
\end{proof}

Ample line bundles admit the following useful characterization.

\begin{prop}\label{prop:amplechar}
Let $X$ be a finite-Krull-dimensional Noetherian derived scheme, and let $\Ls$ be a line bundle on $X$. Then, the following conditions are equivalent:
\begin{enumerate}
\item $\Ls$ is ample;
\item for every connective almost perfect sheaf $\Fc$ on $X$, the sheaves $\Ls^{\otimes n} \otimes \Fc$ are globally generated for $n \gg 0$.
\end{enumerate}
\end{prop}
\begin{proof}
By Grothendieck vanishing (Lemma \ref{lem:gvanish}), we can assume that $\Fc$ is coherent. Let us first assume that $\Ls$ is ample, and choose $s_i \in \Gamma(X; \Ls^{\otimes n_i})$, $i=1..r$, in such a way that $X_{s_i}$ form an affine open cover of $X$. As 
$$\pi_0\big(\Gamma(X_{s_i}; \Fc)\big) \cong \Gamma\big(X_{s_i}; \pi_0(\Fc) \big),$$
we can use Lemma \ref{lem:ampleloc} to conclude that, for $d_i \gg 0$, there exist global sections of $\Ls^{\otimes d_i n_i} \otimes \Fc$ generating it at the points of $X_{s_i}.$ Letting $N = n_1 \cdots n_r$, it then follows that $\Ls^{\otimes iN} \otimes \Fc$ is globally generated for $i \gg 0$. Since the last conclusion is also true for the sheaves $\Ls \otimes \Fc, ..., \Ls^{\otimes N-1} \otimes \Fc$, it follows that $\Ls^{\otimes n} \otimes \Fc$ is globally generated for $n \gg 0$.

Suppose then that, for all coherent sheaves $\Fc$, $\Ls^{\otimes n} \otimes \Fc$ is globally generated for $n \gg 0$. Let $x \in X$ be a point, and let $U$ be an affine open neighborhood of $x$ such that $\Ls \vert_U$ is trivial. Let $Z \hook X$ be the complement of $U$ with the reduced scheme structure, and consider the cofiber sequence
$$\Ic \to \Oc_X \to \Oc_Z.$$
If  $\Ls^{\otimes n} \otimes \Ic$ is globally generated, then it is possible to find such a global section $s$ that its image in $s' \in \Gamma(X; \Ls^{\otimes n})$ does not vanish at $x$. It follows that $X_{s'} \subset U$ is an affine open subset containing $x$.
\end{proof}

Next, we consider the relative version of ampleness.

\begin{defn}\label{def:relample}
Let $f: X \to Y$ be a morphism of finite-Krull-dimensional Noetherian derived schemes. Then, a line bundle $\Ls$ on $X$ is \emph{relatively ample} (or \emph{$f$-ample}) if, for every affine open set $\Spec(B) \to Y$, $\Ls$ restricts to an ample line bundle on the inverse image $f^{-1} \Spec(B)$.
\end{defn}

If $Y$ is affine, then a line bundle $\Ls$ on $X$ is relatively ample over $Y$ if and only if $\Ls$ is ample. Combining this observation with the following lemma allows us to check relative ampleness on an affine cover of $Y$.

\begin{lem}\label{lem:relampleloc}
Let $f: X \to Y$ be a morphism of finite-Krull-dimensional Noetherian derived schemes, let $(U_i)_{i \in I}$ be an open cover of $Y$, and let $\Ls$ be a line bundle on $X$.  Then, $\Ls$ is $f$-ample if and only if for all $i \in I$, we have that $\Ls_i := \Ls \vert_{f^{-1} U_i}$ is $f_i$-ample, where $f_i: f^{-1} U_i \to U_i$ is the restriction of $f$.
\end{lem}
\begin{proof}
The only if direction is clear, so let us prove that all $\Ls_i$ being $f_i$-ample implies the relative ampleness of $\Ls$. In other words we have to prove the following: if $Y = \Spec(B)$ is affine, $b_1, ..., b_r \in B$ are functions such that $(Y_{b_i})_{i=1}^r$ is an affine open cover of $Y$ with $\Ls \vert_{X_{b_i}}$ ample for all $i$, then $\Ls$ is ample. But this is obvious: if $x \in X$, then $x$ lies inside one of the open subsets $X_{b_i}$, and we may use Lemma \ref{lem:ampleloc} to find a global section $s$ of $\Ls^{\otimes n}$ with $X_{s} \subset X_{b_i}$ affine and containing $x$.
\end{proof}

As a consequence, we prove that relatively ample line bundles behave well in compositions and in derived pullbacks.

\begin{prop}\label{prop:relamplecomppull}
Let $f: X \to Y$ be a morphism of finite-Krull-dimensional Noetherian derived schemes having a relatively ample line bundle $\Ls$. Then,
\begin{enumerate}
\item if $f$ factors through $f': X \to Y'$, $\Ls$ is $f'$-ample;

\item if $g: Y \to Z$ admits a relatively ample line bundle $\Ms$, then $\Ls \otimes f^* \Ms^{\otimes n}$ is $(g \circ f)$-ample for $n \gg 0$;
\item if 
$$
\begin{tikzcd}
X' \arrow[]{r}{f'} \arrow[]{d}{h'} & Y' \arrow[]{d}{h} \\
X \arrow[]{r}{f} & Y
\end{tikzcd}
$$
is Cartesian with $X'$ and $Y'$ Noetherian, then $h'^* \Ls$ is $f'$-ample.
\end{enumerate}
\end{prop}
\begin{proof}
\begin{enumerate}
\item This follows from Lemma \ref{lem:relampleloc} and from the stability of ample line bundles under restrictions to open subschemes.

\item The proof is the same as the classical proof, see \cite{stacks} Tags 0C4K (1) and 0892.

\item Without loss of generality we can assume that both $Y$ and $Y'$ are affine. It follows that $h'$ is affine, and as ample line bundles are stable under affine pullbacks, the claim follows. \qedhere
\end{enumerate}
\end{proof}

It is also true that relative ampleness can be checked on the truncation, at least for proper morphisms. 

\begin{prop}\label{prop:proprelample}
Let $A$ be a Noetherian derived ring of finite Krull dimension, and let $X \to \Spec(A)$ be a proper morphism of derived schemes that is locally of finite presentation. Then, for a line bundle $\Ls$ on $X$, the following conditions are equivalent:
\begin{enumerate}
\item $\Ls_\cl$ is ample on $X_\cl$;
\item given an almost perfect quasi-coherent sheaf $\Fc$ on $X$ and $i \in \Zb$, we have that
$$\pi_i \big( \Gamma(X; \Fc \otimes \Ls^{\otimes n}) \big) \cong \Gamma_\cl(X; \pi_i(\Fc \otimes \Ls^{\otimes n}))$$
for $n \gg 0$;
\item $\Ls$ is ample.
\end{enumerate}
\end{prop}
\begin{proof}
Let us first assume that the truncation $\Ls_\cl$ is ample. By Grothendieck vanishing, it is enough to ensure that condition $2.$ holds for $\Fc$ coherent and connective. We will proceed by induction on $n$, where $n$ is the largest $i$ such that $\pi_i(\Fc)$  is nontrivial. The base case $n=0$ is easy: in this situation
$$\Fc \otimes_{\Oc_X} \Ls^{\otimes m} \simeq \Fc \otimes_{\Oc_{X_\cl}} \Ls_\cl^{\otimes m},$$
and the claim follows from the classical result of Serre on the vanishing of sheaf cohomology (see, e.g., \cite{stacks} Tag 0B5U). In general we have a cofibre sequence
$$\Fc' \to \Fc \to \pi_0(\Fc).$$
Applying the induction assumption to $\Fc'$ and $\pi_0(\Fc)$, the result holds for them for $m$ large enough. Then for any such $m$, and $n>m$, it follows from the associated cohomology long exact sequence that
\begin{align*}
&\pi_i\big(\Gamma(X; \Fc \otimes \Ls^{\otimes n})\big) \cong \\
&
\begin{cases}
\pi_i\big(\Gamma(X; \Fc' \otimes \Ls^{\otimes n})\big) \cong \Gamma_\cl(X; \pi_i(\Fc \otimes \Ls^{\otimes n})) & \text{for $i > 0$;} \\
\Gamma_\cl(X; \pi_0(\Fc \otimes \Ls^{\otimes n}))\big) & \text{for $i = 0$;} \\
0 & \text{for $i < 0$.}
\end{cases}
\end{align*}
This shows that $\Ls$ satisfies condition 2.

Suppose then that $\Ls$ satisfies condition 2. Applying this to all truncated coherent sheaves, we can immediately conclude that $\Ls_\cl$ is ample using the result of Serre cited above. As global sections of $\Ls^{\otimes n}_\cl$ lift to global sections of $\Ls^{\otimes n}$ for $n \gg 0$, the ampleness of $\Ls_\cl$ implies that of $\Ls$.
\end{proof}

\subsection{Quasi-projective morphisms}\label{ssect:qproj}

Here, we define and study the basic properties of quasi-projective morphisms. Let us begin with explaining what we mean by a quasi-projective morphism.

\begin{defn}\label{def:qproj}
A locally almost finitely presented morphism $f: X \to Y$ of finite-Krull-dimensional Noetherian derived schemes is \emph{quasi-projective} if it admits a relatively ample line bundle. A proper quasi-projective morphism is called \emph{projective}\footnote{In \cite{annala-base-ind-cob} we did not assume that $f$ is locally of almost finite presentation, only of finite type. However, in order to ensure the compatibility of our algebraic cobordism with $K$-theory, we need to fix a class of projective morphisms with the property that pushforward of an almost perfect object is almost perfect. For this, we need to assume that $f$ is locally of almost finite presentation \cite{SAG}. This definition has also the added benefit that a quasi-projective derived scheme over a Noetherian derived ring $A$ is automatically Noetherian.}.
\end{defn}

In light of Proposition \ref{prop:proprelample}, we may generalize the notion of projective morphism further to get rid of the Noetherian hypothesis.

\begin{defn}
Let $f: X \to Y$ be a proper morphism of derived schemes that is locally of almost finite presentation. Then, $f$ is \emph{projective} if there exists a line bundle $\Ls$ on $X$ such that $\Ls_\cl$ is $f_\cl$-ample.
\end{defn}

We have chosen to use the above definitions for (quasi-)projectivity, as it has the following good properties.

\begin{prop}\label{prop:qprojstab}
Quasi-projective morphisms between Noetherian derived schemes of finite Krull dimension are stable under compositions and pullbacks in which the fiber product is Noetherian. Projective morphisms are stable under pullbacks and compositions.
\end{prop}
\begin{proof}
The first claim follows immediately from Proposition \ref{prop:relamplecomppull}. The second claim follows immediately from the analogous fact in classical algebraic geometry, see e.g. \cite{stacks} Tags 0C4K and 0893.
\end{proof}

Quasi-projective morphisms have the benefit over more general morphisms that they often admit convenient global factorizations.

\begin{thm}\label{thm:qprojfact}
Let $f: X \to Y$ be a quasi-projective morphism of finite-Krull-dimensional Noetherian derived schemes. Then, if $Y$ admits an ample line bundle, $f$ admits a factorization
$$X \stackrel i \hook \Pb^n_Y \stackrel \pi \to Y,$$
where $i$ is a locally closed embedding and $\pi$ is the natural projection. If $f$ is projective, then $i$ is a closed embedding.
\end{thm}
\begin{proof}
By assumption, $X$ admits an ample line bundle $\Ls$; let $b_1,...,b_r$ be a sequence of global sections of $\Ls$ such that $X_{b_i}$ is an affine open cover for $X$. As ample line bundles truncate to ample line bundles, there exists $n > 0$ and a generating sequence of global sections $s_0,...,s_r$ of $\Ls^{\otimes n}_\cl$ such that the induced morphism $X_\cl \hook \Pb^r \times Y_\cl$ is a locally closed embedding (\cite{stacks} Tag 0B42). Then by Lemma \ref{lem:ampleloc}, there exists such an $m>0$ that the global sections $s_{ij} := b_i^{\otimes m} \otimes s_j$ of $\Ls^{\otimes (n+m)}_\cl$ lift to global sections $\tilde s_{ij}$ of $\Ls^{\otimes (n+m)}$. As the morphism $i_\cl: X_\cl \hook \Pb^n \times Y_\cl$ induced by the $s_{ij}$ is still a locally closed embedding, the map $i: X \hook \Pb^n \times Y$ induced by the $s_{ij}$ is a locally closed embedding of derived schemes. Thus, we have proven the first claim. The second claim follows from the fact that $i$ is proper.
\end{proof}

\section{Derived complete intersections}\label{sect:dci}

This section is dedicated to the study of complete intersections in derived algebraic geometry. We start with Section \ref{ssect:qsm}, in which we recall the basic facts concerning quasi-smooth and smooth morphisms in derived algebraic geometry. In Section \ref{ssect:dci}, we study the corresponding absolute notion, i.e., that of a derived complete intersection ring. The results of this subsection were originally stated in the background section of \cite{annala-base-ind-cob}.

\subsection{Smooth and quasi-smooth morphisms}\label{ssect:qsm}

In this section, we study smooth and quasi-smooth morphisms and related notions in derived algebraic geometry. We say that a morphism $f: X \to Y$ of derived schemes is \emph{locally of finite Tor-dimension}, if, for each affine open $\Spec(A) \subset Y$ and an affine open $\Spec(B) \subset X$ mapping to $\Spec(A)$, $B$ has finite Tor-dimension as a derived $A$-module. 

Let us begin with the following definitions.

\begin{defn}\label{defn:sm}
Let $f: X \to Y$ be a morphism of derived schemes that is locally of finite Tor-dimension. Then, $f$ is \emph{formally smooth} if the relative cotangent complex $\Lb_{X/Y}$ has Tor-dimension 0\footnote{We have made being locally of finite Tor-dimension part of the definition of a formally (quasi-)smooth morphism in order to ensure that bivariant $K$-theory admits orientations along such morphisms. Perhaps, at least under some fairly weak additional hypothesis, this would follow from the cotangent complex having Tor-dimension $\leq 1$. See e.g. \cite{avramov:1999} for some classical results in this direction. In any case, this condition may be dropped from the definition of a (quasi-)smooth morphism, as it is automatic under the other hypotheses.}. A formally quasi-smooth morphism morphism $f$ is \emph{smooth} if it is locally of (almost) finite presentation.
\end{defn}

\begin{defn}\label{defn:qsm}
Let $f: X \to Y$ be a morphism of derived schemes that is locally of finite Tor-dimension. Then, $f$ is \emph{formally quasi-smooth} if the relative cotangent complex $\Lb_{X/Y}$ has Tor-dimension 1. A formally quasi-smooth morphism is \emph{quasi-smooth} if it is locally of (almost) finite presentation. A \emph{derived regular embedding} is a quasi-smooth closed immersion.
\end{defn}

As the cotangent complex of a quasi-smooth morphism is perfect, quasi-smooth morphisms admit a notion of relative dimension. 

\begin{defn}\label{def:reldim}
Let $f: X \to Y$ be a quasi-smooth morphism, and assume that the relative cotangent complex $\Lb_{X/Y}$ has constant virtual rank $r$. Then, $f$ is said to have \emph{relative virtual dimension $r$}. If $i: Z \hook X$ is a quasi-smooth closed immersion of relative virtual dimension $r$, then we say that it has \emph{virtual codimension} $-r$.
\end{defn}

\begin{rem}\label{rem:reldim}
The relative virtual dimension is well-defined on each connected component of the source. Hence, the only problem for defining it in general is when the source is a disjoint union of derived schemes that have different relative virtual dimensions over the target.
\end{rem}

The classes of morphisms defined above have many good properties, as summarized in the following result.

\begin{prop}\label{prop:qsmsprop}
Formally quasi-smooth morphisms (formally smooth morphisms, quasi-smooth morphisms, smooth morphisms, derived regular embeddings) satisfy the following properties:
\begin{enumerate}
\item they are stable under compositions and derived pullbacks\footnote{Note that formally lci morphisms, lci morphisms and regular embeddings are not stable under classical pullbacks. This is the main reason why it seems necessary to use derived geometry to construct a good theory of bivariant algebraic cobordism.};

\item whenever the relative virtual dimension is well defined, it is stable under pullbacks and additive in compositions;

\item a morphism $f: X \to Y$ of classical schemes is smooth (quasi-smooth, derived regular embedding) if and only if it is smooth (lci, regular embedding);

\item whenever a quasi-smooth morphism factors as
$$X \xto{f} P \xto{p} Y,$$
where $p$ is smooth, then $f$ is quasi-smooth.
\end{enumerate}
\end{prop}
\begin{proof}
All the claims follow directly from the theory of cotangent complex in derived and classical algebraic geometry and the properties of Tor-dimension.
\end{proof}

Derived regular embeddings admit a local characterization as the derived vanishing locus of a sequence of functions on the target. To state this more precisely, we recall the following definition from \cite{khan-rydh}.

\begin{defn}\label{def:dvanloc}
Let $A$ be a derived ring, and let $a_1, ..., a_r \in A$. Then the \emph{derived quotient} $A \modmod ( a_1,...,a_r )$ is defined as the tensor product $A \otimes_{\Zb[x_1,...,x_r]} \Zb$ of derived rings, where the maps $\Zb[x_1,...,x_r] \to A$ and $\Zb[x_1,...,x_r]$ are defined as $x_i \mapsto a_i$ and $x_i \mapsto 0$, respectively.
\end{defn}

\begin{ex}\label{ex:dvanlov}
Derived vanishing loci of vector bundles are locally described as opposites of derived quotients. Moreover, $\Spec(A \modmod ( a_1,...,a_r )) \to \Spec(A)$ is equivalent to the derived vanishing locus of the section $( a_1,...,a_r )$ of the trivial vector bundle $\Oc_{\Spec(A)}^{\oplus r}$.
\end{ex}

The following result is Proposition 2.3.8 of \cite{khan-rydh}.

\begin{prop}\label{prop:regembchar}
A closed embedding $i: Z \to X$ of schemes is a regular embedding if and only if, Zariski locally on $X$, it is equivalent to the inclusion of derived vanishing locus of a sequence of functions on $X$.
\end{prop}

Moreover, derived regular embeddings of codimension 1 are virtual Cartier divisors, i.e., vanishing loci of sections of line bundles (this is a consequence of Proposition 3.2.6 of \cite{khan-rydh}).

\begin{prop}\label{prop:vcartchar}
Let $i: D \hook X$ be a derived regular embedding of codimension $1$. Then there exists a line bundle $\Oc_X(D)$ on $X$, a global section $s_D \in |\Gamma(X, \Oc_X(D))|$, and an equivalence of derived $X$-schemes $D \simeq V_X(s_D)$. 
\end{prop}

The above result allows us to define sum of virtual Cartier divisors.

\begin{defn}\label{def:vcartsum}
Let $D_1$ and $D_2$ be virtual Cartier divisors on a derived scheme $X$. Then their \emph{sum}, $D_1 + D_2$, is the virtual Cartier divisor that is obtained as the derived vanishing locus of the global section $s_{D_1} \otimes s_{D_2}$ of $\Oc_X(D_1) \otimes \Oc_X(D_2)$.
\end{defn}

\subsection{Derived complete intersection schemes}\label{ssect:dci}

Here, we study an absolute counterpart of quasi-smoothness.

\begin{defn}\label{def:dci}
A Noetherian local derived ring $A$ is called \emph{derived complete intersection} if only finitely many of the homotopy modules $\pi_i(A)$ are non zero, and if the cotangent complex $\Lb_{\kappa/A}$ is perfect and concentrated in degrees 1 and 2. A Noetherian derived scheme $X$ is \emph{derived complete intersection} (or \emph{derived regular}) if all of its local derived rings are derived complete intersection. 
\end{defn} 

\begin{ex}\label{ex:dci}
If $A$ is a regular local ring and $a_1,...,a_r$ is a sequence of elements of $A$, then $A \modmod (a_1,...,a_r)$ is a derived complete intersection ring.
\end{ex}

Definition \ref{def:dci} captures the notion of derived local rings that are derived quotients of regular local rings, as shown by the next result.

\begin{lem}\label{lem:dcitoreg}
If a derived complete intersection ring $B$ admits a surjection $A \to B$ from a regular local ring $A$, then there exists elements $a_1,...,a_r$ and an equivalence $B \simeq A \modmod(a_1,...,a_r)$ of $A$-algebras. In other words, a closed immersion from a derived regular scheme to a regular scheme is a derived regular embedding.
\end{lem}
\begin{proof}
Let $\kappa$ be the residue field of $B$ and consider the sequence $A \to B \to \kappa$. Then, from the fundamental triangle
$$\kappa \otimes_B \Lb_{B/A} \to \Lb_{\kappa/A} \to \Lb_{\kappa / B}$$
we can deduce that $\kappa \otimes_B \Lb_{B/A}$ is concentrated in degree 1. Thus $\Lb_{B/A}$ has Tor-dimension $1$ and the claim follows from \ref{prop:regembchar}.
\end{proof}

This allows us to prove the following useful characterization of derived regular schemes.

\begin{prop}\label{prop:dregchar}
Let $X$ be a Noetherian scheme that admits a finite-type morphism $f: X \to Y$ to a regular Noetherian scheme $Y$. Then $X$ is derived complete intersection if and only if $X \to Y$ is quasi-smooth. 
\end{prop}
\begin{proof}
Locally, $f$ factors as $X \stackrel i \hook \Ab^n \times Y \to Y$, where the first morphism is a closed embedding. As $f$ is quasi-smooth if and only if $i$ is a derived regular embedding (Proposition \ref{prop:qsmsprop}), the claim follows from Lemma \ref{lem:dcitoreg}. 
\end{proof}

\begin{war}\label{war:dregchar}
In \cite{annala-base-ind-cob}, Proposition \ref{prop:dregchar} was stated without any finite-type assumptions. The proof of this is not correct, and the statement may be false. In the paper, the result was used to identify the absolute algebraic bordism of $X$ with any bivariant group of form $\Omega(X \to Y)$ where $Y$ is regular ($Y \cong \Spec(\Zb)$ being a canonical choice). Using the above result, algebraic bordism may be identified with such bivariant groups only when $X \to Y$ is of finite type.   
\end{war}

We end with the following useful result.

\begin{prop}\label{prop:dregqsm}
Let $X$ be a derived regular scheme, and $V \to X$ a formally quasi-smooth morphism. If $V$ is Noetherian, then $V$ is a derived regular scheme.
\end{prop}
\begin{proof}
Without loss of generality, we may assume that $X \simeq \Spec(A)$ and $V \simeq \Spec(B)$, where $A$ is a derived complete intersection ring and $B$ is a derived local ring with residue field $\kappa$. Considering the fundamental triangle
$$\kappa \otimes_B \Lb_{B/A} \to \Lb_{\kappa / A} \to \Lb_{\kappa/B},$$
the claim follows if we can show that $\Lb_{\kappa / A}$ is concentrated in degrees $\leq 1$. Note that the map $A \to \kappa$ factors through $A \to \kappa '$, where $\kappa '$ is the residue field of $A$. Considering the cofiber sequence
$$\kappa \otimes_{\kappa'} \Lb_{\kappa' / A} \to \Lb_{\kappa / A} \to \Lb_{\kappa / \kappa'},$$
we observe that it suffices to show that $\Lb_{\kappa / \kappa'}$ is concentrated in degrees $\leq 1$. But this is always the case, as one can see by factoring the field extension into a separable extension followed by a purely inseparable one, and computing the cotangent complex.
\end{proof}

\section{Derived blowups}\label{sect:blowpup}

One of the main technical tools used in the study of derived algebraic cobordism is the derived blowup construction of Khan and Rydh \cite{khan-rydh}. Here, we recall the basic facts related to this construction. The results following Theorem \ref{thm:dblowupprop} are compiled from \cite{annala-yokura, annala-chern}.

\begin{defn}\label{def:vcartdivover}
Let $Z \hookrightarrow X$ be a derived regular embedding. Then, for any $X$-scheme $S$, a \emph{virtual Cartier divisor} on $S$ \emph{lying over $Z$} is the datum of a commutative diagram
\begin{center}
\begin{tikzcd}
D \arrow[hook]{r}{i_D} \arrow[]{d}{g} & S \arrow[]{d}{} \\
Z \arrow[hook]{r}{} & X
\end{tikzcd}
\end{center}
such that
\begin{enumerate}
\item $i_D$ is a derived regular embedding of virtual codimension 1;
\item the truncation is Cartesian square in the ordinary category of schemes;
\item the canonical morphism
$$g^* \Nc^\vee_{Z/X} \to \Nc^\vee_{D/S}$$
is surjective.
\end{enumerate}
\end{defn}

The universal property of derived blowups can be phrased in terms of virtual Cartier divisors lying over the center.

\begin{defn}\label{def:dblowup}
Let $Z \hookrightarrow X$ be a derived regular embedding. Then the \emph{derived blowup} $\bl_Z(X)$ is the derived $X$-scheme representing virtual Cartier divisors lying over $Z$. In other words, given a derived $X$-scheme $S$, the space of $X$-morphisms 
$$S \to \bl_Z(X)$$
is naturally identified with the space of virtual Cartier divisors of $S$ that lie over $Z$.
\end{defn}

Derived blowups satisfy the following basic properties.

\begin{thm}\label{thm:dblowupprop}
Let $i: Z \hookrightarrow X$ be a derived regular embedding. Then,
\begin{enumerate}
\item the derived blowup $\bl_Z(X)$ exists and is unique up to contractible space of choices;

\item the structure morphism $\pi: \bl_Z(X) \to X$ is projective, quasi-smooth, and induces an equivalence 
$$\bl_Z(X) - \Ec \to X - Z,$$
where $\Ec$ is the universal virtual Cartier divisor on $\bl_Z(X)$ lying over $Z$ (also called the \emph{exceptional divisor});

\item the derived blowup $\bl_Z(X) \to X$ is stable under derived base change;

\item in the blowup square
$$
\begin{tikzcd}
\Ec \arrow[d]{}{g} \arrow[hook, r]{}{i_\Ec} & \bl_Z(X) \arrow[d] \\
Z \arrow[hook, r] & X
\end{tikzcd}
$$
the induced surjection $\Nc_{Z/X}^\vee \to \Nc^\vee_{\Ec/\bl_{Z}(X)}$ identifies $\Ec$ with $\Pb(\Nc_{Z/X}^\vee)$; hence $\Nc^\vee_{\Ec/\bl_{Z}(X)} \simeq \Oc(1)$;

\item if $Z \stackrel i \hookrightarrow X \stackrel j \hookrightarrow Y$ is a sequence of quasi-smooth closed embeddings, then there exists a canonical derived regular embedding $\tilde j: \bl_Z(X) \hookrightarrow \bl_Z(Y)$ called the \emph{strict transform};


\item if $Z$ and $X$ are classical schemes, then there exists a natural equivalence 
$$\bl_Z(X) \simeq \bl_{Z}^\mathrm{cl}(X),$$
where the right hand side is the classical blowup.
\end{enumerate}
\end{thm}
\begin{proof}
Every claim expect the projectivity of $\pi: \bl_Z(X) \to X$ follows directly from \cite{khan-rydh} Theorem 4.1.5 (see also \emph{ibid.} 4.3.4 and 4.3.5 for the fourth claim). As Khan--Rydh prove $\pi$ to be proper, it suffices to find a line bundle on the blowup that  truncates to a relatively ample line bundle in the usual sense. Clearly, $\Oc(-\Ec)$ is such a line bundle.
\end{proof}

For the sequel, it will be useful to understand the behavior of strict transforms on the exceptional divisors.

\begin{lem}\label{lem:linearblowup}
Consider a sequence $Z \hook X \hook Y$ of derived regular embeddings.  Then the derived pullback of the strict transform $\bl_Z(X) \hookrightarrow \bl_Z(Y)$ to the exceptional divisor of $\bl_Z(Y)$ is equivalent to the projectivized inclusion $\Pb(\Nc_{Z/X}) \hook \Pb(\Nc_{Z/Y})$.
\end{lem}
\begin{proof}
Recall that the strict transform is induced by the outer square in
$$
\begin{tikzcd}
\Pb(\Nc_{Z/X}) \arrow[hookrightarrow]{r} \arrow[->]{d}{g} & \bl_Z(X) \arrow[->]{d} \\ 
Z \arrow[hookrightarrow]{r} \arrow[->]{d}{\mathrm{Id}_Z} & X \arrow[->]{d} \\
Z \arrow[hookrightarrow]{r} & Y
\end{tikzcd}
$$
where the upper square is the blowup square associated to $\bl_Z(X)$. By the basic functoriality properties of the cotangent complex, the morphism $g^* \Nc_{Z/Y}^\vee \to \Nc_{\Pb(\Nc_{Z/X})/\bl_Z(X)}^\vee$ is naturally identified with the composition
$$g^* \Nc_{Z/Y}^\vee \to g^* \Nc_{Z/X}^\vee \to \Nc_{\Pb(\Nc_{Z/X})/\bl_Z(X)}^\vee,$$
which clearly gives rise to the desired morphism.
\end{proof}

The above result allows us to prove the following useful result about blowing up intersections of derived regular embeddings. 

\begin{prop}\label{prop:dblowupofint}
Let $X$ be a Noetherian derived scheme, and let $Z_1 \hook X$ and $Z_2 \hook X$ be derived regular embeddings. Let us denote by $Z_{12} \hook X$ the inclusion of the intersection of $Z_1$ and $Z_2$ inside $X$. Then,
\begin{enumerate}
\item the exceptional divisor 
$$\Pb(\Nc_{Z_{12} / X}) \simeq \Pb(\Nc_{Z_1/X} \vert_{Z_{12}} \oplus \Nc_{Z_2/X}\vert_{Z_{12}}) \hook \bl_{Z_{12}}(X)$$
meets the strict transform $\bl_{Z_{12}}(Z_i) \hook \bl_{Z_{12}}(X)$ in $\Pb(\Nc_{Z_{3-i}/X} \vert_{Z_{12}})$;

\item the strict transforms $\bl_{Z_{12}}(Z_1) \hook \bl_{Z_{12}}(X)$ and $\bl_{Z_{12}}(Z_2) \hook \bl_{Z_{12}}(X)$ have an empty intersection.
\end{enumerate}
\end{prop}
\begin{proof}
\begin{enumerate}
\item As $\Nc_{Z_{12} / Z_i} \simeq \Nc_{Z_{3-i} / X} \vert_{Z_{12}}$, the claim follows immediately from Lemma \ref{lem:linearblowup}.
 
\item Outside the exceptional divisor of $\bl_{Z_{12}}(X)$, the strict transforms are clearly disjoint. The claim then follows from the fact that the linear subbundles $\Pb(\Nc_{Z_{1} / X} \vert_{Z_{12}})$ and $\Pb(\Nc_{Z_{2} / X} \vert_{Z_{12}})$ have an empty intersection inside the bundle $\Pb(\Nc_{Z_{1} / X} \vert_{Z_{12}} \oplus \Nc_{Z_{2} / X} \vert_{Z_{12}})$. \qedhere
\end{enumerate}
\end{proof}

Next, we analyze the behavior of pullbacks of divisors in derived blowups. This will be of fundamental importance later on, as it leads to the derived blow-up formula in cobordism, which in turn is employed in almost all cobordism computations. We start with the following preliminary result.

\begin{lem}\label{lem:vdivsumint}
Let $D_1$ and $D_2$ be virtual Cartier divisors on a derived scheme $X$. Then there exists a natural commutative square
$$
\begin{tikzcd}
D_1 \times_X D_2 \arrow[r] \arrow[d]{}{g} & D_2 \arrow[d] \\
D_1 \arrow[r] & D_1 + D_2
\end{tikzcd}
$$
of derived $X$-schemes. Moreover, the map $\Lb_{D_1/D_1 + D_2}[-1] \to \Nc^\vee_{D_1 \times_X D_2 / D_2}$ is a surjection.
\end{lem}
\begin{proof}
The maps $D_i \to D_1 + D_2$ are induced by the paths $s_1 \otimes \gamma_2$ and $\gamma_1 \otimes s_2$ in the spaces of global sections $|\Gamma(D_i; \Oc(D_1 + D_2))|$, where $\gamma_i$ is the canonical path $s_i \sim 0$ in $|\Gamma(D_i; \Oc(D_i))|$. Moreover, the product $\gamma_1 \otimes \gamma_2$ gives a homotopy in $|\Gamma(D_1 \times_X D_2; \Oc(D_1 + D_2))|$ between the paths $\gamma_1 \otimes s_2$ and $s_1 \otimes \gamma_2$, hence equipping the square in the statement with this homotopy makes it commutative, proving the first claim. The second claim follows from the fact that the composition 
$$g^* \Nc^\vee_{D_1 / X} \to g^*\Lb_{D_1/D_1+D_2}[-1] \to \Nc^\vee_{D_1 \times_X D_2 / D_2}$$
is an isomorphism, so the latter morphism must be a surjection. Note that the shifted cotangent complex $\Lb_{D_1/D_1+D_2}[-1]$ is connective because $D_1$ is a closed derived subscheme of $D_1 + D_2$.
\end{proof}

We are now ready to prove the following important result.

\begin{lem}\label{lem:strictransvspullback}
Let $Z \hook D \hook X$ be a sequence of derived regular embeddings, where $D$ is a virtual Cartier divisor on $X$. Then, denoting $\pi: \bl_Z(X) \to X$, we have an equivalence 
$$\pi^*(D) \simeq \Ec + \bl_Z(D)$$ 
of virtual Cartier divisors on $X$, where $\pi^*(D)$ denotes the pullback of $D \hook X$ to the blowup.
\end{lem}
\begin{proof}
One checks locally and by reducing to the classical case that, for any quasi-coherent sheaf $\Fc$ on $X$, the unit map $\Fc \to \pi_*(\pi^*(\Fc))$ is an equivalence \cite{khan:dcdh}. Hence, recalling that $\Ec$ is the projective bundle $\Pb(\Nc_{Z/X})$ over $Z$, we deduce that the cofiber sequence
$$\Oc(D-\Ec) \xto{s_\Ec} \Oc(D) \to \Oc_\Ec(D)$$
pushes forward to
$$\Ic_{Z/X}(D) \to \Oc(D) \to \Oc_Z(D),$$
where $\Ic_{Z/X}$ is the fiber of $\Oc_X \to \Oc_Z$. The sequence $Z \hook D \hook X$ provides us with a commutative diagram
$$
\begin{tikzcd}
\Oc_X \arrow[d]{}{s_D} \arrow[rd]{}{0} & \\
\Oc_X(D) \arrow[r] & \Oc_Z(D)
\end{tikzcd}
$$
which by the above analysis corresponds to a global section of $s_{D'}$ of $\Oc(D-\Ec)$, where $D'$ is, by definition, the derived vanishing locus of $s_{D'}$.

As $\pi^*(D) = D' + \Ec$, it is enough for us to identify $D'$ with the strict transform $\bl_Z(D)$. Denoting by $\Ec'$ the intersection of $\Ec$ and $D'$ inside $\bl_Z(X)$, we may form the commutative diagram
$$
\begin{tikzcd}
\Ec' \arrow[d,hook] \arrow[r,hook] & D' \arrow[d,hook]  \\
\Ec \arrow[d]{}{g'} \arrow[r,hook] & \pi^*(D) \arrow[d] \arrow[r,hook] & \bl_Z(X) \arrow[d] \\
Z \arrow[r,hook] & D \arrow[r,hook] & X,
\end{tikzcd}
$$
where the commutativity of the square in the bottom-left corner follows from the general properties of fiber squares. We claim that the square 
$$
\begin{tikzcd}
\Ec' \arrow[d]{}{g} \arrow[r] & D' \arrow[d,hook]  \\
Z \arrow[r,hook] & D 
\end{tikzcd}
$$
exhibits a virtual Cartier divisor on $D'$ lying over $Z$. Indeed the square truncates into a classical Cartesian square as a composition of two squares with this property, and the canonical morphism $g^*\Nc^\vee_{Z/D} \to \Nc^\vee_{\Ec'/D'}$ is a surjection, because $g'^* \Nc^\vee_{Z/D} \to \Nc^\vee_{\Ec / \pi^*(D)}$ is (and because of the second claim of Lemma \ref{lem:vdivsumint}).

We have obtained a canonical morphism $D' \to \bl_Z(D)$. We claim that this is an equivalence. Using naturality in pullbacks, it is enough to verify the claim in the cases where $X = \Ab^n = \Spec (\Zb[x_1,...,x_n])$, $D$ is the derived vanishing locus of $x_1$, and $Z = \{0\}$ is the derived vanishing locus of $x_1,...,x_n$. But this is classical, so we are done.
\end{proof}

Next, we study derived blowups along derived vanishing loci of global sections of vector bundles, which admit an alternative universal property as the universal derived scheme on which the section can be replaced by a nowhere-vanishing one.  

Suppose that $X$ is a derived scheme, $E$ is a vector bundle on $X$, and $s$ a global section of $E$. Let us denote by $s'$ the global section of $Q$ on $\Pb(E)$ corresponding to $s$ under the equivalence of Lemma \ref{lem:pbsect}. On the derived vanishing locus $Z := V_{\Pb(E)}(s')$, there exists a canonical null-homotopy of $s'$, and therefore the section $s$ lifts canonically to a global section $\tilde s$ of $\Oc(-1) \vert_{Z}$. As $V_Z(\tilde s)$ is naturally identified as $V_{\Pb(E)}(s)$, we may form the commutative square
$$
\begin{tikzcd}
V_Z(\tilde s) \arrow[d]{}{g} \arrow[r]{}{i} & Z \arrow[d] \\
V_X(s) \arrow[r] & X,
\end{tikzcd}
$$
where $g$ is the restriction of $\Pb(E) \to X$ over $V_X(s)$. 

\begin{lem}\label{lem:dblowuppres}
The above square exhibits a virtual Cartier divisor on $Z$ lying over $V_X(s)$. Moreover, the induced map $Z \to \bl_{V_X(s)}(X)$ is an isomorphism.
\end{lem}
\begin{proof}
The outer square of the commutative diagram
$$
\begin{tikzcd}
V_Z(\tilde s) \arrow[hook,r]{}{i} \arrow[d]{}{g} & Z \arrow[hook, r] \arrow[rd] & \Pb(E) \arrow[d] \\
V_X(s) \arrow[rr] && X
\end{tikzcd}
$$
is Cartesian, and therefore one concludes from the cofiber sequence of cotangent complexes associated to the upper row that the induced map $g^*: \Nc^\vee_{V_X(s) / X} \to \Nc^\vee_{V_Z(\tilde s) / Z}$ is surjective. Next we prove that the induced map $Z \to \bl_{V_X(s)}(X)$ is an equivalence. Since being an equivalence can be checked locally, and as both constructions are stable under pullbacks, it is enough to prove this only in the case when $X = \Ab^n$, $E = \Oc_{\Ab^n}^{\oplus n}$, and $s = (x_1,...,x_n)$. But this is a classical result, so the claim follows.
\end{proof}

Hence the derived blow ups of vector bundle section admit the following universal property as ``resolution schemes''\footnote{A \emph{resolution scheme} may be associated to an arbitrary morphism of vector bundles $E \to F$ by considering the derived vanishing locus of a natural map of vector bundles on a relative Grassmannian. However, when the ranks of $E$ and $F$ are more than $1$, these schemes do not admit a description in terms of derived blowup of a derived regular embedding.}.

\begin{thm}\label{thm:dblowupres}
Let $X$ be a derived scheme, $E$ a vector bundle on $X$, and $s$ a global section of $E$. Then, for each $X$-scheme $S$, the space of $X$-morphisms $S \to \bl_{V_X(s)}$ is equivalent to the space of commutative diagrams
$$
\begin{tikzcd}
\Oc_Y \arrow[d] \arrow[rd]{}{s \vert_Y} & \\
\Ls \arrow[r, hook]{}{} & E \vert_Y
\end{tikzcd}
$$
of vector bundles, where the horizontal map is an inclusion of a line bundle into the vector bundle $E$.
\end{thm}
\begin{proof}
Combining the universal properties of projective bundles and derived vanishing loci, we see that the space of $X$-morphisms $S \to \bl_{V_X(s)}$ is equivalent to the space of commutative diagrams
$$
\begin{tikzcd}
& \Oc_Y \arrow[d]{}{s \vert_Y} \arrow[rd]{}{0} & \\
\Ls \arrow[r,hook] & E \vert_Y \arrow[r] & Q,
\end{tikzcd}
$$
where the bottom row is a cofiber sequence. However, such diagrams are equivalent to diagrams in the statement as cofiber sequences of quasi-coherent sheaves are fiber sequences.
\end{proof}

\chapter{Algebraic cobordism and its basic properties}\label{ch:cob}

In this Chapter, we construct the bivariant algebraic cobordism $\Omega$, following the work done in \cite{annala-cob, annala-yokura, annala-chern, annala-base-ind-cob}. The basic idea of the construction is rather easy: what one needs to do is to combine the construction of derived algebraic bordism of Lowrey--Schürg \cite{lowrey--schurg} with Yokura's universal bivariant theory in order to obtain a bivariant theory, the associated homology theory of which recovers derived algebraic bordism. However, instead of using the original Lowrey--Schürg relations to define our theory, we use slightly more geometric set of relations, closely related to the double point relations of Levine--Pandharipande \cite{levine-pandharipande}.

\section{Universal precobordism}\label{sect:precobcons}
\subsection{The bivariant functoriality}\label{ssect:cobbivfunc}

Here, we record the bivariant functorialities that will be relevant in this thesis. The most general bivariant functoriality we are going to work with is the following.

\begin{defn}\label{def:genfun}
We denote by $\Fc_a := (\Cc_a, \Cs_a, \Is_a, \Ss_a)$ the bivariant functoriality in which
\begin{enumerate}
\item $\Cc_a$ is the full sub-$\infty$-category of derived schemes consisting of qcqs derived schemes\footnote{In \cite{annala-base-ind-cob} we restricted our attention to the full subcategory of those finite-Krull-dimensional Noetherian derived schemes that admit an ample line bundle (or an ample family of line bundles). However, this is problematic as the Noetherian property is not stable under fiber products! Hence, for formal reasons, we have to define $\Cc_a$ as above, even though we only care about the bivariant groups of morphisms of pleasant enough Noetherian schemes.};

\item $\Cs_a$ consists of projective morphisms;

\item $\Is_a$ consists of all Cartesian squares;

\item $\Ss_a$ consists of all formally quasi-smooth morphisms.
\end{enumerate}
\end{defn}


At times it will be useful to restrict our degree of generality.

\begin{defn}\label{def:Afun}
Let $S$ be a finite-Krull-dimensional Noetherian derived scheme that admits an ample line bundle. Then, $\Fc_S := (\Cc_S, \Cs_S, \Is_S, \Ss_S)$ is the bivariant functoriality in which
\begin{enumerate}
\item $\Cc_S$ is the full sub-$\infty$-category of derived $S$-schemes consisting of those objects that are quasi-projective over $S$;

\item $\Cs_S$ consists of projective $S$-morphisms;

\item $\Is_S$ consists of all Cartesian squares;

\item $\Ss_S$ consists of all quasi-smooth $S$-morphisms\footnote{As all morphisms in $\Cc_S$ are locally of almost finite presentation, every formally quasi-smooth morphism in $\Cc_S$ is quasi-smooth.}.
\end{enumerate}
If $S = \Spec(A)$, then, then the notation $\Fc_A$ is also used.
\end{defn}

Note that any nice enough oriented bivariant theory with a functoriality as above has strong orientations along smooth morphisms, as explained below.

\begin{lem}\label{lem:autopoincare}
The class of smooth morphisms forms a class of very specialized morphisms for $\Fc_S$. In particular, if $\Bb$ is a centrally and stably oriented bivariant theory with functoriality $\Fc_S$, then 
\begin{enumerate}
\item the orientations along smooth morphisms are strong;

\item $\Bb$ satisfies Poincaré duality for smooth objects (Proposition \ref{prop:poincare}).
\end{enumerate}
\end{lem}
\begin{proof}
The results follow immediately from the basic properties of smooth morphisms in derived algebraic geometry, and from the results of Section \ref{ssect:poincare}.
\end{proof}

We end this subsection by the following remark concerning partial grading.

\begin{rem}\label{rem:bivgrad}
The universal additive theory $\Ab_{\Fc_a}$ admits a partial grading (Definition \ref{def:bivpartgrad}) along morphisms that are locally of almost finite presentation. Indeed, if $X \to Y$ is such a map, then an element $[V \to X] \in \Ab_{\Fc_a}(X \to Y)$, where $V$ is connected, the morphism $V \to X$ is projective, and the composition $V \to Y$ is (formally) quasi-smooth, has degree $d$, where $-d$ is the relative virtual dimension of $V$ over $Y$. The corresponding graded bivariant group shall be denoted by $\Ab^*_{\Fc_a}(X \to Y)$. Moreover, the universal additive theory $\Ab_{\Fc_S}$ admits a grading, and the theory shall be denoted by $\Ab^*_{\Fc_S}$. The sets of relations imposed on these theories in this section in order to obtain bivariant algebraic cobordism are \emph{homogeneous} in the sense that algebraic cobordism inherits the above (partial) grading. 
\end{rem}

\subsection{Universal precobordism}

Here, we construct and study the universal precobordism $\PCob$, which is an intermediate step in the construction of bivariant algebraic cobordism. It is obtained from the universal additive theory with the desired functoriality, by enforcing the derived analogue of Levine--Pandharipande's double point relations. Importantly, we will show that the Euler classes of line bundles satisfy a formal group law, which provides $\PCob$ with an $\Lb$-linear structure, where $\Lb$ is the Lazard ring \cite{lazard:1955}.

\begin{defn}\label{def:univprecob}
The \emph{universal precobordism} $\PCob$ is defined as the quotient of $\Ab_{\Fc_a}$ by the \emph{bivariant ideal $\Ic_\mathrm{dpt}$} of \emph{derived double point relations}, where $\Ic_\mathrm{dpt}(X \to Y)$ is the subgroup of $\Ab_{\Fc_a}(X \to Y)$ generated by elements of form
$$[W_0 \to X] - [D_1 \to X] - [D_2 \to X] + [\Pb_{D_1 \cap D_2}(\Oc(D_1) \oplus \Oc) \to X],$$
associated to a projective morphism $W \to \Pb^1 \times X$ such that the composition $W \to \Pb^1 \times Y$ is formally quasi-smooth, where
\begin{enumerate}
\item $W_0$ is the fiber of $W \to \Pb^1 \times X$ over $\{0\} \times X$;

\item the fiber $W_\infty$ of $W \to \Pb^1 \times X$ over $\{\infty\} \times X$ is the sum of virtual Cartier divisors $D_1$ and $D_2$ on $W$;

\item $D_1 \cap D_2$ is the fiber product of $D_i$ over $W$.
\end{enumerate}
It is straightforward to check that $\Ic_\mathrm{dpt}$ is indeed a bivariant ideal. Similarly, but starting with $\Ab^*_{\Fc_S}$, one constructs the \emph{universal $S$-precobordism} $\PCob^*_S$. 
\end{defn}

\begin{rem}\label{rem:precobvsAprecob}
Clearly, there are natural isomorphisms
$$\PCob^*_S(X \xto{\tilde f} Y) \cong \PCob^*_S(X \xto{f} Y)$$
compatible with all the bivariant operations and orientations, where $f$ is obtained from the $S$-morphism $\tilde f$ by forgetting the structure morphisms and the homotopy witnessing the commutativity of the triangle
$$
\begin{tikzcd}
X \arrow[d]{}{\pi_X} \arrow[r]{}{f} & Y \arrow[ld]{}{\pi_Y} \\
S.
\end{tikzcd}
$$
From now on, when we prove that a result holds for $\PCob$, then, in order to avoid unnecessary repetition, we will not explicitly mention that the analogous result holds for $\PCob^*_S$ as well. The same applies for the precobordism groups with bundles constructed below.
\end{rem}

Following the classical paper of Lee and Pandharipande \cite{pandharipande:2012}, we will construct ``algebras'' over $\PCob$, in which the cycles are equipped with line or vector bundles on the source. These theories will be useful for algebraic cobordism computations, as we shall see in the sequel. Let us begin with the following modification of the universal additive bivariant theory.

\begin{defn}\label{def:univaddthywbundles}
Let $\Fc$ be either $\Fc_a$ or $\Fc_S$. Then we define the \emph{additive bivariant theories with bundles} $(\Ab_\Fc^{\bullet,*}, \bullet_\oplus)$ and $(\Ab_\Fc^{\bullet,*}, \bullet_\otimes)$ as the bivariant theory where $\Ab_\Fc^{\bullet,*}(X \to Y)$ is the Abelian group generated by equivalence classes $[V \to X; E]$, where $V \to X$ is a projective morphism such that the composition $V \to Y$ is (formally) quasi-smooth, and where $E$ is a vector bundle on $X$, modulo relations of form
$$[V_1 \sqcup V_2 \to X; E_1 \sqcup E_2] = [V_1 \to X; E_1] + [V_2 \to X; E_2].$$
If $E$ has rank $r$, then $[V \to X; E] \in \Ab^{\bullet, r}_\Fc(X \to Y)$. The bivariant operations are defined as follows:
\begin{enumerate}
\item \emph{bivariant pushforward}: if $X \to  Y$ factors through a projective map $g: X \to X'$, then the bivariant pushforward is defined by linearly extending the formula
$$g_*([V \xto{f} X; E]) := [V \xto{g \circ f} X'; E];$$

\item \emph{bivariant pullback}: if
$$
\begin{tikzcd}
X' \arrow[]{d}{g'} \arrow[]{r} & Y' \arrow[]{d}{g} \\
X \arrow[]{r} & Y
\end{tikzcd}
$$
is Cartesian, then the bivariant pullback is defined by linearly extending the formula
$$g^*([V \xto f X; E]) := [V' \xto{f'} X'; g''^*(E)],$$
where $f'$ is the pullback of $f$ along $g'$ and $g''$ is the natural map $V' \to V$;

\item \emph{bivariant products}: the two bivariant products are defined by bilinearly extending the formulas
$$[V \to X; E] \bullet_\oplus [W \to X; F] := [V' \to X; E \vert_{V'} \oplus F \vert_{V'}]$$
and
$$[V \to X; E] \bullet_\otimes [W \to X; F] := [V' \to X; E \vert_{V'} \otimes F \vert_{V'}],$$
where the morphism $V' \to X$ is as in the diagram
$$
\begin{tikzcd}
V' \arrow[]{r} \arrow[]{d} & X' \arrow[]{r} \arrow[]{d} & W \arrow[]{d} \\
V \arrow[]{r} & X \arrow[]{r} & Y \arrow[]{r} & Z;
\end{tikzcd}
$$
\end{enumerate}
One checks that these operations are well defined similarly to the proof of Lemma \ref{lem:univthyisbiv}. Note that the second degree is additive in $\bullet_\oplus$ but multiplicative in $\bullet_\otimes$. Moreover, we define the \emph{additive theory with line bundles} as the subtheory $(\Ab^{\bullet,1}_\Fc, \bullet_\otimes)$ of $(\Ab^{\bullet,*}_\Fc, \bullet_\otimes)$.

These theories are stably oriented by the elements 
$$\theta_\oplus(f) := [X \to X; 0] \in \Ab^{\bullet, 0}_\Fc(X \to Y)$$
and
$$\theta_\otimes(f) := [X \to X; \Oc_X] \in \Ab^{\bullet, 1}_\Fc(X \to Y)$$
respectively, where $f: X \to Y$ is a formally quasi-smooth morphism. Moreover, all three theories are additive in the sense of Definition \ref{def:add}.
\end{defn}

The universal property of $\Ab_\Fc$ induces a Grothendieck transformations $\Ab_\Fc \to (\Ab^{\bullet, *}_\Fc, \bullet_\oplus)$ and $\Ab_\Fc \to (\Ab^{\bullet, *}_\Fc, \bullet_\otimes)$ which are defined by the formulas 
$$[V \to X] \mapsto [V \to X; 0]$$
and 
$$[V \to X] \mapsto [V \to X; \Oc_V],$$
respectively. Moreover, there are also the \emph{forgetful Grothendieck transformations} $(\Ab^{\bullet, *}_\Fc, \bullet_\oplus) \to \Ab_\Fc$ and $(\Ab^{\bullet, *}_\Fc, \bullet_\otimes) \to \Ab_\Fc$, which forget the vector bundle.

\begin{defn}\label{def:precobwbundles}
The \emph{precobordism theories with bundles} $(\PCob^{\bullet, *}, \bullet_\oplus)$ and $(\PCob^{\bullet, *}, \bullet_\otimes)$ are defined as the quotients of $(\Ab^{\bullet, *}_{\Fc_a}, \bullet_\oplus)$ and $(\Ab^{\bullet, *}_{\Fc_a}, \bullet_\oplus)$, respectively, by the bivariant ideal $\Ic^{\bullet,*}_\mathrm{dpt}$ of \emph{derived double point relations with bundles}, where $\Ic^{\bullet,*}_\mathrm{dpt}(X \to Y)$ is the subgroup of $\Ab^{\bullet,*}_{\Fc_a}(X \to Y)$ generated by elements of form
$$[W_0 \to X; E \vert_{W_0}] - [D_1 \to X; E \vert_{D_1}] - [D_2 \to X; E \vert_{D_2}] + [\Pb \to X; E \vert_{\Pb}],$$
where $\Pb = \Pb_{D_1 \cap D_2}(\Oc(D_1) \oplus \Oc)$, associated to a projective morphism $W \to \Pb^1 \times X$ such that the composition $W \to \Pb^1 \times Y$ is formally quasi-smooth, and
\begin{enumerate}
\item $E$ is a vector bundle on $W$;

\item $W_0$ is the fiber of $W \to \Pb^1 \times X$ over $\{0\} \times X$;

\item the fiber $W_\infty$ of $W \to \Pb^1 \times X$ over $\{\infty\} \times X$ is the sum of virtual Cartier divisors $D_1$ and $D_2$ on $W$;

\item $D_1 \cap D_2$ is the fiber product of $D_i$ over $W$.
\end{enumerate}
It is straightforward to check that $\Ic^{\bullet,*}_\mathrm{dpt}$ is indeed a bivariant ideal. Similarly, but starting with $\Ab^*_{\Fc_S}$, one constructs the \emph{$S$-precobordism theories with bundles} $(\PCob^{*,*}_S, \bullet_\oplus)$ and $(\PCob^{*,*}_S, \bullet_\otimes)$. 
\end{defn}

We will denote the induced Grothendieck transformations $\PCob \to (\PCob^{\bullet,*}, \bullet_\oplus)$ and $\PCob \to (\PCob^{\bullet,*}, \bullet_\otimes)$ by $\iota_\oplus$ and $\iota_\otimes$, respectively. The same notation is used in the case of $S$-precobordism groups as well.

\subsection{Euler classes}

Euler classes of vector bundles will play a fundamental role in our study of algebraic  cobordism.

\begin{defn}\label{def:eulerclass}
Let $\Bb$ be an oriented bivariant theory with functoriality $\Fc_a$ or $\Fc_S$. Then, if $X \in \Cc$ and $E$ is a vector bundle on $X$, the \emph{Euler class} of $X$ is defined as 
$$e(E) := i_{0!}(1_{V_X(0)}) \in \Bb^\bullet(X),$$
where $i_0: V_X(0) \hook X$ is the derived vanishing locus of the zero-section of $E$.
\end{defn}

Clearly, in case the vector bundle $E$ has rank $r$, $e(E) \in \PCob^r(X)$. The following result will be useful in practice.

\begin{lem}\label{lem:eulerclasssect}
Let $X$ be a derived scheme, $E$ a vector bundle on $X$, and let $s$ be a global section of $E$. Then
$$e(E) = [V_X(s) \to X] \in \PCob^*(X).$$
\end{lem}
\begin{proof}
By definition, $e(E) = [V_X(0) \to X] \in \PCob^*(X)$. Let $x_0 \otimes s$ be a section of $E(1)$ on $\Pb^1 \times X$. Then, the map $V_{\Pb^1 \times X}(x_0 \otimes s) \to \Pb^1 \times X$ realizes the desired equation in $\PCob^*(X)$\footnote{Here, we take $D_2 = \emptyset$.}.
\end{proof}

In particular, the Euler class of a vector bundle admitting a nowhere-vanishing section (e.g., a trivial line bundle), is trivial. Euler classes are multiplicative in cofiber sequences of vector bundles.

\begin{prop}\label{prop:eulerclassmult}
Let $X$ be a derived scheme. Then, if $E' \to E \to E''$ is a cofiber sequence of vector bundles on $X$, we have that $e(E) = e(E') \bullet e(E'') \in \PCob^*(X)$.
\end{prop}
\begin{proof}
Let $s$ be a global section of $E$, and let $s''$ be the induced global section of $E''$. Then, on $V_X(s'')$, the section $s$ admits a natural lift to a section $s'$ of $E'$, and moreover $V_{V_X(s'')}(s')$ is naturally equivalent to $V_X(s)$. Thus, denoting by $i$ the inclusion $V_X(s'') \hook X$, we have
\begin{align*}
e(E) &= i_!(e(E' \vert_{V_X(s'')})) \\
&= e(E') \bullet i_!(1_{V_X(s'')}) & (\text{projection formula}) \\
&= e(E') \bullet e(E''),
\end{align*}
as desired.
\end{proof}

\subsection{The extended double point relation}

There exists an analogue of Levine--Pandharipande's extended double point relations (see \cite{levine-pandharipande} Lemma 16) for $\PCob$. Before proving this, we record the derived blowup formula, which will be useful in cobordism computations.

\begin{lem}\label{lem:dblowupform}
Let $Z \hook X$ be a derived regular embedding. Then,
$$1_X = [\bl_Z(X) \to X] + [\Pb_Z(\Nc_{Z / X} \oplus \Oc) \to X] - [\Pb_\Ec(\Oc(-1) \oplus \Oc) \to X] \in \PCob^*(X),$$
where $\Ec \hook \bl_Z(X)$ is the exceptional divisor. 
\end{lem}
\begin{proof}
Consider $W := \bl_{\infty \times Z}(\Pb^1 \times X)$. Then, by Lemma \ref{lem:strictransvspullback}, the fiber of $W \to \Pb^1$ over $\infty$ is the sum of $\bl_Z(X)$ and the exceptional divisor $\Ec' \hook W$, leading to the above formula. As $\Oc(\Ec')\vert_{\Ec'} \simeq \Oc(-1)$, the third term of the right hand side is of the correct form.
\end{proof}

We are now ready to prove the extended double point formula. Let us start by fixing notation. Let $X$ be a derived scheme and let $\Ls_1$ and $\Ls_2$ be line bundles on $E$. For simplicity, let us denote $\Ls_3 := \Ls_1 \otimes \Ls_2$. Then, 
\begin{align*}
\Pb_1 &:= \Pb_X(\Ls^\vee_2 \oplus \Oc); \\
\Pb_2 &:= \Pb_{\Pb(\Ls_1 \oplus \Ls_3)}(\Oc(-1) \oplus \Oc); \\
\Pb_3 &:= \Pb(\Ls_1 \oplus \Ls_3 \oplus \Oc).
\end{align*} 

\begin{lem}\label{lem:edpt}
Let everything be as above. Then
\begin{align*}
e(\Ls_3) &= e(\Ls_1) + e(\Ls_2) - e(\Ls_1) \bullet e(\Ls_2) \bullet [\Pb_1 \to X] \\
&+ e(\Ls_1) \bullet e(\Ls_2) \bullet e(\Ls_3) \bullet ([\Pb_2 \to X] - [\Pb_3 \to X])
\end{align*}
in $\PCob^*(X)$.
\end{lem}
\begin{proof}
Let $D_1, D_2$ and $D_3$ be the derived vanishing loci of arbitrary sections of $\Ls_1, \Ls_2$ and $\Ls_3$, respectively. Let $\pi:X' \to X$ be the derived blowup of $X$ at $D_1 \cap D_3$, let $D'_1$ and $D'_3$ be the strict transforms of $D_1$ and $D_3$, repectively, and let $D'_2$ be the pullback of $D_2$ to $X'$. Note that $D'_1$ and $D'_3$ do not intersect, and that the virtual Cartier divisors $D'_3$ and $D'_1 + D'_2$ are rationally equivalent.

Next, denote by $X''$ the blowup of $X'$ at $D'_2 \cap D'_3$, by $D''_2$ and $D''_3$ the strict transforms of $D'_2$ and $D'_3$, respectively, and by $D''_1$ the pullback of $D'_1$ to $X''$. As $D'_1$ and $D'_3$ do not intersect, the natural map $D''_1 \to D'_1$ is an equivalence. As the virtual Cartier divisors $D''_3$ and $D''_1 + D''_2$ are rationally equivalent and disjoint, they give a rise to a morphism $X'' \to \Pb^1$ with $D''_3$ the fiber over $0$ and $D''_1 + D''_2$ the fiber over $\infty$. The induced derived double point relation is 
$$[D''_3 \to X] = [D''_1 \to X] + [D''_2 \to X] - [\Pb_{D''_1 \cap D''_2}(\Oc(D''_1) \oplus \Oc) \to X] \in \PCob(X).$$
We will investigate each term of the equation. 
\begin{enumerate}
\item $[D''_3 \to X]$: As $D''_3$ is the composition of two derived blowups along virtual Cartier divisors, $[D''_3 \to X] = [D_3 \to X] = e(\Ls_3)$ 

\item $[D''_1 \to X]$: As $D''_1$ is the derived blowup of $D_1$ along a virtual Cartier divisor, we have that $[D''_1 \to X] = [D_1 \to X] = e(\Ls_1)$.

\item $[D''_2 \to X]$: Arguing as above, $D''_2$ is equivalent to the derived blowup $\bl_{D_1 \cap D_2 \cap D_3}(D_2)$. Hence, applying the derived blowup formula, we compute that
\begin{align*}
[D''_2 \to X] &= [D_2 \to X] - [\Pb_{D_1 \cap D_2 \cap D_3}(\Oc(D_1) \oplus \Oc(D_3) \oplus \Oc) \to X] \\
&+ [\Pb_{\Pb_{D_1 \cap D_2 \cap D_3}(\Oc(D_1) \oplus \Oc(D_3))}(\Oc(-1) \oplus \Oc)] \\
&= e(\Ls_2) + e(\Ls_1) \bullet e(\Ls_2) \bullet e(\Ls_3) \bullet ([\Pb_2 \to X] - [\Pb_3 \to X]).
\end{align*}

\item $[\Pb_{D''_1 \cap D''_2}(\Oc(D''_1) \oplus \Oc) \to X]$: As $D'_1$ and $D'_3$ do not intersect, the natural map $D''_1 \cap D''_2 \to D'_1 \cap D'_2$ is an equivalence. On the other hand, as
\begin{align*}
D'_2 \times_{X'} D_1 &\simeq (D_2 \times_X X') \times_{X'} D_1 \\
&\simeq D_2 \times_{X} D_1,
\end{align*}
$D'_1 \cap D'_2 \simeq D_1 \cap D_2$. As
\begin{align*}
\Oc(D''_1) \vert_{D''_1} &\simeq \Oc(D'_1) \vert_{D'_1} \\
&= \Oc(\pi^*(D_1) - \Ec) \vert_{D'_1} \\
&= \Oc(D_1 - D_3) \vert_{D_1},
\end{align*}
where $\Ec$ is the exceptional divisor of $X'$, we have that
\begin{align*}
[\Pb_{D''_1 \cap D''_2}(\Oc(D''_1) \oplus \Oc) \to X] &= [\Pb_{D_1 \cap D_2}(\Oc(D_1 - D_3) \oplus \Oc) \to X] \\
&= e(\Ls_1) \bullet e(\Ls_2) \bullet [\Pb_1 \to X].
\end{align*}
\end{enumerate}
Combining the above equations yields the desired formula.
\end{proof}

As the first application of the above result, we prove the nilpotence of Euler classes.

\begin{prop}\label{prop:eulerclassnilp}
Let $X$ be a finite-Krull-dimensional Noetherian scheme that admits an ample line bundle. Then, for each vector bundle $E$ on $X$, the Euler class $e(E)$ is nilpotent in $\PCob(X)$. 
\end{prop}
\begin{proof}
If $E$ admits a sequence of sections $s_1,...,s_r$ such that they simultaneously vanish nowhere, then $e(E)^r = 0$. In particular, since $X$ is quasi-compact, such a finite sequence can be found for $E$ that is globally generated.

Suppose then that $\Ls$ is a line bundle the dual of which is globally generated. Then, it follows from Lemma \ref{lem:edpt} by setting $\Ls_1 = \Ls$ and $\Ls_2 = \Ls^\vee$, that
$$0 = e(\Ls) + e(\Ls^\vee) - e(\Ls) \bullet e(\Ls^\vee) \bullet [\Pb_1 \to X],$$
or, in other words, that
$$e(\Ls) = -\frac{e(\Ls^\vee)}{1 - e(\Ls^\vee) \bullet [\Pb_1 \to X]},$$
where the right hand side is well defined because $e(\Ls^\vee)$ is nilpotent. From this formula, it follows that also $e(\Ls)$ is nilpotent.

As $X$ admits an ample line bundle, any line bundle $\Ls$ on $X$ is equivalent to one of form $\Ls_1 \otimes \Ls_2^\vee$, where $\Ls_1$ and $\Ls_2$ are globally generated. Applying Lemma \ref{lem:edpt} again, we compute that 
$$e(\Ls) = \frac{e(\Ls_1) + e(\Ls_2^\vee) - e(\Ls_1) \bullet e(\Ls_2^\vee) \bullet [\Pb^1 \to X]}{1 - e(\Ls_1) \bullet e(\Ls_2^\vee) \bullet ([\Pb_2 \to X]-[\Pb_3 \to X])},$$
and the nilpotence of $e(\Ls)$ follows from that of $e(\Ls_1)$ and $e(\Ls_2^\vee)$.

Finally, let $E$ be a vector bundle on $X$. As $X$ admits an ample line bundle, there exists a line bundle $\Ls$ such that $E^\vee \otimes \Ls$ is globally generated. Let $s_1, ..., s_r$ be a sequence of morphisms $E \to \Ls$ that correspond to a sequence of global sections of $E^\vee \otimes \Ls$ the simultaneous vanishing locus of which is empty. In other words, they give rise to a surjective map $E^{\oplus r} \to \Ls$ of vector bundles, the kernel of which we denote by $F$. As Euler classes are multiplicative, we have that
$$e(E)^r = e(F) \bullet e(\Ls),$$
and the nilpotence of $e(E)$ follows from that of $e(\Ls)$. 
\end{proof}

\subsection{The formal group law of universal precobordism}

In order to prove that the Euler class of a tensor product of line bundles can be computed by means of a formal group law, we will need to be able to compute classes of certain projective bundles in the precobordism ring. To do so, we will compute certain elements in the precobordism theories with bundles.

It will be useful to fix the following notation. 

\begin{defn}\label{def:fundtower}
Given a derived scheme $X$ and a line bundle $\Ls$ on $X$, we recursively define the following pairs:
$$\big(P_0(X,\Ls), \Ms_0(X,\Ls) \big) := \big(X, \Ls \big)$$
and 
$$\big(P_{i+1}(X,\Ls), \Ms_{i+1}(X,\Ls) \big) := \big(\Pb_{P_i(X,\Ls)} (\Ms_i(X,\Ls) \oplus \Oc), \Ms_i(X,\Ls)(H_{i+1}) \big),$$
where $\Oc(H_{i+1})$ is the hyperplane bundle on $\Pb_{P_i(X,\Ls)} \big(\Ms_i(X,\Ls) \oplus \Oc\big)$. If $X = \Spec(\Zb)$ and $\Ls = \Oc$, then these pairs are denoted more simply by $(P_i, \Ms_i)$. Moreover, we define elements
$$\beta_i := [P_i \to \Spec(\Zb), \Ms_i]$$
in $\PCob^{i,1}(\Spec(\Zb))$, with the convention that $\beta_i = 0$ if $i$ is negative.
\end{defn}

\begin{lem}\label{lem:lbundform}
Let $X$ be a finite-Krull-dimensional Noetherian derived scheme that admits an ample line bundle, and let $\Ls$ be a line bundle on $X$. Then, the equation
$$[X \to X; \Ls] = \sum_{i \geq 0} e(\Ls)^i \bullet_\otimes \big(\beta_i - \beta_{i-1}\bullet_\otimes [\Pb(\Ls \oplus \Oc) \to X; \Oc]\big)$$
holds in $\PCob^{*,1}(X)$.  
\end{lem}
\begin{proof}
Consider the derived blowup $W$ of $\Pb^1 \times X$ at $\infty \times D$, where $D$ is the derived vanishing locus of a section of $\Ls$. Moreover, let $\Ls'$ be the line bundle $\Ls(-\Ec)$ on $W$, where $\Ec$ is the exceptional divisor of $W$. The line bundle $\Ls'$ restricts to
\begin{enumerate}
\item $\Ls$ on the strict transform of $0 \times X$;
\item $\Oc_X$ on the strict transform of $\infty \times X$;
\item $\Ls(1)$ on the exceptional divisor $\Ec \simeq \Pb_D(\Ls \oplus \Oc)$;
\item $\Oc_X$ on the intersection of the strict transform of $\infty \times X$ and $\Ec$.
\end{enumerate}
Hence, the derived scheme $P_i(W, \Ls')$ together with its natural map to $\Pb^1 \times X$ and the line bundle $\Ms_i(W, \Ls')$, realizes the equation
\begin{align*}
&[P_i(X, \Ls) \to X; \Ms_i(X, \Ls)] \\
&= [P_i(X, \Oc) \to X; \Ms_i(X, \Oc)]\\
&+ [P_{i+1}(D, \Ls \vert_D) \to X; \Ms_{i+1}(D, \Ls \vert_D)]\\
&- [P_i(\Pb_D(\Ls\vert_D \oplus \Oc), \Oc) \to X; \Ms_i(\Pb_D(\Ls\vert_D \oplus \Oc), \Oc) ] \\
&= \beta_i + e(\Ls) \bullet_\otimes [P_{i+1}(X, \Ls) \to X; \Ms_{i+1}(X, \Ls)] \\
&- \beta_{i} \bullet_\otimes e(\Ls) \bullet_\otimes [\Pb(\Ls \oplus \Oc) \to X; \Oc].
\end{align*}
Since $e(\Ls)$ is nilpotent, the claim follows from the above formulas.
\end{proof}

In order to use the above results in order to compute the classes of projective bundles, we need to be able to construct precobordism elements given a precobordism element equipped with a bundle.

\begin{defn}\label{def:projtrans}
The \emph{projectivization transformations} $\Pb,\wtil\Pb:\PCob^{\bullet,*}(X \to Y) \to \PCob(X \to Y)$ are defined by the formulas
$$[V \to X; E] \mapsto [\Pb_V(E) \to X]$$
and
$$[V \to X; E] \mapsto [\Pb_{\Pb_V(E)}(\Oc(-1) \oplus \Oc) \to X]$$
respectively. Clearly these maps are well-defined group homomorphisms. One defines $\Pb, \wtil\Pb: \PCob_S^{*,*} \to \PCob^*_S$ in a similar fashion.
\end{defn}

These transformations are not Grothendieck transformations. However, in a sense, they are maps of $\PCob$-bimodules.

\begin{prop}\label{prop:projtransprops}
Let $\epsilon$ be one of the transformations $\Pb, \wtil \Pb$. Then, 
\begin{enumerate}
\item $\epsilon$ commutes with pushforwards;
\item $\epsilon$ commutes with pullbacks;
\item if $\alpha \in \PCob(X \to Y)$ and $\beta \in \PCob^{\bullet,*}(Y \to Z)$, then
$$\epsilon\big(\iota_\oplus(\alpha) \bullet_\oplus \beta \big) = \alpha \bullet \epsilon(\beta)$$
and
$$\epsilon\big(\iota_\otimes(\alpha) \bullet_\otimes \beta \big) = \alpha \bullet \epsilon(\beta);$$
\item if $\alpha \in \PCob^{\bullet,*}(X \to Y)$ and $\beta \in \PCob(Y \to Z)$, then
$$\epsilon\big(\alpha \bullet_\oplus \iota_\oplus(\beta) \big) = \epsilon(\alpha) \bullet \beta$$
and
$$\epsilon\big(\alpha \bullet_\otimes \iota_\otimes(\beta) \big) = \epsilon(\alpha) \bullet \beta.$$
\end{enumerate}
\end{prop}
\begin{proof}
The claim follows from the fact that forming projective bundles is compatible with pullbacks.
\end{proof}

As a first application, we prove the following formula.

\begin{lem}\label{lem:pbform1}
Let $X$ be a Noetherian derived scheme that admits an ample line bundle. Then, for each line bundle $\Ls$ on $X$, the equation
$$[\Pb(\Ls \oplus \Oc) \to X] = \frac{\sum_{i \geq 0} \beta'_{i+1} e(\Ls)^i}{\sum_{i \geq 0} \beta'_i e(\Ls)^i}$$
holds in $\PCob^*(X)$, where $\beta'_i := [P_i \to \Spec(\Zb)] \in \PCob^*(\Spec(\Zb))$, and where $P_i$ as in Definition \ref{def:fundtower}. 
\end{lem}
\begin{proof}
Rephrasing Lemma \ref{lem:lbundform}, we obtain
$$[X \to X; \Ls] = \sum_{i \geq 0} e(\Ls)^i \bullet_\oplus \big(\beta_i - \beta_{i-1} \bullet_\oplus \iota_\oplus([\Pb(\Ls \oplus \Oc) \to X])\big) \in \PCob^{*,1}(X).$$
Moreover, if we denote by $[\Oc]$ the class $[\Spec(\Zb) \to \Spec(\Zb); \Oc]$ in the ring $\PCob^{*,1}(\Spec(\Zb))$, then $\Pb(\beta_i \bullet_\oplus [\Oc]) = \beta'_{i+1}$ unless $i$ is negative, and hence, computing $\Pb([\Oc] \bullet_\oplus [X \to X; \Ls])$ using the two sides of the above equation, we obtain the formula 
$$[\Pb(\Ls \oplus \Oc) \to X] = \sum_{i \geq 0} e(\Ls)^i \beta'_{i+1} - \sum_{i \geq 1} e(\Ls)^i \bullet \beta'_{i} [\Pb(\Ls \oplus \Oc) \to X],$$
from which the desired equation follows by solving $[\Pb(\Ls \oplus \Oc) \to X]$.
\end{proof}

Secondly, we need to compute the classes $[\Pb_2 \to X]$ and $[\Pb_3 \to X]$ that appear in the extended double point formula (Lemma \ref{lem:edpt}). 

\begin{lem}\label{lem:pbform2}
There exists power series $G(x,y), H(x,y) \in \PCob^\bullet(\Spec(\Zb))[[x,y]]$ such that, if $X$ is a Noetherian derived scheme with an ample line bundle, then,
\begin{enumerate}
\item $[\Pb_{\Pb(\Ls_1 \oplus (\Ls_1 \otimes \Ls_2))}(\Oc(-1) \oplus \Oc) \to X] = G\big(e(\Ls_1), e(\Ls_2)\big)$;
\item $[\Pb(\Ls_1 \oplus (\Ls_1 \otimes \Ls_2) \oplus \Oc) \to X] = H\big(e(\Ls_1), e(\Ls_2)\big)$.
\end{enumerate}
\end{lem}
\begin{proof}
\begin{enumerate}
\item The claim follows by expressing
\begin{align*}
&[X \to X; \Ls_1 \oplus (\Ls_1 \otimes \Ls_2)] \\
&= [X \to X; \Ls_1] \bullet_\otimes  ([X \to X; \Oc] \bullet_\oplus [X \to X; \Ls_2])
\end{align*}
as a power series in $e(\Ls_1)$ and $e(\Ls_2)$ and with coefficients in $\PCob^{*,*}(\Spec(\Zb))$ which can be done by using Lemmas \ref{lem:lbundform} and \ref{lem:pbform1}, and then applying the transformation $\wtil \Pb$.

\item The claim follows by expressing
\begin{align*}
& [X \to X; \Ls_1 \oplus (\Ls_1 \otimes \Ls_2) \oplus \Oc] \\
=& [X \to X; \Ls_1] \bullet_\oplus  ([X \to X; \Ls_1] \bullet_\otimes [X \to X; \Ls_2]) \bullet \oplus [X \to X; \Oc]
\end{align*}
as a power series in $e(\Ls_1)$ and $e(\Ls_2)$ and with coefficients in $\PCob^{*,*}(\Spec(\Zb))$ which can be done by using Lemmas \ref{lem:lbundform} and \ref{lem:pbform1}, and then applying the transformation $\Pb$. \qedhere
\end{enumerate}
\end{proof}

We are now ready to prove the main result of this section.

\begin{thm}\label{thm:fgl}
There exists a power series $F(x,y) \in \PCob^*(\Spec(\Zb))[[x,y]]$ such that, given a Noetherian derived scheme $X$ admitting an ample line bundle, and line bundles $\Ls_1, \Ls_2$ on $X$, then
\begin{align*}
e(\Ls_1 \otimes \Ls_2) &= F\big(e(\Ls_1), e(\Ls_2)\big) \\
&=: \sum_{i,j \geq 0} a_{ij} e(\Ls_1)^i \bullet e(\Ls_2)^j \in \PCob^*(X).
\end{align*}
Moreover, the formal power series $F$ is a commutative formal group law. 
\end{thm}
\begin{proof}
As $\Pb(\Ls^\vee \oplus \Oc) \simeq \Pb(\Oc \oplus \Ls)$, we may combine Lemmas \ref{lem:edpt}, \ref{lem:pbform1} and \ref{lem:pbform2} to obtain the formula
\begin{align*}
&e(\Ls_1 \otimes \Ls_2) \\
&= e(\Ls_1) + e(\Ls_2) - e(\Ls_1) \bullet e(\Ls_2) \bullet \frac{\sum_{i \geq 0} \beta'_{i+1} e(\Ls_2)^i}{\sum_{i \geq 0} \beta'_i e(\Ls_2)^i} \\
&+ e(\Ls_1) \bullet e(\Ls_2) \bullet e(\Ls_1 \otimes \Ls_2) \bullet \Big(G\big(e(\Ls_1), e(\Ls_2)\big) - H\big(e(\Ls_1), e(\Ls_2)\big)\Big),
\end{align*}
from which the desired formula for $e(\Ls_1 \otimes \Ls_2)$ may be obtained. That $F(x,y)$ is a commutative formal group law, follows from the basic properties of tensor products of line bundles. 
\end{proof}

There exists a canonical linear structure on $\PCob$ over the Lazard ring.

\begin{cor}\label{cor:llinprecob}
Let $\Lb$ be the Lazard ring. Then, the ring homomorphism $\Lb \to \PCob^*(\Spec(\Zb))$  classifying the formal group law of Theorem \ref{thm:fgl} equips $\PCob^*$ with a canonical $\Lb$-linear structure. \qed
\end{cor}

\section{Algebraic cobordism}\label{sect:bivAcob}

Here, we concisely recall the definition of bivariant algebraic $S$-cobordism, where $S$ is a finite-Krull-dimensional Noetherian derived scheme that admits an ample line bundle. There exists also a base independent, and in a sense, simpler, version of the definition we consider below, which was studied in \cite{annala-base-ind-cob}, but we will not recall it here.

Suppose that we are provided with the data of a virtual Cartier divisor $D \hook W$, and an equivalence
$$D \simeq n_1 D_1 + \cdots + n_r D_r$$
of virtual Cartier divisors on $W$ with $n_i > 0$. Then, denoting by $+_{F}$ the formal addition given by the formal group law $F$ and by $[n]_{F} \cdot$ the formal multiplication (iterated formal addition), the formal power series 
$$[n_1]_{F} \cdot x_1 +_{F} \cdots +_{F} [n_r]_{F} \cdot x_r$$ 
in $r$ variables has a unique expression of form 
$$\sum_{I \subset \{1,2,...,r \}} \textbf{x}^I F^{n_1,...,n_r}_{I}(x_1,...,x_r),$$
where 
$$\textbf{x}^I = \prod_{i \in I} x_i$$
and $F_{I}^{n_1,...,n_r}(x_1,...,x_r)$ contains only variables $x_i$ such that $i \in I$. Note that $F_{\emptyset}^{n_1,...,n_r}(x_1,...,x_r) = 0$. Using this notation, we make the following definition.

\begin{defn}\label{def:zetaclass}
Let everything be as above, and suppose that $W$ is a quasi-smooth and quasi-projective derived $S$-scheme. Then, we define
$$
\zeta_{W,D,D_1,...,D_r} := \sum_{I \subset \{1,2,...,r\}} \iota^I_*\Biggl(F^{n_1,...,n_r}_{I} \Bigl( e\bigl(\Oc(D_1)\bigr), ...,e\bigl(\Oc(D_r)\bigr) \Bigr) \bullet 1_{D_I/pt} \Biggr)
$$
in $\PCob^S_*(D)$, where $\iota^I$ is the inclusion $D_I \hookrightarrow D$ of the derived intersection $\cap_{i \in I} D_i$ inside $W$.
\end{defn}

In the definition of bivariant algebraic $S$-cobordism, the following special case plays an important role.

\begin{defn}\label{def:Asncd}
Let $W$ be a smooth and quasi-projective (derived) $S$-scheme. Then an \emph{$S$-snc divisor} on $W$ is the data of an effective Cartier divisor $D \hookrightarrow W$, Zariski connected effective Cartier divisors $D_1,...,D_r \hookrightarrow W$ with 
$$D_I := \cap_{i \in I} D_i$$
smooth and of the expected relative dimension over $S$ for all $I \subset \{1,2,...,r\}$, such that there exist positive integers $n_1,...,n_r$ and an equivalence
$$D \simeq n_1 D_1 + \cdots + n_r D_r$$
as Cartier divisors on $W$. Somewhat abusively, we will use the shorthand notation
$$\zeta_{W,D} := \zeta_{W,D, D_1, ..., D_r}$$
whenever $D \hook W$ and $D_i$ form an $S$-snc divisor. 
\end{defn}

We are now ready to recall the definition of bivariant algebraic cobordism over $S$.

\begin{defn}\label{def:bivAcob}
\emph{Bivariant algebraic $S$-cobordism theory} is the additive stably oriented bivariant theory with functoriality $\Fs_S$ defined as the quotient
$$\Omega^*_S := \PCob^*_S / \langle \Rc^\mathrm{snc}_S \rangle,$$
where $\Rc^\mathrm{snc}_S$ is the bivariant subset of \emph{snc relations} consisting of 
$$\zeta_{W,D} - \theta(\pi_D) \in \PCob^*_S(D \to pt),$$ 
where $D \hook W$ ranges over all $S$-snc divisors, and $\langle \Rc^\mathrm{snc}_S \rangle$ is the bivariant ideal generated by $\Rc^\mathrm{snc}_S$.
\end{defn}

\begin{rem}
Suppose that $S$ is a finite-Krull-dimensional regular Noetherian scheme that admits an ample line bundle. Then, the algebraic $S$-bordism group $\Omega^S_*(X) := \Omega^{-*}_S(X \to S)$ is generated by cycles of form $[V \xto{f} X]$, where $f$ is projective and $V$ is a derived complete intersection scheme (Proposition \ref{prop:dregchar}). Moreover, it can be shown that the defining relations do not depend on $S$ (see \cite{annala-base-ind-cob}). We expect that $\Omega^S_*(X)$ is a correct model for the \emph{absolute algebraic bordism} of $X$, i.e., the cobordism analogue of the theory of Chow-groups and $K$-theory of coherent sheaves. Some evidence of this is provided by Corollary \ref{cor:bivcobextendsLM}, which states that, for a field $k$ of characteristic 0, $\Omega^k_*$ essentially recovers the algebraic bordism theory of Levine and Morel \cite{levine-morel}.
\end{rem}

\chapter{Projective bundle formula and applications}\label{ch:pbf}

Here, we prove the projective bundle formula for any bivariant theory that is the quotient of universal precobordism $\PCob$. Of course, the analogous results holds then for quotients of the universal $S$-precobordism theories $\PCob^*_S$, in particular for the bivariant algebraic $S$-cobordism $\Omega^*_S$. As an application we construct Chern classes of vector bundles and use them to prove that the zeroth algebraic $K$-theory can be recovered from the cobordism ring. Originally, these results appeared in \cite{annala-yokura, annala-chern, annala-base-ind-cob} in various degrees of generality.

\section{Projective bundle formula}\label{sect:pbf}

In this section, we prove the projective bundle formula. Throughout the section $\Bb$ will denote a bivariant theory that is obtained as a quotient of $\PCob$. 

\subsection{Projective bundle formula for trivial bundles}\label{ssect:trivpbf}

The purpose of this section is to prove the projective bundle formula for trivial projective bundles. The following two preliminary results will be useful.

\begin{lem}\label{lem:pbinj}
Let $X \to Y$ be a morphism of derived schemes, and and let $E$ and $F$ be vector bundles. Then the bivariant pushforward
$$i_*: \Bb\big(\Pb(E) \to Y\big) \to \Bb\big(\Pb(E \oplus F) \to Y\big)$$
is an injection.
\end{lem}
\begin{proof}
We prove the result by showing that $\Pb(E)$ is almost a retract of $\Pb(E \oplus F)$. Let us denote by $\pi: P \to \Pb(E \oplus F)$ the blowup of $\Pb(E \oplus F)$ at $\Pb(F)$. As $\Pb(E)$ and $\Pb(F)$ do not intersect inside $\Pb(E \oplus F)$, we have a Cartesian square
$$
\begin{tikzcd}
\Pb(E) \arrow[d]{}{\Id} \arrow[r]{}{i'} & P \arrow[d]{}{\pi} \\
\Pb(E) \arrow[r]{}{i} & \Pb(E \oplus F)
\end{tikzcd}
$$
On the other hand, by Proposition \ref{prop:pbundvan} and Theorem \ref{thm:dblowupres}, $P$ has a universal property, by which it is the universal $\Pb(E \oplus F)$-scheme on which the natural composition $\Oc(-1) \hook E \oplus F \to E$, factors through a line bundle $\Oc(-1 + \Ec)$ that embeds into $E$. In particular, the inclusion $\Oc(-1 + \Ec) \hook E$ induces a projective map $p: P \to \Pb(E)$, and the composition $p \circ i'$ is the identity. As
\begin{align*}
p_* \circ \pi^* \circ i_* &= p_* \circ i'_* \\
&= \Id,
\end{align*}
it follows that $i_*$ is an injection.
\end{proof}

\begin{lem}\label{lem:pbclasspres}
Let $X \to Y$ be a morphism of derived schemes, and let $E$ a vector bundle on $X$. Let $V \to X$ be a projective morphism such that the composition $V \to Y$ is formally quasi-smooth, and let $s_1,s_2: \Ls \hook E$ be two inclusions of vector bundles on $V$, where $\Ls$ is a line bundle. Then, the maps $s_1$ and $s_2$ induce projective morphisms $f_1, f_2: V \to \Pb(E)$ over $X$, and 
$$[V \stackrel{f_1}{\to} \Pb(E)] = [V \stackrel{f_2}{\to} \Pb(E)] \in \Bb(\Pb(E) \to Y).$$
\end{lem}
\begin{proof}
Let us define a morphism
$$s := x_0 s_1 + x_1 s_2: \Ls \to E(1)$$
of vector bundles on $\Pb^1 \times V$, and let $Z$ be the derived vanishing locus of $s$. By Theorem \ref{thm:dblowupres}, $s$ factors through a natural inclusion $\Ls(\Ec) \hook E(1)$ on the blowup $\bl_Z(\Pb^1 \times V)$, and therefore we obtain a projective morphism $\bl_Z(\Pb^1 \times V) \to \Pb_{\Pb^1 \times X}(E(1))$ such that the composition $\bl_Z(\Pb^1 \times V) \to Y$ is formally quasi-smooth. As the derived vanishing locus $Z$ of $s^\vee$ is disjoint from $0 \times V$ and $\infty \times V$, the above data induces the relation
$$[V \xto{f_1} \Pb(E)] = [V \xto{f_2} \Pb(E)],$$
as desired.
\end{proof}

Consider then the group
$$\Bb(\Pb^\infty \times X \to Y) := \colim_{n \geq 0} \Bb\big(\Pb^n \times X \to Y\big),$$
where the colimit is taken along a sequence of pushforwards along linear embeddings, and let $i_n$ be the natural map $\Bb\big(\Pb^n \times X \to Y\big) \to \Bb(\Pb^\infty \times X \to Y)$. It is clear that $\Bb(\Pb^\infty \times X \to Y)$ admits the structure of a $\Bb^\bullet(X) \dash \Bb^\bullet(Y)$-bimodule. 

By Lemma \ref{lem:pbclasspres}, a precobordism cycle in $\Bb(\Pb^\infty \times X \to Y)$ is uniquely determined by $[V \to X; \Ls]$, where $V \to X$ is a projective morphism such that the composition $V \to Y$ is formally quasi-smooth, and where $\Ls$ (the pullback of $\Oc(1)$) is a globally generated line bundle. Hence, $\PCob(\Pb^\infty \times X \to Y)$ is closely related to the bivariant group with line bundles $\PCob^{\bullet, 1}(X \to Y)$ we studied earlier. Accordingly, we have the following result.

\begin{lem}\label{lem:gglbundform}
Let $X$ be a qcqs derived scheme and let $\Ls$ be a globally generated line bundle. Then,
$$[X \to X; \Ls] = \sum_{i \geq 0} e(\Ls)^i (\gamma_i - \gamma_{i-1} [\Pb(\Ls \oplus \Oc) \to X])\in \Bb(\Pb^\infty \times X \to X),$$
where $\gamma_i = [P_i \times X \to X; \Ms_i]$ in the notation of Definition \ref{def:fundtower}.
\end{lem}
\begin{proof}
This follows from the proof of Lemma \ref{lem:lbundform} with the difference that the line bundle $\Ls' := \Ls(-\Ec)$ on the derived blowup $W := \bl_{\infty \times D}(\Pb^1 \times X)$ should be replaced by $\Ls(1-\Ec)$. Indeed, as $\infty \times D$ is the derived vanishing locus of a section of $\Oc(1) \oplus \Ls$, by Theorem \ref{thm:dblowupres}, there exists an inclusion of vector bundles $\Oc(\Ec) \hook \Oc(1) \oplus \Ls$ on $W$, and therefore a surjection $\Ls \oplus \Oc(1) \to \Ls(1-\Ec)$. It follows that $\Ls(1 - \Ec)$ and all the line bundles $\Ms_i(W, \Ls')$ are globally generated. The proof of Lemma \ref{lem:lbundform} runs through, because $e(\Ls)$ is nilpotent by the global generation hypothesis.
\end{proof}

Hence, we have proven the following result.

\begin{lem}\label{lem:pbgen}
Let $X$ be a qcqs derived scheme. Then, the $\gamma_i$ generate $\Bb(\Pb^\infty \times X \to X)$ as a $\Bb^{\bullet}(X)$-module. Moreover, the same is true for the elements $\rho_i := [\Pb^i \times X \to X; \Oc(1)]$.
\end{lem}
\begin{proof}
The group $\Bb(\Pb^\infty \times X \to X)$ is generated by cycles of form $[V \xto f X; \Ls]$, where $f$ is projective and quasi-smooth, and where $\Ls$ is a globally generated line bundle. As
$$[V \xto f X; \Ls] = (\Id_{\Pb^\infty} \times f)_* \big([V \to V; \Ls] \bullet \theta(f)\big),$$ 
where $(\Id_{\Pb^\infty} \times f)_*: \Bb(\Pb^\infty \times V \to X) \to \Bb(\Pb^\infty \times X \to X)$ is induced by the pushforwards along $\Id_{\Pb^n} \times f: \Pb^n \times V \to \Pb^n \times X$, and as
$$[P_i \times V \to X; \Ms_i] = \gamma_i [V \to X],$$
the first claim follows from Lemma \ref{lem:gglbundform}.

The second claim follows by computing $[\Pb^i \times X \to X; \Oc(1)]$ using Lemma \ref{lem:gglbundform}: indeed, this yields formula for $\rho_i$ as a $\Bb^\bullet(X)$-linear combination of $\gamma_j$, where $j \leq i$, and where the coefficient of $\gamma_i$ is $1$, showing that $\gamma_j$ can be expressed in terms of $\rho_j$.
\end{proof}

We are now ready to compute $\Bb(\Pb^\infty \times X \to Y)$ in terms of $\Bb(\Pb^\infty \times X \to Y)$.

\begin{thm}\label{thm:pinfcomp}
Let $X \to Y$ be a morphism of qcqs derived schemes. Then, the map
$$\bigoplus_{n = 0}^\infty \Bb(X \to Y) \to \Bb(\Pb^\infty \times X \to Y),$$
where the $n^{th}$ morphism is given by the composition
$$i_n(\theta(\pr_2) \bullet -): \Bb(X \to Y) \to \Bb(\Pb^n \times X \to Y) \to \Bb(\Pb^\infty \times X \to Y),$$
is an isomorphism of $\Bb^\bullet(X) \dash \Bb^\bullet(Y)$-bimodules.
\end{thm}
\begin{proof}
The map is clearly a homomorphism of $\Bb^\bullet(X) \dash \Bb^\bullet(Y)$-bimodules, so we only need to show that it is bijective. Observe that the $n^{th}$ map sends $[V \to X] \in \Bb(X \to Y)$ to $[\Pb^n \times V \to X; \Oc(1)] \in \Bb(\Pb^\infty \times X \to Y)$ and therefore the surjectivity follows from Lemma \ref{lem:pbgen}, and from the fact that, for a projective morphism $f: V \to X$ such that the composition $g: V \to Y$ is quasi-smooth, and a globally generated line bundle $\Ls$ on $V$, we have that
$$[V \to X; \Ls] = (\Id_{\Pb^\infty} \times f)_*\big([V \to V; \Ls] \bullet \theta(g)\big),$$
where $(\Id_{\Pb^\infty} \times f)_*$ is as in the proof of Lemma \ref{lem:pbgen}.

To prove the injectivity, let 
$$\alpha := \sum_{n=0}^r i_n(\theta(\pr_2) \bullet \alpha_n),$$
where $\alpha_i \in \Bb(X \to Y)$ are such that $\alpha_r \not = 0$. Then, denoting by $\pr_{2*}$ the ``pushforward'' $\Bb(\Pb^\infty \times X \to Y) \to \Bb(X \to Y)$, we observe that
$$\pr_{2*} \big( e(\Oc(1))^r \bullet \alpha) = \alpha_r,$$
which proves the desired injectivity.
\end{proof}

%
%

As a corollary, we prove projective bundle formula for trivial projective bundles.

\begin{cor}\label{cor:trivpbf}
Let $X \to Y$ be a map of qcqs derived schemes. Then, the map
$$\bigoplus_{i=0}^n \Bb(X \to Y) \to \Bb(\Pb^n \times X \to Y),$$
where the $i^{th}$ map is $e(\Oc(1))^i \bullet \theta(\pr_2) \bullet -$, is an isomorphism.
\end{cor}
\begin{proof}
Since, by Lemma \ref{lem:pbinj}, the natural map $i_n: \Bb(\Pb^n \times X \to Y) \to \Bb(\Pb^\infty \times X \to Y)$ is an injection, and as $i_n \circ (e(\Oc(1))^i \bullet \theta(\pr_2) \bullet -) = i_{n-i}(\theta(\pr_2) \bullet -)$, it suffices to show that 
$$i_{j}(\theta(\pr_2) \bullet -):\bigoplus_{j=0}^n  \Bb(X \to Y) \to \Bb(\Pb^\infty \times X \to Y)$$
surjects onto the image of $i_n$. But this is clear, because by Theorem \ref{thm:pinfcomp}, the map surjects onto the subgroup of those elements that are killed by $e(\Oc(1))^{n+1}$, which contains the image of $i_n$.
\end{proof}

\subsection{Lifting classes to projective bundles}\label{ssect:etaclass}

Before proving the general projective bundle formula, we need to construct an element $\eta_{X,E} \in \Bb^\bullet(\Pb(E))$ that pushes forward to $1_X \in \Bb^\bullet(X)$. By the following result, the existence of such an element implies several useful results. 

\begin{lem}\label{lem:pbpullinj}
Let $f: X' \to X$ be a projective and quasi-smooth map of derived schemes, and suppose that there exists $\eta \in \Bb^\bullet(X')$ such that $f_!(\eta) = 1_X$. Then, 
\begin{enumerate}
\item if $g: X'' \to X'$ is another projective and quasi-smooth map of derived schemes and $\eta' \in \Bb^\bullet(X'')$ satisfies $g_!(\eta') = 1_{X'}$, then 
$$(f \circ g)_! (\eta' \bullet g^*(\eta)) = 1_X;$$

\item the ``pullback'' morphism
$$\theta(f) \bullet - : \Bb(X \to Y) \to \Bb(X' \to Y)$$
is an injection for all maps $X \to Y$;

\item the pullback morphism
$$f^*: \Bb^\bullet(X) \to \Bb^\bullet(X')$$
is an injection;

\item if $\beta \in \Bb^\bullet(X')$ is the inverse of $f^*(\alpha)$, where $\alpha \in \Bb^\bullet(X)$, then $f_!(\beta \bullet \eta)$ is the inverse of $\alpha$ in $\Bb^\bullet(X)$.
\end{enumerate}
\end{lem}
\begin{proof}
\begin{enumerate}
\item Indeed, we may use Proposition \ref{prop:orcohformal} to compute that
\begin{align*}
(f \circ g)_! (\eta' \bullet g^*(\eta)) &= f_! \big( g_!(\eta') \bullet \eta \big) \\
&= f_! \big( \eta \big) \\
&= 1_X,
\end{align*}
as desired.

\item Indeed, for any $\alpha \in \PCob^*(X \to Y)$ we compute that 
\begin{align*}
f_* \big( \eta \bullet \theta(f) \bullet \alpha \big) &= f_* \big( \eta \bullet \theta(f)\big) \bullet \alpha & (A_{12}) \\
&= f_! ( \eta ) \bullet \alpha \\
&= \alpha,
\end{align*}
from which the claim follows.

\item Indeed, we may use Proposition \ref{prop:orcohformal} to compute that, for any $\alpha \in \PCob^*(X)$, we have that
\begin{align*}
f_!(f^*(\alpha) \bullet \eta) &= \alpha \bullet f_!(\eta) \\
&= \alpha,
\end{align*}
from which the claim follows.

\item Indeed, we may use Proposition \ref{prop:orcohformal} to compute that
\begin{align*}
1_X &= f_!(\eta) \\
&= f_!(f^*(\alpha) \bullet \beta \bullet \eta) \\
&= \alpha \bullet f_!(\beta \bullet \eta),
\end{align*}
as desired. \qedhere
\end{enumerate}
\end{proof}

The construction of $\eta_{X,E}$ is based on the following lemma.

\begin{lem}\label{lem:vbundform}
Let $X$ be a finite-Krull-dimensional Noetherian derived scheme that admits an ample line bundle, and let $E$ be a vector bundle on $X$. Then,
$$1_X = \sum_{i=0}^\infty e(E)^i \bullet [\Pb(E \oplus \Oc) \to X]^i \bullet \big([\bl_Z(X) \to X] - [\Pb_\Ec(\Oc(-1) \oplus \Oc) \to X]\big)$$
in $\PCob^*(X)$, where $i: Z \hook X$ is the derived vanishing locus of a section of $E$, and $\Ec \hook \bl_Z(X)$ is the exceptional divisor.
\end{lem}
\begin{proof}
Let $W:=\bl_{\infty \times Z}(\Pb^1 \times X)$. Then the fiber of the canonical map $W \to \Pb^1$ is sum of virtual Cartier divisors $\bl_Z(X)$ and $\Pb_Z(E \oplus \Oc)$ intersecting at $\Ec$. Hence, the $i^{th}$ fiber power of $\Pb_W(E \oplus \Oc)$ over $W$ realizes the relation
\begin{align*}
& [\Pb(E \oplus \Oc) \to X]^i \\
&= [\Pb(E \oplus \Oc) \to X]^i \bullet [\bl_Z(X) \to X] + i_! \big([\Pb_Z(E \oplus \Oc) \to Z]^{i+1}\big)\\
&- [\Pb(E \oplus \Oc) \to X]^i \bullet [\Pb_\Ec(\Oc(-1) \oplus \Oc) \to X] \\
&= [\Pb(E \oplus \Oc) \to X]^i \bullet [\bl_Z(X) \to X] + e(E) \bullet [\Pb(E \oplus \Oc) \to X]^{i+1} \\
&- [\Pb(E \oplus \Oc) \to X]^i \bullet [\Pb_\Ec(\Oc(-1) \oplus \Oc) \to X] \in \PCob^*(X).
\end{align*}
As $e(E)$ is nilpotent (Proposition \ref{prop:eulerclassnilp}), the above formulas may be combined to obtain the desired equation.
\end{proof}

\begin{defn}\label{def:etaclass}
Let $X$ be a finite-Krull-dimensional Noetherian derived scheme that admits an ample line bundle, and let $E$ be a vector bundle on $X$. Then we define
$$\eta_{X,E} := \frac{e(Q) - e(E) \bullet [\Pb_{\Pb(E)}(\Oc(-1) \oplus \Oc) \to \Pb(E)]}{1 - e(E) \bullet [\Pb_{\Pb(E)}(E \oplus \Oc) \to \Pb(E)]} \in \PCob^*(\Pb(E))$$
where 
$$\Oc(-1) \to E \to Q$$ 
is the tautological cofiber sequence of vector bundles on $\Pb(E)$.
\end{defn}

The class $\eta_{X,E}$ has the desired property, as we now verify.

\begin{prop}\label{prop:etaclass}
Let $X$ be a finite-Krull-dimensional Noetherian derived scheme that admits an ample line bundle, let $E$ be a vector bundle on $X$, and denote by $\pi$ the structure morphism $\Pb(E) \to X$. Then,
$$\pi_!(\eta_{X,E}) = 1_X \in \PCob^*(X).$$
\end{prop}
\begin{proof}
By Theorem \ref{lem:dblowuppres}, we have $[\bl_Z(X) \hook \Pb(E)] = e(Q)$ and $$[\Ec \hook \Pb(E)] = e(\Oc(-1)) \bullet e(Q) = e(E).$$ Hence, the claim follows immediately from Lemma \ref{lem:vbundform}.
\end{proof}

\subsection{Projective bundle formula}\label{ssect:pbf}

Here, we prove the projective bundle formula for general projective bundles. Let $X \to Y$ be a map of derived schemes and suppose that $E$ is a vector bundle on $X$. Consider the diagram
$$
\begin{tikzcd}
& \Bb(\Pb^\infty \times X \to Y) \arrow[d]{}{j_E}\\
\Bb(\Pb(E) \to Y) \arrow[r]{}{i'_E} & \colim_n \Bb(\Pb(E \oplus \Oc^{\oplus n}) \to Y).
\end{tikzcd}
$$
By Lemma \ref{lem:pbinj}, $i'_E$ and $j_E$ are injective. In fact, under relatively mild hypotheses, $j_E$ is an isomorphism. 

\begin{lem}
Let $X$ be a finite-Krull-dimensional Noetherian derived scheme that admits an ample line bundle. Then $j_E$ is an isomorphism.
\end{lem}
\begin{proof}
By Lemma \ref{lem:pbclasspres}, the generating cycles of $\colim_n \Bb(\Pb(E \oplus \Oc^{\oplus n})$ are uniquely described by $[V \xto f X; \Ls]$, where $f: V \to X$ is a projective morphism such that the composition $g: V \to Y$ is formally quasi-smooth, and $\Ls$ is a line bundle that admits a surjection from $E^\vee \oplus \Oc^{\oplus N}$ for some $N$. As 
$$[V \xto f X; \Ls] = \tilde f_* ([V \to V; \Ls] \bullet \theta(g)),$$
where $\tilde f_*$ is the map 
$$\colim_n \Bb(\Pb_V(E\vert_V \oplus \Oc^{\oplus n}) \to Y) \to \Bb(\Pb(E \oplus \Oc^{\oplus n}) \to Y)$$
induced by the obvious pushforwards, it is enough to show that $[V \to V; \Ls]$ is in the image of $j_{E\vert_V}$.

In order to prove the claim, we mimic the proof of Lemma \ref{lem:lbundform}. Consider a surjection $E^\vee \oplus \Oc^{\oplus N} \to \Ls$ on $V$, and let $D \hook V$ be the derived vanishing locus of a section of $\Ls$. Let $W:=\bl_{\infty \times D}(\Pb^1 \times V)$ and let $\Ls' := \Ls(1 - \Ec)$ on $W$. By Theorem \ref{thm:dblowupres}, there exists a surjection of vector bundles $\Oc(-1) \oplus \Ls^\vee \to \Oc(-\Ec)$ on the blowup $W$, which can be used to construct a surjection
$$E^\vee \oplus \Oc^{\oplus N+2} \to \Ls \oplus \Oc(1) \to \Ls(1-\Ec).$$
Moreover, given a line bundle $\Ms$ that admits a surjection $\Ec^\vee \oplus \Oc^{\oplus r} \to \Ms$, then, on $\Pb(\Ms \oplus \Oc)$, the tautological surjection $\Ms^\vee \oplus \Oc \to \Oc(1)$ can be used to construct a surjection 
$$E^\vee \oplus \Oc^{\oplus r+1} \to \Oc \oplus \Ms \to \Ms(1).$$
Hence, using the notation of Definition \ref{def:fundtower}, the line bundles $\Ms_i(W,\Ls')$ on the derived schemes $P_i(W, \Ls')$ admit a surjection from $E^\vee \oplus \Oc^{\oplus m}$ for $m \gg 0$.

Hence, as in the proof of Lemma \ref{lem:lbundform} the natural maps $P_i(W, \Ls') \to \Pb^1 \times V$, together with a lift $P_i(W, \Ls') \to \Pb_{\Pb^1 \times V} (E \oplus \Oc^{\oplus m})$ for $m \gg 0$ such that $\Oc(1)$ pulls back to $\Ms_i(W,\Ls')$, realize relations
\begin{align*}
&[P_i(V, \Ls) \to V; \Ms_i(V, \Ls)] \\
&= [P_i \times V \to V; \Ms_i] + e(\Ls) \bullet [P_{i+1}(V, \Ls) \to V; \Ms_{i+1}(V, \Ls)] \\
&- [P_{i+1} \times V \to V; \Ms_{i+1}] \bullet e(\Ls) \bullet [\Pb(\Ls \oplus \Oc) \to V]
\end{align*}
in $\colim_n \Bb(\Pb(E \oplus \Oc^{\oplus n}) \to V)$, which, by the nilpotence of $e(\Ls)$, may be used to express $[V \to V; \Ls]$ as a $\PCob^*(V)$-linear combination of $[P_i \times V \to V; \Ms_i]$. Since the line bundles $\Ms_i$ are globally generated, the classes $[P_i \times V \to V; \Ms_i]$ lie in the image of $j_{E \vert_V}$. This implies that $[V \to V; \Ls]$ lies in the image of $j_{E|_V}$.
\end{proof}

\begin{defn}\label{def:fundemb}
Let $X \to Y$ be a map derived schemes, and suppose that $X$ is Noetherian, of finite Krull dimension, and admits an ample line bundle. Let $E$ be a vector bundle on $X$. Then we define the \emph{fundamental embedding}
$$i_E: \Bb(\Pb(E) \to Y) \hook \Bb(\Pb^\infty \times X \to Y)$$
as the composition $j^{-1}_E \circ i'_E$. Clearly, $i_E$ is an injection.

Hence, for any $\alpha \in \Bb(\Pb(E) \to Y)$
$$i_E(\alpha) = \sum_{j=0}^N i_j\big(\theta(\pr_2) \bullet u_j(\alpha) \big)$$
for some uniquely determined $u_j(\alpha) \in \Bb(X \to Y)$ called the \emph{coefficients} of $\alpha$; moreover, these coefficients uniquely determine $\alpha$.
\end{defn}

The coefficients of bivariant elements play a fundamental role in the proof of projective bundle formula. We begin our attack by making the following observation.

\begin{lem}\label{lem:rcoeffs1}
Let $X$ be a finite-Krull-dimensional Noetherian derived scheme that admits an ample line bundle, and let $E$ be a vector bundle of rank $r$ on $X$. Then, there exists $d_i(E) \in \Bb^\bullet(X)$ such that
$$e(\Oc(1))^r = \sum_{i=1}^{r} d_i(E) \bullet e(\Oc(1))^{r-i} \in \Bb^\bullet(\Pb(E)).$$
\end{lem}
\begin{proof}
Let $\pi$ be the natural projection $\Fl(E) \to \Pb(E)$ and let $\eta \in \Bb^\bullet(\Fl(E))$ be such that $\pi_!(\eta) = 1_{\Pb(E)}$. Such an $\eta$ can be found by \ref{prop:etaclass} as $\pi$ is a tower of projective bundles.
Let 
$$0 = E_0 \subset E_1 \subset \cdots \subset E_{r-1} \subset E_r = E$$
be a filtration of $E$ on $\Fl(E)$ with line bundle quotients $\Ls_i := E_{i}/E_{i-1}$. Then, since
\begin{align*}
e(\Ls_1(1)) \bullet \cdots \bullet e(\Ls_r(1)) &= e(E(1)) \\
&= 0
\end{align*} 
it follows that also
\begin{align*}
e(\Ls_1^\vee(-1)) \bullet \cdots \bullet e(\Ls^\vee_r(-1)) &= \inv_F \big(e(\Ls_1(1))\big) \bullet \cdots \bullet \inv_F \big( e(\Ls_r(1))\big) \\
&= 0,
\end{align*}
where $\inv_F$ is the \emph{formal inverse} power series for the formal group law of Theorem \ref{thm:fgl}. Hence,
$$\prod_{i=1}^r \big( F(e(\Ls_i^\vee(-1)), e(\Oc(1))) - e(\Oc(1)) \big) = \prod_{i=1}^r \big( e(\Ls_i^\vee) - e(\Oc(1)) \big)$$
vanishes. In conclusion,
$$\sum_{i=0}^r (-1)^{r-i} s_i\big(e(\Ls_1^\vee),...,e(\Ls_r^\vee)\big) \bullet e(\Oc(1))^{r-i} = 0 \in \Bb^\bullet(\Fl(E)),$$
where $s_i$ are the elementary symmetric polynomials, and the claim follows with 
$$d_i(E) := (-1)^{r-i+1} \pi_!\Big( s_i\big(e(\Ls_1^\vee),...,e(\Ls_r^\vee)\big) \bullet \eta \Big)$$
from the projection formula.
\end{proof}

The above result can be used to show that the first $r$ coefficients of a bivariant class $\alpha \in \Bb(\Pb(E) \to Y)$, $r$ being the rank of $E$, completely determine $\alpha$. 

\begin{lem}\label{lem:rcoeffs2}
Let $X$ be a finite-Krull-dimensional Noetherian derived scheme that admits an ample line bundle, and let $E$ be a vector bundle of rank $r$ on $X$. Then, if $\alpha \in \Bb(\Pb(E) \to Y)$ is such that $u_i(\alpha) = 0$ for all $i < r$, we have that $\alpha = 0$.
\end{lem}
\begin{proof}
Clearly
$$e(\Oc(1))^r \bullet \iota_E(\alpha) = \sum_{j} i_{j}\big(\theta(\pr_2) \bullet u_{j+r}(\alpha) \big).$$
On the other hand, by Lemma \ref{lem:rcoeffs1}, we can also compute that
\begin{align*}
e(\Oc(1))^r \bullet \iota_E(\alpha) &= \sum_{i=1}^r d_i(E) \bullet e(\Oc(1))^{r-i} \bullet \iota_E(\alpha)
\end{align*}
and therefore
$$u_{r+i}(\alpha) = \sum_{j=0}^{r-1} d_{r-j}(E) \bullet u_{j+i}(\alpha).$$
Hence, if $u_{i}(\alpha) = 0$ for all $i < r$, it follows that all the coefficients $u_i(\alpha)$ vanish, and therefore $\alpha=0$.
\end{proof}

Hence, the map
$$
\begin{bmatrix}
u_0(-) \\
\vdots \\
u_{r-1}(-)
\end{bmatrix}:
\Bb^\bullet(\Pb(E) \to X) \to \Bb^\bullet(X)^{\oplus r}
$$
is an injection of $\Bb^\bullet(X)$-modules. The following observations allow us to find a $\Bb^\bullet(X)$-linear basis for $\Bb^\bullet(\Pb(E) \to X)$.

\begin{lem}\label{lem:coeffobs}
Let $X$ be a finite-Krull-dimensional Noetherian derived scheme that admits an ample line bundle, and let $E$ be a vector bundle of rank $r$ on $X$. Denote by $\pi$ the structure morphism $\Pb(E) \to X$. Then $u_{r-1}(\theta(\pi)) \in \Bb^\bullet(X)$ is a unit and $u_{i}(\theta(\pi)) \in \Bb^\bullet(X)$ is nilpotent for all other $i$.
\end{lem}
\begin{proof}
We start with the case of a split vector bundle, i.e.,
$$E \simeq \bigoplus_{i=1}^r \Ls_i.$$
By construction of $i_E$, it is enough to understand the image of $\theta(\pi)$ inside
$$\Bb^\bullet(\Pb(\Oc^\infty \oplus E) \to X) :=  \colim_{n \geq 0} \Bb^\bullet(\Pb(\Oc^n \oplus E) \to X).$$
If $E_I = \bigoplus_{i \in I} \Ls_i$ for some $I \subset \{1,...,r\}$ and $n \geq 0$, we will denote by 
$$[\Pb(\Oc^{\oplus n} \oplus E') \hook \Pb(\Oc^\infty \oplus E)]$$
the class of the obvious linear embedding. 

Let $E' := \bigoplus_{i=2}^r \Ls_i$ and consider $\Pb(\Oc \oplus \Ls_1 \oplus E')$. As $\Ls_1(1)$ has a section with derived vanishing locus $\Pb(\Oc \oplus E')$ and $\Oc(1)$ has a section with derived vanishing locus $\Pb(\Ls_1 \oplus E')$ we may use the $\Lb$-linear structure of $\Bb^*$ provided by the formal group law (Theorem \ref{thm:fgl}) to compute that
\begin{align*}
&[\Pb(\Oc \oplus E') \hook \Pb(\Oc^\infty \oplus E)] \\
=& [\Pb(\Ls_1 \oplus E') \hook \Pb(\Oc^\infty \oplus E)] \\
+& e(\Ls_1) \bullet [\Pb(\Oc \oplus \Ls_1 \oplus E') \hook \Pb(\Oc^{\oplus\infty} \oplus E)] \\
+& \sum_{i,j \geq 1} a_{ij} e(\Oc(1))^{i-1} \bullet e(\Ls_1)^j \bullet [\Pb(\Ls_1 \oplus E') \hook \Pb(\Oc^{\oplus \infty} \oplus E)].
\end{align*}
From the above, we can solve
\begin{align*}
&[\Pb(\Ls_1 \oplus E') \hook \Pb(\Oc^\infty \oplus E)] \\
&= H\big(e(\Oc(1)), e(\Ls_1)\big) \bullet \Big( [\Pb(\Oc \oplus E') \hook \Pb(\Oc^\infty \oplus E)]  \\
&- e(\Ls_1) \bullet [\Pb(\Oc \oplus \Ls_1 \oplus E') \hook \Pb(\Oc^{\oplus\infty} \oplus E)]\Big),
\end{align*}
where 
$$H(t,x) := {1 \over 1 +  \sum_{i,j \geq 1} a_{ij} t^{i-1} \bullet x^j} \in \Lb[[t,x]].$$
Hence 
\begin{align*}
&[\Pb(\Ls_1 \oplus E') \hook \Pb(\Oc^\infty \oplus E)] \\
=& (-1)^{i_1-1} \sum_{i_1 \geq 1} H\big(e(\Oc(1)), e(\Ls_1)\big)^{i_1} \bullet e(\Ls_1)^{i_1 - 1} \\
& \bullet [\Pb(\Oc^{\oplus i_1} \oplus E') \hook \Pb(\Oc^{\oplus\infty} \oplus E)].
\end{align*}
The other $\Ls_i$ can be dealt with in a similar fashion, and therefore $[\Pb(E) \hook \Pb(\Oc^\infty \oplus E)]$ equals to
\begin{align*}
&\sum_{i_1 \cdots i_r \geq 1} \Bigg( \prod_{k=1}^r (-1)^{i_k-1} \Big( H\big(e(\Oc(1)), e(\Ls_k)\big)^{i_k} \bullet e(\Ls_k)^{i_k - 1} \Big) \\
&\bullet [\Pb(\Oc^{\oplus i_1 + \cdots + i_r}) \hook \Pb(\Oc^{\oplus\infty} \oplus E)] \Bigg).
\end{align*}
From the above we can solve that $u_i(\theta(\pi)) = H\big(e(\Ls_1),...,e(\Ls_r)\big)$, where 
$$H_i(x_1,...,x_r) \in \Lb[[x_1,...,x_r]]$$
are universal symmetric power series with $\Lb$-coefficients. Moreover $H_{r-1}$ has constant coefficient 1, while the other $H_i$ have no constant coefficients, proving the claim in the case of a split vector bundle.

We are left with the general case. As one easily checks that coefficients are compatible with pullbacks, we may pull $E$ back along the structure morphism $\Fl(E) \to X$, $\Fl(E)$ being the flag bundle, where $E$ admits a filtration 
$$0 = E_0 \subset E_1 \subset \cdots \subset E_{r-1} \subset E_{r} = E$$
with line bundle quotients $\Ls_i := E_{i}/E_{i-1}$. Since $\Fl(E) \to X$ is a tower of projective bundles, Lemma \ref{lem:pbpullinj} and Proposition \ref{prop:etaclass} imply that the coefficient $u_i(\theta(\pi))$ is invertible or nilpotent in $\Bb^\bullet(X)$ if and only if it is invertible or nilpotent in $\Bb^\bullet(\Fl(E))$, respectively. Denoting by $E'$ the sum $\bigoplus_{i=1}^r \Ls_i$, by $\pi'$ the structure map $\Pb(E') \to X$, and by $\psi_1$ and $\psi_2$ the natural inclusions $\Pb(E) \hook \Pb(E \oplus E')$ and $\Pb(E') \hook \Pb(E \oplus E')$ respectively, we may compute that
\begin{align*}
i_E(\theta(\pi)) &= i_{E \oplus E'} \big(\psi_{1*}(\theta(\pi))\big) \\
&= i_{E \oplus E'} \big(e(E'(1))\big) \\
&= i_{E \oplus E'} \big(e(\Ls_1(1)) \bullet \cdots \bullet e(\Ls_r(1))\big) \\
&= i_{E \oplus E'} \big(e(E(1))\big) \\
&= i_{E \oplus E'} \big(\psi_{2*}(\theta(\pi'))\big) \\
&= i_{E'}(\theta(\pi')),
\end{align*}
so the general case follows from the split case.
\end{proof}

\begin{lem}\label{lem:coeffobs2}
Let $X$ be a finite-Krull-dimensional Noetherian derived scheme that admits an ample line bundle, and let $E$ be a vector bundle of rank $r$ on $X$. Then, for all $i \leq r-1$ and $j \leq r-1$, there exist unique $\alpha_{j,i}(E) \in \Bb^\bullet(X)$, which satisfy
\begin{align*}
u_k \Bigg(\sum_{j=0}^{r-1} e(\Oc(1))^j \bullet \theta(\pi) \bullet \alpha_{j,i}(E) \Bigg) &= 
\begin{cases}
1_X & \text{if k=i;} \\
0 & \text{if $k<r$ and $k \not = i$,}
\end{cases}
\end{align*} 
where $\pi$ is the structure morphism $\Pb(E) \to X$. Moreover, $[\alpha_{j,i}(E)]$ is an invertible matrix with coefficients in $\Bb^\bullet(X)$.
\end{lem}
\begin{proof}
Consider the matrix $A(E)$, whose entries are
\begin{align*}
A(E)_{j,i} :=& u_j\big(e(\Oc(1))^i \bullet \theta(\pi) \big)  \\
=& u_{i+j}\big(\theta(\pi)\big).
\end{align*}
The matrix $A(E)$ can be interpreted as follows: given
$$\overline{\alpha} = 
\begin{bmatrix}
\alpha_0 \\
\vdots \\
\alpha_{r-1}
\end{bmatrix} 
\in \Bb^*(X)^{\oplus r},
$$
then $A(E) \overline{\alpha}$ is the vector of first $r$ coefficients of $\sum_{i=0}^{r-1} e(\Oc(1))^i \bullet \theta(\pi) \bullet \alpha_i$. By Lemma \ref{lem:coeffobs}, $A(E)$ has units on the anti-diagonal and nilpotent elements everywhere else, and therefore, by basic linear algebra, $A(E)$ is invertible. Then $[\alpha_{j,i}(E)] := A(E)^{-1}$ is the unique matrix with the desired property.
\end{proof}

The following lemma establishes the fact that for any choice for the first $r$ coefficients, there exists a bivariant element in $\Bb(X \to Y)$ realizing them.

\begin{lem}\label{lem:pbfsurj}
Let $X \to Y$ be a morphism of derived schemes, suppose that $X$ is Noetherian, of finite Krull dimension, and admits an ample line bundle, and let $E$ be a vector bundle of rank $r$ on $X$. Moreover, let $u_0,...,u_{r-1} \in \Bb(X \to Y)$. Then
$$\alpha := \sum_{i=0}^{r-1} \sum_{j=0}^{r-1} e(\Oc(1))^j \bullet \theta(\pi) \bullet \alpha_{j,i}(E) \bullet u_i$$
is the unique element of $\Bb(\Pb(E) \to X)$ with $u_i(\alpha) = u_i$ for all $i < r$.
\end{lem}
\begin{proof}
This follows immediately from Lemma \ref{lem:coeffobs2}.
\end{proof}

We finally have all the ingredients to prove the projective bundle formula.

\begin{thm}\label{thm:pbf}
Let $X \to Y$ be a morphism of derived schemes, suppose that $X$ is Noetherian, of finite Krull dimension, and admits an ample line bundle, and let $E$ be a vector bundle of rank $r$ on $X$. Then, for any bivariant quotient $\Bb$ of the universal precobordism theory $\PCob$, the map
$$e(\Oc(1))^i \bullet \theta(\pi) \bullet - :\bigoplus_{i=0}^{r-1} \Bb(X \to Y) \to \Bb(\Pb(E) \to Y)$$
is an isomorphism, where $\pi$ is the structure morphism $\Pb(E) \to X$.
\end{thm}
\begin{proof}
By Lemma \ref{lem:pbfsurj}, for each $u_0,...,u_{r-1} \in \Bb(X \to Y)$, the element 
$$\alpha := \sum_{j=0}^{r-1} e(\Oc(1))^j \bullet \theta(\pi) \bullet \Bigg( \sum_{i=0}^{r-1} \alpha_{j,i}(E) \bullet u_i \Bigg)$$
has $u_i(\alpha) = u_i$ for $i < r$. By Lemma \ref{lem:coeffobs2}, any element of $\Bb(\Pb(E) \to Y)$ is uniquely determined by its first $r$ coefficients, and by Lemma \ref{lem:coeffobs2}, the matrix $[\alpha_{j,i}(E)]$ is invertible. The desired bijectivity follows from this.
\end{proof}

\section{Chern classes}

The purpose of this section is to define Chern classes of vector bundles on Noetherian derived schemes admitting an ample line bundle.

\begin{defn}\label{def:chernclass}
Let $X$ be a finite-Krull-dimensional Noetherian derived scheme that admits an ample line bundle, and let $E$ be a vector bundle of rank $r$ on $X$. Then, by the projective bundle formula (Theorem \ref{thm:pbf}), the equation
$$\sum_{i=0}^r (-1)^i e(\Oc(1))^{i} \bullet c_{r-i}(E) = 0 \in \PCob^r(\Pb(E^\vee))$$
holds for uniquely defined $c_i(E) \in \PCob^i(X)$ with $c_0(E)=1_X$. The element $c_i(E)$ is called the \emph{$i^{th}$ Chern class} of $E$. 
\end{defn} 

Chern classes satisfy the expected properties, as the following result shows.

\begin{thm}\label{thm:chernclass}
Define the \emph{total Chern class} of a rank $r$ vector bundle $E$ as
$$c(E) = 1 + c_1(E) + \cdots + c_r(E).$$
Then, the Chern classes satisfy the following properties:
\begin{enumerate}
\item \emph{naturality:} $f^* c_i(E) = c_i(f^*E)$, where $f: Y \to X$ is a map of finite-Krull-dimensional Noetherian derived schemes admitting ample line bundles;

\item \emph{Whitney sum formula:} given a cofiber sequence 
$$E' \to E \to E''$$
of vector bundles on $X$, the equation
$$c(E) = c(E') \bullet c(E'') \in \PCob^*(X)$$ 
holds;

\item \emph{normalization:} if $E$ is a rank $r$ vector bundle, then $c_r(E) = e(E)$.
\end{enumerate}
\end{thm}
\begin{proof}
Naturality is trivial. To prove the Whitney sum formula, we can pull back to flag bundles to reduce to the case where both $E'$ and $E''$ admit filtrations
$$0 = E'_0 \subset E'_1 \subset \cdots \subset E'_{r'} = E'$$
and
$$0 = E''_0 \subset E''_1 \subset \cdots \subset E''_{r''} = E''$$
with graded pieces line bundles $\Ls'_i := E'_i/E'_{i-1}$ and $\Ls''_i := E''_i / E''_{i-1}$ respectively. One then argues as in Lemma \ref{lem:rcoeffs1} to show that
$$c_i(E') = s_i\big(e(\Ls'_1), ..., e(\Ls'_{r'}) \big),$$
$$c_i(E'') = s_i\big(e(\Ls''_1), ..., e(\Ls''_{r''}) \big)$$
and 
$$c_i(E) = s_i\big(e(\Ls'_1), ..., e(\Ls'_{r'}), e(\Ls''_1), ..., e(\Ls''_{r''}) \big)$$
proving the Whitney sum formula. Normalization is trivial for line bundles, and follows for general vector bundles from the Whitney sum formula after pulling back to a flag bundle.
\end{proof}

\section{Cohomological Conner--Floyd theorem}\label{sect:cf}

Here, we prove the following analogue of a well-known theorem of Conner and Floyd \cite{conner:1966,conner:1969}.

\begin{thm}\label{thm:cf}
Let $X$ be a finite-Krull-dimensional Noetherian derived scheme that admits an ample line bundle. Then, the map
$$\eta_K : \Zb_m \otimes_\Lb \PCob^*(X) \to K^0(X)$$
defined by the formula
$$[V \xto{f} X] \mapsto [f_* \Oc_V]$$
is an isomorphism, where $\Zb_m$ is the integers, considered as an $\Lb$-algebra via the multiplicative formal group law $x + y - xy$, and where $K^0(X)$ denotes the $K$-theory of perfect complexes on $X$.
\end{thm}

Our result generalizes that of Levine and Morel \cite{levine-morel}, which holds for smooth varieties over a field of characteristic 0. There exists another generalization of the aforementioned result by Dai \cite{dai:2010}, which states that $G$-theory (also known as the $K$-theory of coherent sheaves and $K$-theory homology) of a quasi-projective scheme over a field of characteristic 0 can be recovered from its algebraic bordism by enforcing the multiplicative formal group law. 

\subsection{Bivariant $K$-theory}

Here, we concisely recall bivariant algebraic $K$-theory in the context of derived algebraic geometry. In classical algebraic geometry, this theory was introduced in \cite{fulton-macpherson}.

\begin{defn}\label{def:bivKthy}
Let $f: X \to Y$ be a map of derived schemes. Then, a quasi-coherent sheaf $\Fc \in \QCoh(X)$ is \emph{$f$-perfect}, or \emph{relatively perfect}, if $\Fc$ is an almost perfect object of $\QCoh(X)$, and has a finite Tor-amplitude as an $f^{-1} \Oc_Y$-module\footnote{In other words, for each affine open $\Spec(A) \subset Y$ and an affine open $\Spec(B) \subset X$ mapping to $\Spec(A)$, the derived $B$-module $\Fc \vert_{\Spec(B)}$ has finite Tor-amplitude as a derived $A$-module.}. 

The \emph{bivariant $K$-theory group} $K^0(X \xto{f} Y)$ is defined as the Abelian group generated by the equivalence classes of $f$-perfect objects, modulo the relations $[\Fc] = [\Fc'] + [\Fc'']$ whenever
$$\Fc' \to \Fc \to \Fc''$$
is a cofiber sequence of $f$-perfect objects in $\QCoh(X)$. 
\end{defn} 

The bivariant operations are defined as follows.

\begin{enumerate}
\item \emph{Bivariant pushforward}: Suppose $X \to Y$ factors through a proper morphism $g: X \to X'$ that is locally of almost finite presentation. Then, given $[\Fc] \in K^0(X \to Y)$, we define
$$f_*([\Fc]) := [f_* \Fc] \in K^0(X' \to Y).$$
Note that $f_* \Fc$ is relatively perfect over $Y$ because $f_*$ preserves almost perfect objects (\cite{SAG} Theorem 3.7.0.2), and because for any affine open $U \subset Y$ and an affine open $V \subset X'$ mapping to $U$, $f_*(\Fc)(V) = \Fc(f^{-1}V)$ has finite Tor-amplitude as a finite limit of modules of finite Tor-amplitude.

\item \emph{Bivariant pullback}: Let 
$$
\begin{tikzcd}
X' \arrow[d]{}{g'} \arrow[r] & Y' \arrow[d]{}{g} \\
X \arrow[r] & Y
\end{tikzcd}
$$
be a Cartesian square, and suppose $[\Fc] \in K^0(X \to Y)$. Then, we define
$$g^*([\Fc]) := [g'^* \Fc] \in K^0(X' \to Y').$$

\item \emph{Bivariant product}: Let $X \xto{f} Y \xto{g} Z$ be a sequence of morphisms of derived schemes. Then, for any $[\Fc] \in K^0(X \to Y)$ and $[\Gc] \in K^0(Y \to Z)$, we define
$$[\Fc] \bullet [\Gc] := [\Fc \otimes f^* \Gc] \in K^0(X \to Z).$$
Note that $\Fc \otimes f^* \Gc$ is almost perfect as a tensor product of almost perfect objects (\cite{HA} Proposition 7.2.4.11 and \cite{HTT} 5.5.8.6). That it has finite Tor-amplitude over $Z$ follows from the basic properties of tensor products of derived modules.
\end{enumerate}

Moreover, for a formally quasi-smooth morphism $f: X \to Y$, we define
$$\theta(f) := [\Oc_X] \in K^0(X \to Y).$$

\begin{prop}\label{prop:bivKthy}
With the above definitions, $K^0$ is a stably oriented bivariant theory with functoriality $\Fc_a$, as defined in Section \ref{ssect:cobbivfunc}. Hence, there exists a unique orientation preserving Grothendieck transformation
$$\eta_K: \PCob^* \to K^0,$$
which, on a morphism $X \to Y$, is defined by the formula
$$[V \xto{f} X] \mapsto [f_* \Oc_V].$$
\end{prop}
\begin{proof}
The proof that bivariant algebraic $K$-theory gives a bivariant theory on classical schemes is essentially contained in \cite{berthelot:1971}, as noted in \cite{fulton-macpherson}. We omit the proof in the derived setting.
\end{proof}

\subsection{Proof of Theorem \ref{thm:cf}}

Here, we prove the Conner--Floyd theorem. An important fact that we are going to use is that, for a derived scheme $X$ admitting an ample line bundle, the Grothendieck ring $K^0(X)$ has a presentation as the Abelian group generated by equivalence classes $[E]$ of vector bundles, modulo the relations $[E] = [E'] + [E'']$ whenever $E' \to E \to E''$ is a cofiber sequence of vector bundles. For the proof of this in the context of classical algebraic geometry, see e.g. \cite{TT}, and in the context of derived geometry, see \cite{annala-qpnote}. With this fact in mind, we make the following definition.

\begin{defn}\label{def:multcc}
Let $X$ be a finite-Krull-dimensional Noetherian derived scheme that admits an ample line bundle. Then, we define the map
$$\ch_m: K^0(X) \to \Zb_m \otimes_\Lb \PCob^*(X)$$
by the formula
$$[E] \mapsto \rank(E) - c_1(E^\vee),$$
which is well defined due to the additivity of first Chern classes in cofiber sequences of vector bundles.
\end{defn}

The map $\ch_m$ satisfies many pleasant properties which will be useful later.

\begin{lem}\label{lem:multccprops}
Let $X$ be a finite-Krull-dimensional Noetherian derived scheme that admits an ample line bundle. Then,
\begin{enumerate}
\item $\ch_m: K^0(X) \to \Zb_m \otimes_\Lb \Omega^*(X)$ is a map of rings;

\item if $f$ is a map $X \to Y$ of finite-Krull-dimensional Noetherian derived schemes that admit an ample line bundle, then $\ch_m \circ f^* = f^* \circ \ch_m$;

\item $\eta_K \circ \ch_m = \Id: K^0(X) \to K^0(X)$;

\item $\ch_m$ sends Chern classes to Chern classes.
\end{enumerate}
\end{lem}
\begin{proof}
\begin{enumerate}
\item Indeed, using the well-known formula
$$c_1(E \otimes F) = c_1(E) + c_1(F)- c_1(E) c_1(F)$$
for $K$-theoretic Chern classes, we compute that
\begin{align*}
\ch_m([E \otimes F]) &= \rank(E) \rank(F) - c_1(E^\vee \otimes F^\vee) \\
&= \rank(E) \rank(F) - c_1(E^\vee) - c_1(F^\vee) + c_1(E^\vee) c_1(F^\vee) \\
&= \ch_m([E])\ch_m([F]),
\end{align*}
proving that $\ch_m$ is compatible with the multiplication. It is trivial that $\ch_m$ is compatible with addition.

\item This follows immediately from the naturality of Chern classes.

\item This follows immediately from the fact that, in $K$-theory, first Chern classes satisfy the formula
$$c_1(E) = \rank(E) - [E^\vee].$$

\item If $\Ls$ is a line bundle, then
\begin{align*}
\ch_m(c_1(\Ls)) &= \ch_m([\Oc_X]) - \ch_m([\Ls^\vee]) \\
&= 1_X - (1_X - c_1(\Ls)) \\
&= c_1(\Ls),
\end{align*}
proving the claim for line bundles. The general claim follows from this by the splitting principle, and from the fact that $\ch_m$ is a map of rings. \qedhere
\end{enumerate}
\end{proof}

The final ingredient required for the proof is that $\ch_m$ commutes with Gysin pushforwards. Due to the global factorization results of projective morphisms between Noetherian schemes admitting an ample line bundles, this splits into two cases: pushforwards along the projection $\Pb^n \times X \to X$ and pushforwards along derived regular embeddings. The following result takes care of the first case.

\begin{lem}\label{lem:ccgysinproj}
Let $X$ be a finite-Krull-dimensional Noetherian derived scheme that admits an ample line bundle. Then, $\ch_m$ commutes with Gysin pushforward along $\pr_2: \Pb^n \times X \to X$.
\end{lem}
\begin{proof}
Since for any $K$-theory class $\alpha \in K^0(X)$, we have 
$$
\pi_*\big(c_1(\Oc(1))^i \bullet \pr_2^*(\alpha) \big) =
\begin{cases}
\alpha & \text{if $i \leq n$};\\
0 & \text{otherwise},
\end{cases}
$$
it suffices to show that the analogous formula holds for $\Zb_m \otimes_\Lb \Omega^*$. By projection formula, we conclude that we want to show that
$$\pi_*\big(c_1(\Oc(1))^i\big) 
\begin{cases}
1_X & \text{if $i \leq n$};\\
0 & \text{otherwise}.
\end{cases}
$$
In other words, we aim to show that for all $m \geq 0$, $[\Pb^m \times X \to X] = 1_X \in \PCob^*(X)$.

To do so, we use an argument from the Proof of Lemma 4.2.3 in \cite{levine-morel}. In order to make the notation less heavy, we will denote for each vector bundle $E$ on $\Pb^n \times X$ by $[\Pb_{\Pb^n}(E)]$ the element $[\Pb_{\Pb^n \times X}(E) \to X] \in \PCob^*(X)$; moreover, the tautological hyperplane bundle will be denoted by $\Oc(H)$ in order to distinguish it from the $\Oc(1)$ of $\Pb^n$.

Analyzing the Chern class of $\Oc(1+H)$ on $\Pb_{\Pb^n}(\Oc^{\oplus m} \oplus \Oc(1)^{\oplus l})$, we arrive at the formula
\begin{align*}
&[\Pb_{\Pb^n}(\Oc^{\oplus m} \oplus \Oc(1)^{\oplus l - 1})]\\
&= [\Pb_{\Pb^n}(\Oc^{\oplus m - 1} \oplus \Oc(1)^{\oplus l})] + [\Pb_{\Pb^{n-1}}(\Oc^{\oplus m} \oplus \Oc(1)^{\oplus l})]  \\
&- [\Pb_{\Pb^{n-1}}(\Oc^{\oplus m - 1} \oplus \Oc(1)^{\oplus l})]. 
\end{align*}
Using this repeatedly, and using the convention that $\Pb^{-1} = \emptyset$, we compute that
\begin{align*}
[\Pb_{\Pb^n}(\Oc)] &= [\Pb_{\Pb^n}(\Oc(1))] + [\Pb_{\Pb^{n-1}}(\Oc \oplus \Oc(1))] - [\Pb_{\Pb^{n-1}}(\Oc(1))] \\
&= \sum_{i = 0}^{n} [\Pb_{\Pb^{n-i}}(\Oc(1)^{\oplus i + 1})] -  [\Pb_{\Pb^{n-i-1}}(\Oc(1)^{\oplus i + 1})] \\
&= \sum_{i = 0}^{n} [\Pb^i_X] ( [\Pb^{n-i}_X] - [\Pb^{n-i-1}_X]).
\end{align*}
We are now ready to argue by induction: if it is known that, for all $i<n$, $[\Pb^i_X] = 1_X$, then most of the terms in the above equation vanish, leaving us only with
$$[\Pb^n_X] = [\Pb^n_X] - 1_X + [\Pb^n_X],$$
from which we can deduce that $[\Pb^n_X] = 1_X$, as desired.
\end{proof}

Next, we prove that $\ch_m$ commutes with Gysin pushforwards along regular embeddings. We begin with the following observation.

\begin{lem}\label{lem:ccgysinregemb1}
Let $X$ be a finite-Krull-dimensional Noetherian derived scheme that admits an ample line bundle, and let $\pi: P \to X$ be a quasi-projective morphism admitting a section $i: X \hook P$, which can be identified as the derived vanishing locus of a section of a vector bundle $E$ on $P$. Then $\ch_m$ commutes with $i_!$. 
\end{lem}
\begin{proof}
Indeed, for all $\alpha \in K^0(X)$, 
\begin{align*}
i_!(\ch_m(\alpha)) &= i_!\big(i^*(\pi^*(\ch_m(\alpha)))\big) \\
&= \pi^*(\ch_m(\alpha)) \bullet i_!(1_X) \\
&= \pi^*(\ch_m(\alpha)) \bullet e(E) \\
&= \ch_m\big(\pi^*(\alpha) \bullet e(E) \big) \\
&= \ch_m(i_!(\alpha)),
\end{align*} 
as desired.
\end{proof}
We are now ready to prove the case of a general regular embedding.

\begin{lem}\label{lem:ccgysinregemb2}
Let $X$ be a finite-Krull-dimensional Noetherian derived scheme that admits an ample line bundle, and let $i: Z \hook X$ be a derived regular embedding. Then $\ch_m$ commutes with $i_!$.
\end{lem}
\begin{proof}
Consider the Cartesian diagram
$$
\begin{tikzcd}
\infty \times Z \arrow[r]{}{i_\infty} \arrow[d]{}{j_\infty} & \Pb_Z( \Oc \oplus \Nc_{Z / X}) \arrow[d]{}{\iota} \\
\Pb^1 \times \Zb \arrow[r]{}{F} \arrow[d, <-]{}{j_0} & \bl_{\infty \times Z}(\Pb^1 \times X) \arrow[d, <-]{}{s} \\
0 \times Z \arrow[r]{}{i_0} & \{0 \} \times X,
\end{tikzcd}
$$
where $F$ is the strict transform of $\Pb^1 \times Z \hook \Pb^1 \times X$, and denote by $q$ the canonical map $\bl_{\infty \times Z}(\Pb^1 \times X) \to X$. By Lemma \ref{lem:ccgysinregemb1}, $\ch_m$ commutes with $\iota_{\infty !}$ and $s_!$. Moreover,
\begin{equation}\label{eq:dnc}
\iota_! \circ \iota^* \circ F_! = s_! \circ s^* \circ F_!
\end{equation}
since $\iota_! \circ \iota^*$ and $s_! \circ s^*$ correspond to multiplications by the Euler classes of $\Oc(1) \simeq \Oc(\Ec + \bl_{Z}(X))$ and $\Oc(\Ec)$, respectively, and because the strict transforms of $\Pb^1 \times Z$ and $\infty \times X $ do not intersect inside the blowup. Additionally, as $\iota_! \circ \iota^*$ coincides with multiplication by $e(\Oc(\Ec))$, it commutes with $\ch_m$.

For $\alpha \in K^0(Z)$, denote by $\tilde \alpha$ its pullback to $K^0(\Pb^1 \times Z)$. Then,
\begin{align*}
i_!(\ch_m(\alpha)) &= q_! \circ \iota_! \circ i_{\infty !}(\ch_m(\alpha)) \\
&= q_! \circ \iota_! \big(\ch_m(i_{\infty !}(\alpha))\big) \\
&= q_! \circ \iota_! \big(\ch_m(i_{\infty !} \circ j_\infty^* (\tilde\alpha))\big) \\
&= q_! \circ \iota_! \big(\ch_m(\iota^* \circ F_! (\tilde\alpha))\big) \\
&= q_! \circ \iota_! \circ \iota^* \big(\ch_m(F_! (\tilde\alpha))\big) \\
&= q_!  \big(\ch_m(\iota_! \circ \iota^* \circ F_! (\tilde\alpha))\big) \\
&= q_!  \big(\ch_m(s_! \circ s^* \circ F_! (\tilde\alpha))\big) & (\ref{eq:dnc}) \\
&= q_! \circ s_! \big(\ch_m(s^* \circ F_! (\tilde\alpha))\big) \\
&= \ch_m(q_! \circ s_! \circ s^* \circ F_! (\tilde\alpha)) & (q \circ s = \Id) \\
&= \ch_m(q_! \circ s_! \circ j_{0!} (\alpha)) \\
&= \ch_m(i_!(\alpha)),
\end{align*}
as desired.
\end{proof}

We have shown that $\ch_m$ commutes with Gysin pushforwards along projective and quasi-smooth morphisms of Noetherian derived schemes admitting an ample line bundle. We are now ready to prove the main result of this section.

\begin{proof}[Proof of Theorem \ref{thm:cf}]
Let $X$ be a finite-Krull-dimensional Noetherian derived scheme that admits an ample line bundle. We have already shown that the composition $\eta_K \circ \ch_m$ is the identity morphism $K^0(X) \to K^0(X)$, and therefore it suffices to show that the composition $\ch_m \circ \eta_K$ is the identity. As both $\eta_K$ and $\ch_m$ preserve the identity elements and commute with Gysin pushforwards from projective quasi-smooth $X$-schemes, we conclude that
\begin{align*}
\ch_m \circ \eta_K ([V \xto{f} X]) &= \ch_m \circ \eta_K(f_!(1_V)) \\
&= f_!(1_V) \\
&= [V \xto{f} X],
\end{align*} 
and as such cycles generate $\PCob^*(X)$, we are done. 
\end{proof}

\chapter{Universal property of $\PCob^*_S$}\label{ch:univprop}

Universal properties are useful both as tools to construct transformations between  theories, and, more abstractly, as explanations of the ``role'' of a theory. By construction, universal precobordism is obtained from the universal additive bivariant theory by imposing the derived double point cobordism relation, and therefore it is obvious that $\PCob$ is the universal stably oriented additive bivariant theory satisfying this relation. What is not clear is that $\PCob$ has another universal property: namely, it is the universal stably oriented additive bivariant theory that satisfies the projective bundle formula for trivial projective bundles. This universal property was first discovered during a joint project with Ryomei Iwasa, but it has not been explicitly written down before. It served as a motivation for the concept of \emph{pbf-local sheaves} defined in \cite{annala-iwasa}, which we hope will some day yield a useful model for non-$\Ab^1$-invariant higher algebraic cobordism.

Throughout this chapter, $S$ will be a fixed finite-Krull-dimensional Noetherian derived scheme that admits an ample line bundle. We will use the restricted functoriality $\Fc_S$ (Section \ref{ssect:cobbivfunc}) instead of the more general functoriality $\Fc_a$, because it will be important for us that every derived scheme we deal with admits an ample line bundle.

\section{Bivariant universal property of $\PCob^*_S$}

In this section, we prove a universal property for $\PCob^*_S$ as a bivariant theory with functoriality $\Fc_S$.

\subsection{Formal group law of a bivariant theory with a weak projective bundle formula}\label{ssect:pbftofgl}

Here, we show that, in a bivariant theory satisfying a weak version of the projective bundle formula, the Euler class of a tensor product of line bundles may be computed by the means of a formal group law. We begin with the following definition.

\begin{defn}\label{def:pbfthy}
Let $\Bb$ be a stably oriented additive bivariant theory with functoriality $\Fc_S$. Then, $\Bb$ \emph{satisfies the weak projective bundle formula} if, for every $X \in \Cc_S$, and all $n \geq 0$, the map
$$\Bb^\bullet(X)[t]/(t^{n+1}) \to \Bb^\bullet(\Pb^n \times X)$$
is an isomorphism, where $\Bb^\bullet(\Pb^n \times X)$ is considered as an $\Bb^\bullet(X)$ algebra via $\pr_2^*$, and where $t$ maps to $e(\Oc(1)) \in \Bb^\bullet(\Pb^n \times X)$.
\end{defn}

In particular, the $n=0$ case implies that $c_1(\Oc_X) = 0 \in \Bb(X)$. Theories with weak projective bundle formula satisfy a naive algebraic cobordism relation.

\begin{lem}\label{lem:hfib}
Let $X \in \Cc_S$, and denote by $\iota_0$ and $\iota_\infty$ the natural inclusions $X \hook \Pb^1 \times X$ at 0 and $\infty$, respectively. Then, if $\Bb$ satisfies the weak projective bundle formula,
$$i_0^* = i_\infty^*: \Bb^\bullet(\Pb^1 \times X) \to \Bb^\bullet(X).$$
\end{lem}
\begin{proof}
As both $i_0^*$ and $i_\infty^*$ annihilate $e(\Oc(1))$, the claim follows from the fact that they are sections of $\pr_2: \Pb^1 \times X \to X$.
\end{proof}

It follows that Euler classes may be expressed in terms of derived vanishing loci. 

\begin{prop}\label{prop:pbfeuler}
Suppose that $\Bb$ satisfies the weak projective bundle formula. Then, given $X \in \Cc_S$, a vector bundle $E$ on $X$, and a global section $s$ of $E$,
$$e(E) = j_{s!}(1_{V_X(s)}) \in \Bb^\bullet(X),$$
where $j_s$ is the inclusion $V_X(s) \hook X$. 
\end{prop}
\begin{proof}
Considering the Cartesian diagram
$$
\begin{tikzcd}
V_X(0) \arrow[r, hook]{}{j_0} \arrow[d] & X \arrow[d,hook]{}{i_0} \\
V_{\Pb^1 \times X}(x_0 \otimes s) \arrow[r,hook]{}{j} & \Pb^1 \times X \\
V_X(s) \arrow[r,hook]{}{j_s} \arrow[u, hook] & X \arrow[u, hook]{}[swap]{i_\infty}
\end{tikzcd}
$$
we compute that 
\begin{align*}
e(E) &= j_{0!}(1_{V_X(0)}) \\
&= i^*_0\Big(j_!\big(1_{V_{\Pb^1 \times X}(x_0\otimes s)}\big)\Big) & (\text{Proposition \ref{prop:orcohformal}})\\
&= i^*_\infty\Big(j_!\big(1_{V_{\Pb^1 \times X}(x_0\otimes s)}\big)\Big) & (\text{Lemma \ref{lem:hfib}}) \\
&= j_{s!}(1_{V_X(s)}) & (\text{Proposition \ref{prop:orcohformal}})
\end{align*}
as desired.
\end{proof}

Next, we study Euler classes of tensor products. 

\begin{defn}\label{def:pbffgl}
Suppose that $\Bb$ satisfies the weak projective bundle formula. Then, there exists a unique formal power series
$$F_\Bb(x,y) \in \Bb^\bullet(S)[[x,y]]$$
such that
$$e(\Oc(1,1)) = F_\Bb\big(e(\Oc(1,0)), e(\Oc(0,1))\big) \in \Bb^\bullet(\Pb^n_S \times_{S} \Pb^m_S)$$
for all $n,m \geq 0$. 

Similarly, there exists a unique formal power series
$$G_\Bb(x,y) \in \Bb^\bullet(S)[[x,y]]$$
such that 
$$e(\Oc(1,-1)) = G_\Bb\big(e(\Oc(1,0)), e(\Oc(0,1))\big) \in \Bb^\bullet(\Pb^n_S \times_S \Pb^m_S)$$
for all $n,m \geq 0$. 
\end{defn}

\begin{lem}\label{lem:pbffgl}
Suppose that $\Bb$ satisfies the weak projective bundle formula. Then,
\begin{enumerate}
\item for all $X \in \Cc_S$ and globally generated line bundles $\Ls_1$ and $\Ls_2$ on $X$, 
$$e(\Ls_1 \otimes \Ls_2) = F_\Bb\big(e(\Ls_1), e(\Ls_2)\big)$$
and
$$e(\Ls_1 \otimes \Ls_2^\vee) = G_\Bb\big(e(\Ls_1), e(\Ls_2)\big);$$

\item $F_\Bb(x,y)$ is a commutative formal group law;

\item $G_\Bb(F_\Bb(x,y),y) = x$;

\item $G_\Bb(x,y) = F_{\Bb}(x, i_\Bb(y))$, where $i_\Bb(y) \in \Bb^\bullet(S)[[y]]$ is the formal inverse power series of the formal group law $F_\Bb$. 
\end{enumerate}
\end{lem}
\begin{proof}
\begin{enumerate}
\item This follows immediately from the naturality of Euler classes in pullbacks.

\item This follows immediately from the properties of tensor products of line bundles.

\item Indeed, consider $\Id \times \Delta: \Pb^n_S \times_S \Pb^n_S \to \Pb^n_S \times_S \Pb^n_S \times_S \Pb^n_S$. Then,
$$e(\Oc(1,1,-1)) = G_\Bb\Big(F_\Bb\big(e(\Oc(1,0,0)), e(\Oc(0,1,0))\big), e(\Oc(0,0,1)) \Big),$$
and pulling this back along $\Id \times \Delta$, we obtain the formula
$$e(\Oc(1,0)) = G_\Bb\Big(F_\Bb\big(e(\Oc(1,0)), e(\Oc(0,1))\big), e(\Oc(0,1)) \Big).$$
Taking the limit as $n$ goes to $\infty$, we obtain the desired formula. 

\item As $F_\Bb$ is a formal group law, there exists a unique formal power series $G_\Bb(x,y)$ satisfying  the equation $G_\Bb(F_\Bb(x,y),y) = x$, and this is given precisely by $G_\Bb(x,y) = F_{\Bb}(x, i_\Bb(y))$. \qedhere
\end{enumerate}
\end{proof}

In particular, Euler classes of tensor products of arbitrary line bundles (not just of globally generated ones) can be computed via the formal group law.

\begin{prop}\label{prop:pbffgl}
Suppose that $\Bb$ satisfies the weak projective bundle formula. Then, for all $X \in \Cc_S$ and all line bundles $\Ls, \Ls'$ on $X$, the formula 
$$e(\Ls \otimes \Ls') = F_\Bb\big(e(\Ls), e(\Ls')\big) \in \Bb^\bullet(X)$$
holds.
\end{prop}
\begin{proof}
As $X$ admits an ample line bundle, it is possible to find equivalences
$\Ls \simeq \Ls_1 \otimes \Ls^\vee_2$
and
$\Ls' \simeq \Ls'_1 \otimes \Ls'^\vee_2$
where $\Ls_i$ and $\Ls'_i$ are globally generated. As 
$$\Ls \otimes \Ls' \simeq (\Ls_1 \otimes \Ls'_1) \otimes (\Ls_2 \otimes \Ls'_2)^\vee,$$
we compute, using basic properties of formal group laws, that
\begin{align*}
e(\Ls \otimes \Ls') &= G_\Bb\Big(F_\Bb\big(e(\Ls_1), e(\Ls'_1)\big), F_\Bb\big(e(\Ls_2), e(\Ls'_2)\big)\Big) \\
&=  F_\Bb\Big(F_\Bb\big(e(\Ls_1), e(\Ls'_1)\big), i_\Bb\big(F_\Bb\big(e(\Ls_2), e(\Ls'_2)\big)\big)\Big) \\
&= F_\Bb\Big(F_\Bb\big(e(\Ls_1), i_\Bb(e(\Ls_2))\big), F_\Bb\big(e(\Ls'_1), i_\Bb(e(\Ls'_2))\big)\Big) \\
&= F_\Bb(e(\Ls), e(\Ls')),
\end{align*}
as desired.
\end{proof}

\subsection{Derived double point relations for bivariant theories with a formal group law}\label{ssect:fgltoddpt}

Here, we show that bivariant theories, where Euler classes of tensor products of line bundles can be computed via a Formal group law, satisfy a derived analogue of Lee--Pandharipande's double point cobordism relations (i.e., the relations used in Definition \ref{def:univprecob}). Let us begin with the following definition.

\begin{defn}\label{def:fglthy}
Let $\Bb$ be a stably oriented additive bivariant theory with functoriality $\Fc_S$. Then $\Bb$  \emph{has good Euler classes} if
\begin{enumerate}
\item Euler classes of line bundles are nilpotent;

\item for all $X \in \Cc_S$, line bundles $\Ls$ on $X$, and global sections $s$ of $\Ls$, we have  
$$e(\Ls) = i_!(1_{V_X(s)}),$$
where $i$ is the derived regular embedding $V_X(s) \hook X$;
 
\item there exists a formal group law
$$F_\Bb(x, y) = x + y + \sum_{i,j \geq 1} b_{ij} x^i y^j \in \Bb^\bullet(S)[[x,y]]$$
such that, for all $X \in \Cc_S$, and all line bundles $\Ls_1, \Ls_2$ on $X$, we have
$$e(\Ls_1 \otimes \Ls_2) = F_\Bb\big(e(\Ls_1), e(\Ls_2)\big) \in \Bb^\bullet(X).$$
\end{enumerate}
\end{defn}

By the results of the previous subsection, a theory satisfying the weak projective bundle formula has good Euler classes. Let us continue with the following Lemma, which is essentially Lemma 3.3 from \cite{levine-pandharipande}. It is also closely related with Lemma \ref{lem:pbform1}, which gives another formula for the class $[\Pb(\Ls \oplus \Oc) \to X]$.

\begin{lem}\label{lem:fglpbundform}
Suppose that $\Bb$ has good Euler classes. Then, for all $X \in \Cc_S$ and all line bundles $\Ls$ on $X$, we have
$$\pi_!(1_{\Pb(\Ls \oplus \Oc)}) = - \sum_{i,j \geq 1} b_{ij} e(\Ls)^{i-1} e(\Ls^\vee)^{j-1} \in \Bb^\bullet(X),$$
where $\pi$ is the structure morphism $\Pb(\Ls \oplus \Oc) \to X$.
\end{lem}
\begin{proof}
Consider the inclusion $X \simeq \Pb(\Oc) \to \Pb(\Ls \oplus \Oc)$, which has normal bundle $\Ls$. Then, analyzing the derived blowup $W := \bl_{\infty \times X}(\Pb^1 \times \Pb(\Ls \oplus \Oc))$, we note that the fiber of the natural map $W \to \Pb^1_S$ over
\begin{enumerate}
\item $0$, is equivalent to $\Pb(\Ls \oplus \Oc)$;
\item $\infty$, is the sum of the strict transform $D$ of $\infty \times \Pb(\Ls \oplus \Oc) \hook \Pb^1 \times \Pb(\Ls \oplus \Oc)$ and the exceptional divisor $\Pb(\Ls \oplus \Oc) \simeq \Ec \subset W$, whose intersection is $X$.
\end{enumerate}
Moreover, the restrictions of $\Oc(D)$ and $\Oc(\Ec)$ to $X$ are equivalent to $\Ls^\vee$ and $\Ls$, respectively. Hence, if we denote by $q$ the natural map $W \to X$, we have 
\begin{align*}
\pi_!(1_{\Pb(\Ls \oplus \Oc)}) &= q_!(e(\Oc(1))) \\
&= q_!\Big( F_\Bb \big( e(\Oc(D)), e(\Oc(\Ec)) \big) \Big)
\end{align*} 
from which the desired equation follows immediately.
\end{proof}

Next, we prove that a theory with good Euler classes satisfies the derived double point relations. Since a theory with good Euler classes is by assumption stably oriented and additive, it is a target of a unique orientation preserving Grothendieck transformation from the universal additive theory $\Ab^*_{\Fc_S}$. Hence, it is sensible to use notation $[V \to X]$ to describe an element of $\Bb(X \to Y)$.

\begin{prop}\label{prop:fglthyddpt}
Suppose that $\Bb$ has good Euler classes, and that $X \to Y$ is a morphism in $\Cc_S$. Then, for any projective morphism $\pi: W \to \Pb^1 \times X$ such that the composition $\rho: W \to \Pb^1 \times Y$ is quasi-smooth, the equation
$$[W_0 \to X] = [D_1 \to X] + [D_2 \to X] - [\Pb_{D_1 \cap D_2}(\Oc(D_1) \oplus \Oc) \to X] \in \Bb(X \to Y)$$
holds, where
\begin{enumerate}
\item $W_0$ is the fiber of $W \to \Pb^1 \times X$ over $0 \times X$;

\item the fiber $W_\infty$ of $W \to \Pb^1 \times X$ over $\infty \times X$ is the sum of virtual Cartier divisors $D_1$ and $D_2$ on $W$.
\end{enumerate}
\end{prop}
\begin{proof}
Suppose first that $W = X = Y$. Then, denoting by $j$ the derived regular embedding $D_1 \cap D_2 \hook W$, we compute that 
\begin{align*}
&[W_0 \to W] \\
&= e(\Oc(1)) \\
&= e(\Oc(D_1 + D_2)) \\
&= [D_1 \to W] + [D_2 \to W] \\
&+ \sum_{i,j \geq 1} b_{ij} e(\Oc(D_1))^{i-1} \bullet e(\Oc(D_2))^{j-1} \bullet [D_1 \cap D_2 \to W] \\
&= [D_1 \to W] + [D_2 \to W] + j_!\bigg( \sum_{i,j \geq 1} e(\Oc(D_1))^{i-1} \bullet e(\Oc(D_2))^{j-1} \bigg) \\
&= [D_1 \to W] + [D_2 \to W] + j_!\bigg( \sum_{i,j \geq 1} e(\Oc(D_1))^{i-1} \bullet e(\Oc(-D_1))^{j-1} \bigg) \\
&= [D_1 \to W] + [D_2 \to W] - j_!\big([\Pb_{D_1 \cap D_2}(\Oc(D_1) \oplus \Oc) \to D_1 \cap D_2]\big) \\
&= [D_1 \to W] + [D_2 \to W] - [\Pb_{D_1 \cap D_2}(\Oc(D_1) \oplus \Oc) \to W],
\end{align*}
where, in the fifth equation, we have used the fact that $\Oc(1) \simeq \Oc(D_1 + D_2)$ restricts to trivial line bundle on $W_\infty$.

The general case follows from the above by applying the transformation $\pi_*(- \bullet \theta(\rho))$ to the obtained equation.
\end{proof}

\subsection{Universal property of $\PCob^*_S$}

Here, we show that $\PCob^*_S$ is the universal stably oriented additive bivariant theory that satisfies projective bundle formula.

\begin{thm}\label{thm:precobunivprop}
Let $S$ be a finite-Krull-dimensional Noetherian derived scheme that admits an ample line bundle. Then,
\begin{enumerate}
\item $\PCob^*_S$ is the universal bivariant theory that satisfies the weak projective bundle formula: given another such a bivariant theory $\Bb$, then there exists a unique orientation preserving Grothendieck transformation $\PCob^*_S \to \Bb$;

\item $\PCob^*_S$ is the universal bivariant theory that has good Euler classes: given another such a bivariant theory $\Bb$, then there exists a unique orientation preserving Grothendieck transformation $\PCob^*_S \to \Bb$.
\end{enumerate}
\end{thm}
\begin{proof}
The first claim follows from the second, as a theory that satisfies the weak projective bundle formula has good Euler classes. Since a theory $\Bb$ that has good Euler classes is by hypothesis stably oriented and additive, it admits a unique orientation preserving Grothendieck transformation from the universal additive theory $\eta: \Ab^*_{\Fc_S} \to \Bb$. By Proposition \ref{prop:fglthyddpt}, $\eta$ factors through $\PCob^*_S$, so the claim follows. 
\end{proof}

\section{Homological and cohomological universal properties}

Here we record universal properties for the associated homology and cohomology theories of the universal precobordism $\PCob^*_S$. This makes it easier to study the associated theories, for they no longer have to be considered as a part of a bivariant theory for the universal property to hold.

\subsection{Universal property of the universal precobordism cohomology}

Inspired by Proposition \ref{prop:orcohformal}, we make the following definition.

\begin{defn}\label{def:orcoh}
An \emph{$\Fc_S$-oriented cohomology theory} $\Hb^\bullet$ consists of
\begin{enumerate}
\item a contravariant functor $\Hb^\bullet$ from $\Cc_S$ to the category of (discrete) rings; the contravariant functoriality along a map $f: X \to Y$ is denoted by $f^*$, and it is called the \emph{pullback} along $f$; 

\item additive \emph{Gysin pushforwards} $f_!: \Hb^\bullet(X) \to \Hb^\bullet(Y)$ along projective and quasi-smooth morphisms $f: X \to Y$;
\end{enumerate}
satisfying
\begin{enumerate}
\item \emph{additivity}: the equality 
$$1_{X_1 \sqcup X_2} = \iota_{1!}(1_{X_1}) + \iota_{2!}(1_{X_2})$$
holds in $\Hb^\bullet(X_1 \sqcup X_2)$, where $\iota_i$ are the canonical inclusions $X_i \to X_1 \sqcup X_2$\footnote{In fact, assuming the other hypotheses of a $\Fc_S$-oriented cohomology theory, this is equivalent to requiring the map 
$$
\begin{bmatrix}
\iota_1^* \\
\iota_2^*
\end{bmatrix}
: \Hb^\bullet(X_1 \sqcup X_2) \to \Hb^\bullet(X_1) \times \Hb^\bullet(X_2)$$
to be an isomorphism or rings. The other direction proven similarly to the proof of Proposition \ref{prop:add}, and the other direction is proven by an easy computation.
};

\item \emph{covariant functoriality}: for all projective and quasi-smooth $f: X \to Y$ and $g: Y \to Z$, 
$$(g \circ f)_! = g_! \circ f_!;$$
moreover, $\Id_! = \Id$;

\item \emph{push-pull formula}: if the square
$$
\begin{tikzcd}
X' \arrow[d]{}{f'} \arrow[r]{}{g'} & Y' \arrow[d]{}{f} \\
X  \arrow[r]{}{g} & Y
\end{tikzcd}
$$
is Cartesian and $f$ is both quasi-smooth and projective, then 
$$g^* \circ f_! = f'_! \circ g'^*;$$

\item \emph{projection formula}: if $f: X \to Y$ is projective and quasi-smooth, then, for all $\alpha \in \Hb^\bullet(Y)$ and $\beta \in \Hb^\bullet(X)$, the equality
$$f_! \big(f^*(\alpha) \bullet \beta \big) = \alpha \bullet f_!(\beta)$$
holds.
\end{enumerate}
A \emph{transformation} of $\Fc_S$-oriented cohomology theories is a natural transformation $\eta: \Hb_1^\bullet \to \Hb_2^\bullet$ such that, for all quasi-smooth and projective morphisms $f: X \to Y$, the square
$$
\begin{tikzcd}
\Hb_1^\bullet(X) \arrow[d]{}{f_!} \arrow[r]{}{\eta_X} & \Hb_2^\bullet(X) \arrow[d]{}{f_!} \\
\Hb_1^\bullet(Y) \arrow[r]{}{\eta_Y} & \Hb_2^\bullet(Y)
\end{tikzcd}
$$
commutes.
\end{defn}

The following preliminary result gives an universal property for the associated cohomology theory $\Ab^*_{\Fc_S}$ of the universal additive bivariant theory with functoriality $\Fc_S$ (also denoted by $\Ab^*_{\Fc_S}$).

\begin{lem}
Let $\Hb^\bullet$ be a $\Fc_S$-oriented cohomology theory. Then, there exists a unique transformation of $\Fc_S$-oriented cohomology theories $\eta: \Ab^*_{\Fc_S} \to \Hb^\bullet$.
\end{lem}
\begin{proof}
Since a transformation of $\Fc_S$-oriented cohomology theories preserves identities and commutes with Gysin pushforwards, $\eta$, if it exists, must satisfy the formula
$$\eta([V \xto{f} X]) = f_!(1_V).$$
In order to conclude that this describes a well-defined transformation, we must show that, given projective and quasi-smooth map $f: V_1 \sqcup V_2 \to X$, then
$$f_!(1_{V_1 \sqcup V_2}) = f_{1!}(1_{V_1}) + f_{2!}(1_{V_2}),$$
where $f_i: V_i \to X$ are the restrictions of $f$. However, this follows immediately from the formula
$$1_{V_1 \sqcup V_2} = \iota_{1!}(1_{V_1}) + \iota_{2!}(1_{V_2}) \in \Hb^\bullet(V_1 \sqcup V_2)$$
and the functoriality and additivity of Gysin pushforwards.
\end{proof}

Since the definitions of a bivariant theory satisfying the weak projective bundle formula and of a bivariant theory having good Euler classes only mention the associated cohomology theory, it is easy to export these definitions to the context of $\Fc_S$-oriented cohomology theories. 

\begin{defn}\label{def:orcohpbffgl}
Let $\Hb^\bullet$ be a $\Fc_S$-oriented cohomology theory. Then, $\Hb^\bullet$ \emph{satisfies the weak projective bundle formula} if the map
$$\Hb^\bullet(X)[t]/(t^{n+1}) \to \Hb^\bullet(\Pb^n \times X)$$
is an isomorphism of rings, where $\Hb^\bullet(\Pb^n \times X)$ is considered as an $\Hb^\bullet(X)$-algebra via $\pr^*_2$, and where $t$ maps to $e(\Oc(1)) \in \Hb^\bullet(\Pb^n \times X)$.  The theory $\Hb^\bullet$ \emph{has good Euler classes} if it satisfies the conditions 1., 2., and 3. of Definition \ref{def:fglthy}.
\end{defn}

Following the proof of Propositions \ref{prop:pbfeuler} and \ref{prop:pbffgl}, one arrives at the following result. 

\begin{lem}\label{lem:orcohpbffgl}
A $\Fc_S$-oriented cohomology theory that satisfies the weak projective bundle formula has good Euler classes. \qed 
\end{lem}

Following the proof of Proposition \ref{prop:fglthyddpt}, we arrive at the following result.

\begin{lem}
Let $\Hb^\bullet$ be a $\Fc_S$-oriented cohomology theory that has good Euler classes. Then, for any projective quasi-smooth morphism $\pi: W \to \Pb^1 \times X$, the equation
\begin{align*}
&[W_0 \to X] \\
&= [D_1 \to X] + [D_2 \to X] - [\Pb_{D_1 \cap D_2}(\Oc(D_1) \oplus \Oc) \to X] \in \Hb^\bullet(X \to Y)
\end{align*}
holds, where
\begin{enumerate}
\item $W_0$ is the fiber of $W \to \Pb^1 \times X$ over $0 \times X$;

\item the fiber $W_\infty$ of $W \to \Pb^1 \times X$ over $\infty \times X$ is the sum of virtual Cartier divisors $D_1$ and $D_2$ on $W$. \qed
\end{enumerate}
\end{lem}

And hence, we obtain the following universal properties for $\Hb^\bullet$.

\begin{thm}\label{thm:precobcohunivprop}
Let $S$ be a finite-Krull-dimensional Noetherian derived scheme that admits an ample line bundle. Then,
\begin{enumerate}
\item $\PCob^*_S$ is the universal $\Fc_S$-oriented cohomology theory that satisfies the weak projective bundle formula: given another such a $\Fc_S$-oriented cohomology theory $\Hb^\bullet$, there exists a unique transformation $\PCob^*_S \to \Hb^\bullet$ of $\Fc_S$-oriented cohomology theories;

\item $\PCob^*_S$ is the universal $\Fc_S$-oriented cohomology theory that has good Euler classes: given another such a $\Fc_S$-oriented cohomology theory $\Hb^\bullet$, there exists a unique transformation $\PCob^*_S \to \Hb^\bullet$ of $\Fc_S$-oriented cohomology theories. \qed
\end{enumerate}
\end{thm}

\subsection{Universal property of the universal precobordism homology}

Here, we provide $\PCob^S_{*}$ with universal properties analogous to the bivariant and the cohomological universal properties studied above. Unfortunately, since the bivariant universal properties are rather cohomological in nature, it takes more effort to translate them into homological universal properties than it did to translate them into cohomological ones.

Inspired by Proposition \ref{prop:orhomformal} (cf. \cite{levine-morel} Definition 5.1.3), we make the following definition.

\begin{defn}\label{def:orhom}
Let us denote by $\Cc'_S$ the subcategory of $\Cc_S$ whose objects are those of $\Cc_S$, and whose morphisms are the projective morphisms in $\Cc_S$. Then, a $\Fc_S$-oriented Borel--Moore homology theory consists of
\begin{enumerate}
\item a covariant functor $\Hb_\bullet$ from $\Cc'_S$ to the category of (discrete) Abelian groups; the covariant functoriality along a projective morphism $f: X \to Y$ is denoted by $f_*$, and it is called the \emph{pushforward} along $f$;

\item for each quasi-smooth morphism $f: X \to Y$, a \emph{Gysin pullback} homomorphism 
$$f^!: \Hb_\bullet(Y) \to \Hb_\bullet(X);$$

\item an element $1 \in \Hb_\bullet(S)$, and for each $X,Y \in \Cc_S$ a bi-additive pairing
$$\times: \Hb_\bullet(X) \times \Hb_\bullet(Y) \to \Hb_\bullet(X \times_S Y),$$
referred to as \emph{cross product}, which is associative, admits $1$ as a unit, and is commutative in the sense that the square 
$$
\begin{tikzcd}
\Hb_\bullet(X) \otimes \Hb_\bullet(Y) \arrow[r]{}{\times} \arrow[d]{}{\sigma} & \Hb_\bullet(X \times_{S} Y) \arrow[d]{}{\iota_*} \\
\Hb_\bullet(Y) \otimes \Hb_\bullet(X) \arrow[r]{}{\times} & \Hb_\bullet(Y \times_{S} X)
\end{tikzcd}
$$
commutes, where $\sigma$ is the canonical involution $A \otimes B \to B \otimes A$ and $\iota$ is the canonical involution $X \times_S Y \to Y \times_S X$; for each quasi-smooth $X \in \Cc_S$, we will denote by $1_X$ the \emph{fundamental class} $\pi_X^!(1)$, where $\pi_X$ is the structure morphism $X \to S$;
\end{enumerate}
satisfying
\begin{enumerate}
\item \emph{additivity}: the natural map
$$
\begin{bmatrix}
\iota_{1*} & \iota_{2*}
\end{bmatrix}
: \Hb_\bullet(X_1) \oplus \Hb_\bullet(X_2) \to \Hb_\bullet(X_1 \sqcup X_2)
$$
is an isomorphism, where $\iota_i$ are the canonical inclusions $X_i \to X_1 \sqcup X_2$;

\item \emph{contravariant functoriality}: if $f: X \to Y$ and $g: Y \to Z$ are quasi-smooth, then 
$$(g \circ f)^! = f^! \circ g^!;$$
moreover, $\Id^! = \Id$;

\item \emph{functoriality of cross products}: if $f: X \to X'$ and $g: Y \to Y'$ are projective, then, for all $\alpha \in \Hb_\bullet(X)$ and $\beta \in \Hb_\bullet(Y)$,
$$f_*(\alpha) \times g_*(\beta) = (f \times g)_*(\alpha \times \beta);$$
if $f: X \to X'$ and $g: Y \to Y'$ are quasi-smooth, then, for all $\alpha \in \Hb_\bullet(X')$ and $\beta \in \Hb_\bullet(Y')$,
$$f^!(\alpha) \times g^!(\beta) = (f \times g)^!(\alpha \times \beta);$$

\item \emph{push-pull formula}: if
$$
\begin{tikzcd}
X' \arrow[d]{}{f'} \arrow[r]{}{g'} & Y' \arrow[d]{}{f} \\
X \arrow[r]{}{g} & Y
\end{tikzcd}
$$
is Cartesian, $f$ is projective, and $g$ is quasi-smooth, then
$$g^! \circ f_* = f'_* \circ g'^!.$$
\end{enumerate}
A \emph{transformation} of $\Fc_S$-oriented Borel--Moore homology theories is a natural transformation of functors from $\Cc'_S$ to (discrete) Abelian groups that commutes with Gysin pullbacks.
\end{defn}

The following preliminary result gives an universal property for the associated homology theory $\Ab^{\Fc_S}_*$ of the universal additive bivariant theory with functoriality $\Fc_S$.

\begin{lem}\label{lem:univorhom}
Let $\Hb_\bullet$ be a $\Fc_S$-oriented Borel--Moore homology theory. Then, there exists a unique transformation of $\Fc_S$-oriented Borel--Moore homology theories $\eta: \Ab_*^{\Fc_S} \to \Hb_\bullet$.
\end{lem}
\begin{proof}
Since a transformation of $\Fc_S$-oriented Borel--Moore homology theories preserves fundamental classes and commutes with pushforwards, $\eta$, if it exists, must satisfy the formula
$$\eta([V \xto{f} X]) = f_*(1_V).$$
In order to conclude that this describes a well-defined transformation, we must show that, given projective map $f: V_1 \sqcup V_2 \to X$ with a quasi-smooth source, then
$$f_*(1_{V_1 \sqcup V_2}) = f_{1*}(1_{V_1}) + f_{2*}(1_{V_2}),$$
where $f_i: V_i \to X$ are the restrictions of $f$. By pushing forward, this would follow immediately from the formula
$$1_{V_1 \sqcup V_2} = \iota_{1*}(1_{V_1}) + \iota_{2*}(1_{V_2}) \in \Hb_\bullet(V_1 \sqcup V_2),$$
where $\iota_i$ are the canonical inclusions $V_i \to V_1 \sqcup V_2$. We will prove this formula.

By hypothesis, the map 
$$
\iota_2 :=
\begin{bmatrix}
\iota_{1*} & \iota_{2*}
\end{bmatrix}
: \Hb_\bullet(V_1) \oplus \Hb_\bullet(V_2) \to \Hb_\bullet(V_1 \sqcup V_2)
$$
is an isomorphism, and from the push-pull formula it follows that its inverse is given by
$$
\iota^! :=
\begin{bmatrix}
\iota^!_{1} \\
\iota^!_{2}
\end{bmatrix}
: \Hb_\bullet(V_1 \sqcup V_2) \to \Hb_\bullet(V_1) \oplus \Hb_\bullet(V_2).
$$
As 
\begin{align*}
1_{V_1 \sqcup V_2} &= \iota_*(\iota^!(1_{V_1 \sqcup V_2}))\\
&= \iota_{1*}(1_{V_1}) + \iota_{2*}(1_{V_2}), & (\text{push-pull formula})
\end{align*}
the claim follows.
\end{proof}

The conditions of Definition \ref{def:orhom} imply that a $\Fc_S$-oriented Borel--Moore homology theory $\Hb_\bullet$ is, in a sense, linear over the commutative ring $\Hb_\bullet(S)$.

\begin{prop}\label{prop:orhomlinearstr}
Let $\Hb_\bullet$ be a $\Fc_S$-oriented Borel--Moore homology theory. Then,
\begin{enumerate}
\item the cross product provides each $\Hb_\bullet(X)$ with the structure of a $\Hb_\bullet(S)$-module;

\item pushforwards and Gysin pullbacks are $\Hb_\bullet(S)$-linear;

\item the cross product is $\Hb_\bullet(S)$-bilinear;

\item for all quasi-smooth $X \in \Cc_S$, and all $\alpha \in \Hb_\bullet(Y)$, we have
$$\pr_2^!(\alpha) = 1_X \times \alpha \in \Hb_\bullet(X \times_S Y).$$
\end{enumerate}
\end{prop}
\begin{proof}
The claims 1., 2., and 3. follow immediately from the definition. The fourth claim follows from the definition after observing that $\pr_2 = \pi_X \times \Id_Y$, where $\pi_X: X \to S$ is the structure morphism.
\end{proof}

\begin{defn}\label{def:orhomeulerclass}
If $\Hb_\bullet$ is a $\Fc_S$-oriented Borel--Moore homology theory, $X \in \Cc_S$, and $E$ is a vector bundle on $X$, then we define the \emph{Euler class operator} as
$$e(E) := i_* \circ i^!: \Hb_\bullet(X) \to \Hb_\bullet(X),$$
where $i$ is the derived regular embedding $V_X(0_E) \hook X$, $0_E$ being the zero-section of $E$.
\end{defn}

Euler class operators are natural in the following sense.

\begin{prop}\label{prop:orhomeulerclassprops}
Let $\Hb_\bullet$ be a $\Fc_S$-oriented Borel--Moore homology theory. Then,
\begin{enumerate}
\item if $f: X \to Y$ is a quasi-smooth morphism in $\Cc_S$ and $E$ is a vector bundle on $Y$, then, for all $\alpha \in \Hb_\bullet(Y)$,
$$f^!\big(e(E) (\alpha) \big) = e(f^* E) (f^!(\alpha));$$
\item if $f: X \to Y$ is a projective morphism in $\Cc_S$ and $E$ is a vector bundle on $Y$, then, for all $\beta \in \Hb_\bullet(X)$,
$$f_*\big( e(f^*E) (\beta)\big)  = e(E)(f_*(\beta)).$$
\end{enumerate}
\end{prop}
\begin{proof}
Consider the Cartesian diagram
$$
\begin{tikzcd}
V_X(0_{f^*E}) \arrow[r]{}{i'} \arrow[d]{}{f'} & X \arrow[d]{}{f} \\
V_Y(0_E) \arrow[r]{}{i} & Y.
\end{tikzcd}
$$
Then,
\begin{enumerate}
\item if $f$ is quasi-smooth, we have that for all $\alpha \in \Hb_\bullet(Y)$, 
\begin{align*}
f^!\big(e(E)(\alpha)\big) &= f^!\big(i_*(i^!(\alpha))\big) \\
&= i'_*\big( f'^!(i^!(\alpha)) \big)  & (\text{push-pull formula})  \\
&= i'_*\big( i'^!(f^!(\alpha)) \big) \\
&= e(f^*E)(f^!(\alpha)),
\end{align*}
as desired;

\item if $f$ is projective, we have that for all $\beta \in \Hb_\bullet(X)$,
\begin{align*}
f_*\big(e(f^*E)(\beta)\big) &= f_*\big(i'_*(i'^!(\beta))\big) \\
&= i_*\big(f'_*(i'^!(\beta))\big) \\
&= i_*\big(i^!(f_*(\beta))\big) & (\text{push-pull formula}) \\
&= e(E)(f_*(\beta))
\end{align*}
as desired. \qedhere
\end{enumerate}
\end{proof}

Moreover, they commute with each other.

\begin{prop}\label{prop:orhomeulerclasscomm}
Let $\Hb_\bullet$ be a $\Fc_S$-oriented Borel--Moore homology theory. Then, for each $X \in \Cc_S$, and vector bundles $E$ and $F$ on $X$, we have
$$e(E) \circ e(F) = e(E \oplus F) = e(F) \circ e(E).$$
\end{prop}
\begin{proof}
Considering the Cartesian square
$$
\begin{tikzcd}
V_X(0_{E \oplus F}) \arrow[r, hook]{}{j'} \arrow[d, hook]{}{i'} & V_{X}(0_F) \arrow[d, hook]{}{i} \\
V_X(0_E) \arrow[r, hook]{}{j} & X
\end{tikzcd}
$$
we compute that
\begin{align*}
e(E) \circ e(F) &= j_* \circ j^! \circ i_* \circ i^! \\
&= j_* \circ i'_* \circ j'^! \circ i^! & (\text{push-pull formula})\\
&= (j \circ i')_* \circ (j \circ i')^! \\
&= e(E \oplus F)
\end{align*}
from which the claim follows, as direct summation of vector bundles is commutative.
\end{proof}

Next, we will closely follow Sections \ref{ssect:pbftofgl} and \ref{ssect:fgltoddpt} in order to produce the desired homological universal properties. 

\begin{defn}\label{def:orhompbfthy}
A $\Fc_S$-oriented Borel--Moore homology theory \emph{satisfies the weak projective bundle formula} if, for all $X \in \Cc_S$, and all $n \geq 0$, the map
$$e(\Oc(1))^i\big(\pr_2^!(-)\big):\bigoplus_{i = 0}^{n} \Hb_\bullet(X) \to \Hb_\bullet(\Pb^n \times X)$$ 
is an isomorphism, and $e(\Oc(1))^{n+1}$ is the zero operator.
\end{defn}

\begin{lem}\label{lem:orhomhfib}
Suppose $\Hb_\bullet$ satisfies the weak projective bundle formula, let $X \in \Cc_S$, and let $i_0$ and $i_\infty$ be the inclusions $X \hook \Pb^1 \times X$ at $0$ and $\infty$, respectively. Then,
$$i_0^! = i_\infty^!: \Hb_\bullet(\Pb^1 \times X) \to \Hb_\bullet(X).$$
\end{lem}
\begin{proof}
As $e(\Oc)$ is the zero operator on $\Hb_\bullet(X)$, both $i_0^!$ and $i_\infty^!$ annihilate the image of $e(\Oc(1))^i\big(\pr_2^!(-)\big)$ by Proposition \ref{prop:orhomeulerclassprops}. The claim then follows from the fact that both $i_0$ and $i_\infty$ are sections of $\pr_2$.
\end{proof}

It follows that Euler classes may be expressed in terms of derived vanishing loci. 

\begin{prop}\label{prop:orhompbfeuler}
Suppose that $\Hb_\bullet$ satisfies the weak projective bundle formula. Then, given $X \in \Cc_S$, a vector bundle $E$ on $X$, and a global section $s$ of $E$,
$$e(E) = j_{s*} \circ j_s^!,$$
where $j_s$ is the inclusion $V_X(s) \hook X$. 
\end{prop}
\begin{proof}
Considering the Cartesian diagram
$$
\begin{tikzcd}
V_X(0) \arrow[r, hook]{}{j_0} \arrow[d]{}{i'_0} & X \arrow[d,hook]{}{i_0} \\
V_{\Pb^1 \times X}(x_0 \otimes s) \arrow[r,hook]{}{j} & \Pb^1 \times X \\
V_X(s) \arrow[r,hook]{}{j_s} \arrow[u, hook]{}[swap]{i'_\infty} & X, \arrow[u, hook]{}[swap]{i_\infty}
\end{tikzcd}
$$
we compute that 
\begin{align*}
e(E) &= j_{0*} \circ j_0^! \\
&= j_{0*} \circ j_0^! \circ i_0^! \circ \pr_2^! \\
&= j_{0*} \circ i'^!_0 \circ j^! \circ \pr_2^! \\
&= i_0^! \circ j_* \circ j^! \circ \pr_2^! & (\text{push-pull formula})\\
&= i_\infty^! \circ j_* \circ j^! \circ \pr_2^! & (\text{Lemma \ref{lem:orhomhfib}})\\
&= j_{s*} \circ i'^!_\infty \circ j^! \circ \pr_2^! & (\text{push-pull formula})\\
&= j_{s*} \circ j^!_s \circ i_\infty^! \circ \pr_2^! \\
&= j_{s*} \circ j^!_s,
\end{align*}
as desired.
\end{proof}

Next, we study Euler class operators of tensor products. 

\begin{defn}\label{def:orhompbffgl}
Suppose that $\Hb_\bullet$ satisfies the weak projective bundle formula. Then, there exists a unique formal power series
$$F_\Hb(x,y) = \sum_{i,j \geq 0} b_{ij} x^i y^j \in \Hb_\bullet(S)[[x,y]]$$
such that
$$e(\Oc(1,1))\big(1_{\Pb^n_S \times_S \Pb^m_S}\big) = F_\Hb\big(e(\Oc(1,0)), e(\Oc(0,1))\big)\big(1_{\Pb^n_S \times_S \Pb^m_S}\big)$$
in $\Hb_\bullet(\Pb^n_S \times_S \Pb^m_S)$ for all $n,m \geq 0$. 

Similarly, there exists a unique formal power series
$$G_\Hb(x,y) \in \Hb_\bullet(S)[[x,y]]$$
such that 
$$e(\Oc(1,-1))\big(1_{\Pb^n_S \times_S \Pb^m_S}\big) = G_\Hb\big(e(\Oc(1,0)), e(\Oc(0,1))\big)\big(1_{\Pb^n_S \times_S \Pb^m_S}\big)$$
in $\Hb_\bullet(\Pb^n_S \times_S \Pb^m_S)$
for all $n,m \geq 0$. 
\end{defn}

Before studying the above formal power series further, we make the following observation.

\begin{lem}
Let $\Hb_\bullet$ be a $\Fc_S$-oriented Borel--Moore homology theory, let $X, Y \in \Cc_S$, let $E$ be a vector bundle on $X$, and let $\alpha \in \Hb_\bullet(X)$ and $\beta \in \Hb_\bullet(Y)$. Then,
$$e(\pr_1^* E)\big( \alpha \times \beta \big) = e(E)(\alpha) \times \beta \in \Hb_\bullet \big( X \times_S Y \big)$$
and
$$e(\pr_2^* E)\big( \beta \times \alpha \big) = \beta \times e(E)(\alpha) \in \Hb_\bullet \big( Y \times_S X \big).$$
\end{lem}
\begin{proof}
Since
$$
\begin{tikzcd}
V_{X \times_S Y}(0_{\pr_1^*E}) \simeq V_X(0_E) \times_S Y \arrow[d] \arrow[r,hook]{}{i \times \Id_Y} & X \times_{S} Y \arrow[d]{}{\pr_1} \\
V_X(0_E) \arrow[r,hook]{}{i} & X
\end{tikzcd}
$$
is Cartesian, we have
\begin{align*}
e(\pr_1^*E)(\alpha \times \beta) &= (i \times \Id_Y)_* \big((i \times \Id_Y)^!(\alpha \times \beta)\big) \\
&= i_*\big(i^!(\alpha)\big) \times \beta \\
&= e(E)(\alpha) \times \beta,
\end{align*}
proving the first claim. The second claim is proven in a similar fashion.
\end{proof}

Next, we prove that the power series govern the Euler class operators of tensor product.

\begin{lem}\label{lem:orhompbffgl}
Suppose that $\Hb_\bullet$ satisfies the weak projective bundle formula, that $X \in \Cc_S$, and that $\Ls_1$ and $\Ls_2$ are globally generated line bundles on $X$. Then, 
$$e(\Ls_1 \otimes \Ls_2) = F_\Hb\big(e(\Ls_1), e(\Ls_2)\big)$$
and
$$e(\Ls_1 \otimes \Ls^\vee_2) = G_\Hb\big(e(\Ls_1), e(\Ls_2)\big).$$
\end{lem}
\begin{proof}
The proof of the second claim is similar to that of the first one, so we will only provide a proof for the first claim. Let $f: X \to \Pb^n_S \times_S \Pb^n_S$ be such that $f^* \Oc(1,0) \simeq \Ls_1$ and $f^* \Oc(0,1) \simeq \Ls_2$. Then, denoting by $\Gamma_f$ the graph map $X \to \Pb^n \times \Pb^m \times X$, we compute that, for all $\alpha \in \Hb_\bullet(X)$,
\begin{align*}
e(\Ls_1 \otimes \Ls_2)(\alpha) &= \Gamma_f^!\big( e(\pr_1^* \Oc(1,1))\big(1_{\Pb^n_S \times_S \Pb^m_S} \times \alpha\big) \big) \\
&= \Gamma_f^!\big(e(\Oc(1,1))\big(1_{\Pb^n_S \times_S \Pb^m_S}\big) \times \alpha \big) \\
&= \sum_{i,j \geq 0} b_{ij} \Gamma_f^!\Big( e(\Oc(1,0))^i \circ e(\Oc(0,1))^j \big(1_{\Pb^n_S \times_S \Pb^m_S}\big) \times \alpha \Big) \\
&= \sum_{i,j \geq 0} b_{ij} \Gamma_f^!\Big( e(\pr_1^*\Oc(1,0))^i \circ e(\pr_1^*\Oc(0,1))^j \big(1_{\Pb^n_S \times_S \Pb^m_S} \times \alpha \big) \Big) \\
&= \sum_{i,j \geq 0} b_{ij}  e(\Ls_1)^i \circ e(\Ls_2)^j ( \alpha ), 
\end{align*}
as desired.
\end{proof}

Hence, we obtain the following analogue of Lemma \ref{lem:pbffgl}.

\begin{lem}\label{lem:orhompbffgl2}
Suppose that $\Hb_\bullet$ satisfies the weak projective bundle formula. Then,
\begin{enumerate}
\item $F_\Hb(x,y)$ is a commutative formal group law;

\item $G_\Hb(F_\Hb(x,y),y) = x$;

\item $G_\Hb(x,y) = F_{\Hb}(x, i_\Hb(y))$, where $i_\Hb(y) \in \Hb_\bullet(S)[[y]]$ is the formal inverse power series of the formal group law $F_\Hb$. 
\end{enumerate}
\end{lem}
\begin{proof}
Proof is the same as that of Lemma \ref{lem:pbffgl}.
\end{proof}

In particular, the following result is proven in exactly the same fashion as Proposition \ref{prop:pbffgl}.

\begin{prop}\label{prop:orhompbffgl}
Suppose that $\Hb_\bullet$ satisfies the weak projective bundle formula. Then, for all $X \in \Cc_S$ and all line bundles $\Ls, \Ls'$ on $X$,  
$$e(\Ls \otimes \Ls') = F_\Hb(e(\Ls), e(\Ls'))$$
as operators on $\Hb_\bullet(X)$. \qed
\end{prop}

This result leads naturally to the following definition.

\begin{defn}\label{def:orhomfglthy}
Let $\Hb_\bullet$ be a $\Fc_S$-oriented Borel--Moore homology theory. Then, $\Hb_\bullet$  \emph{has good Euler classes} if
\begin{enumerate}
\item Euler class operators of line bundles are nilpotent;

\item for all $X \in \Cc_S$, line bundles $\Ls$ on $X$, and global sections $s$ of $\Ls$, we have  
$$e(\Ls) = i_* \circ i^!,$$
where $i$ is the derived regular embedding $V_X(s) \hook X$;
 
\item there exists a formal group law
$$F_\Bb(x, y) = x + y + \sum_{i,j \geq 1} b_{ij} x^i y^j \in \Hb_\bullet(S)[[x,y]]$$
such that, for all $X \in \Cc_S$, and all line bundles $\Ls_1, \Ls_2$ on $X$, we have
$$e(\Ls_1 \otimes \Ls_2) = F_\Bb\big(e(\Ls_1), e(\Ls_2)\big)$$
as operators on $\Hb_\bullet(X)$.
\end{enumerate}
\end{defn}

A theory that satisfies the weak projective bundle formula has good Euler classes.

\begin{lem}\label{lem:orhomfglpbundform}
Suppose that $\Hb_\bullet$ has good Euler classes. Then, for all quasi-smooth $X \in \Cc_S$ and all line bundles $\Ls$ on $X$, we have
$$\pi_*(1_{\Pb(\Ls \oplus \Oc)}) = - \sum_{i,j \geq 1} b_{ij} e(\Ls)^{i-1} \circ e(\Ls^\vee)^{j-1}(1_X) \in \Hb_\bullet(X),$$
where $\pi$ is the structure morphism $\Pb(\Ls \oplus \Oc) \to X$.
\end{lem}
\begin{proof}
Consider the inclusion $X \simeq \Pb(\Oc) \to \Pb(\Ls \oplus \Oc)$, which has normal bundle $\Ls$. Then, analyzing the derived blowup $W := \bl_{\infty \times X}(\Pb^1 \times \Pb(\Ls \oplus \Oc))$, we note that the fiber of the natural map $W \to \Pb^1_S$ over
\begin{enumerate}
\item $0$, is equivalent to $\Pb(\Ls \oplus \Oc)$;
\item $\infty$ is the sum of the strict transform $D$ of $\infty \times \Pb(\Ls \oplus \Oc) \hook \Pb^1 \times \Pb(\Ls \oplus \Oc)$ and the exceptional divisor $\Pb(\Ls \oplus \Oc) \simeq \Ec \subset W$, the intersection of which is $X$.
\end{enumerate}
Moreover, the restrictions of $\Oc(D)$ and $\Oc(\Ec)$ to $X$ are equivalent to $\Ls^\vee$ and $\Ls$, respectively. Hence, if we denote by $q$ the natural map $W \to X$, we have 
\begin{align*}
\pi_*(1_{\Pb(\Ls \oplus \Oc)}) &= q_*\big(e(\Oc(1))(1_W)\big) \\
&= q_*\Big( F_\Hb \big( e(\Oc(D)), e(\Oc(\Ec)) \big) (1_W) \Big)
\end{align*} 
from which the desired equation follows immediately.
\end{proof}

The above Lemma allows us to prove that the derived double point formula holds for any $\Fc_S$-oriented Borel--Moore homology theory that has good Euler classes. Given a $\Fc_S$-oriented Borel--Moore homology theory $\Hb_\bullet$, we will use the notation
$$[V \to X] := f_*(1_V) \in \Hb_\bullet(X)$$
for every projective morphism $f: V \to X$ where $V$ is quasi-smooth.

\begin{prop}\label{prop:orhomfglthyddpt}
Suppose that $\Hb_\bullet$ has good Euler classes. Then, for each $X \in \Cc_S$, and  each projective morphism $\pi: W \to \Pb^1 \times X$ where $W$ is quasi-smooth, the equation
$$[W_0 \to X] = [D_1 \to X] + [D_2 \to X] - [\Pb_{D_1 \cap D_2}(\Oc(D_1) \oplus \Oc) \to X] \in \Hb_\bullet(X)$$
holds, where
\begin{enumerate}
\item $W_0$ is the fiber of $W \to \Pb^1 \times X$ over $0 \times X$;

\item the fiber $W_\infty$ of $W \to \Pb^1 \times X$ over $\infty \times X$ is the sum of virtual Cartier divisors $D_1$ and $D_2$ on $W$.
\end{enumerate}
\end{prop}
\begin{proof}
Suppose first that $W = X$. Then, denoting by $j$ the derived regular embedding $D_1 \cap D_2 \hook W$, we compute that 
\begin{align*}
&[W_0 \to W] \\
&= e(\Oc(1))(1_W) \\
&= e(\Oc(D_1 + D_2))(1_W) \\
&= [D_1 \to W] + [D_2 \to W] \\
&+ \sum_{i,j \geq 1} b_{ij} e(\Oc(D_1))^{i-1} \circ e(\Oc(D_2))^{j-1} ([D_1 \cap D_2 \to W]) \\
&= [D_1 \to W] + [D_2 \to W] \\
&+ j_*\bigg( \sum_{i,j \geq 1} b_{ij} e(\Oc(D_1))^{i-1} \circ e(\Oc(D_2))^{j-1}(1_{D_1 \cap D_2}) \bigg) \\
&= [D_1 \to W] + [D_2 \to W] \\
&+ j_*\bigg( \sum_{i,j \geq 1} b_{ij} e(\Oc(D_1))^{i-1} \circ e(\Oc(-D_1))^{j-1}(1_{D_1 \cap D_2}) \bigg) \\
&= [D_1 \to W] + [D_2 \to W] - j_*\big([\Pb_{D_1 \cap D_2}(\Oc(D_1) \oplus \Oc) \to D_1 \cap D_2]\big) \\
&= [D_1 \to W] + [D_2 \to W] - [\Pb_{D_1 \cap D_2}(\Oc(D_1) \oplus \Oc) \to W],
\end{align*}
where, in the fifth equation, we have used the fact that $\Oc(1) \simeq \Oc(D_1 + D_2)$ restricts to the trivial line bundle on $W_\infty$. The general case follows from this by pushing forward to $X$.
\end{proof}

We have proven the following universal property for $\PCob^S_*$. 

\begin{thm}\label{thm:orhomprecobunivprop}
Let $S$ be a finite-Krull-dimensional Noetherian derived scheme that admits an ample line bundle. Then,
\begin{enumerate}
\item $\PCob^S_*$ is the universal $\Fc_S$-oriented Borel-Moore homology theory that satisfies the weak projective bundle formula: given another such a theory $\Hb_\bullet$, there exists a unique transformation of $\Fc_S$-oriented Borel--Moore homology theories $\PCob^S_* \to \Hb_\bullet$;

\item $\PCob^S_*$ is the universal $\Fc_S$-oriented Borel-Moore homology theory that has good Euler classes: given another such a theory $\Hb_\bullet$, there exists a unique transformation of $\Fc_S$-oriented Borel--Moore homology theories $\PCob^S_* \to \Hb_\bullet$. \qed
\end{enumerate}
\end{thm}

\chapter{Comparison with Levine--Morel's algebraic bordism}\label{ch:LMcomp}

Here, we prove that the algebraic bordism homology groups, obtained from the bivariant algebraic cobordism considered in Section \ref{sect:bivAcob}, recover the algebraic bordism groups of Levine and Morel for quasi-projective schemes over a field of characteristic 0. We do this by first extending Lowrey--Schürg's algebraic bordism \cite{lowrey--schurg} to a bivariant theory $d\Omega^*_k$ on $\Fc_k$, and then establishing an equivalence of this bivariant theory with $\Omega^*_k$. Throughout the chapter, $k$ will denote a field.
 
\section{Lowrey--Schürg's algebraic bordism $d\Omega^k_*$ and its bivariant extension}

Here, we recall the construction of Lowrey--Schürg algebraic bordism groups over fields, with the slight adjustments made in \cite{annala-pre-and-cob}. Simultaneously, we extend it to a bivariant theory on $\Fc_k$. Recall that $\Ab^{\Fc_k}_*$ is a $\Fc_k$-oriented Borel--Moore homology theory\footnote{In fact, by Lemma \ref{lem:univorhom}, it is the universal one.}; given such a theory $\Hb_\bullet$, and a collection of subsets $\Sc(X) \subset \Hb_\bullet(X)$, we will denote by $\langle \Sc \rangle_\Hb$ the \emph{$\Fc_k$-oriented-homology ideal} generated by $\Sc$, i.e., the smallest collection of subsets that contains $\Sc$ and that is stable under  sums, pushforwards, Gysin pullbacks, and external products by arbitrary elements.

We will construct the bivariant theory $d \Omega^*_k$ and the homology theory $d \Omega^k_*$ as quotients of $\Lb^* \otimes \Ab_{\Fc_k}^*$ and $\Lb_* \otimes \Ab^{\Fc_k}_*$, respectively, where $\Lb^*$ and $\Lb_*$ are the the Lazard ring, equipped with its cohomological (non-positive) and homological (non-negative) grading, respectively. As $\Lb^* \otimes \Ab^*_{\Fc_k}$, and any quotient $\Bb$ thereof, is generated by orientations, Proposition \ref{prop:homidchar} implies that the $\Fc_k$-oriented-homology ideals of $\Bb_\bullet$ are exactly the restrictions of bivariant ideals of $\Bb$ to $\Bb_\bullet$.

The construction is given in several steps. 

\begin{cons}[\emph{Naive derived cobordism}, cf. \cite{lowrey--schurg} Definition 3.4]\label{cons:naivdercob}
We define
$$d\Omega^{k,\mathrm{naive}}_* := \Lb_* \otimes \Ab^{\Fc_k}_* / \Rc^\mathrm{fib},$$
where $\Rc^\mathrm{fib}$ is the homology ideal of homotopy fiber relations. In other words, $\Rc^\mathrm{fib}(X)$ is the $\Lb$-module generated by elements of form
$$[W_0 \to X] - [W_\infty \to X] \in \Lb_* \otimes \Ab^k_*(X),$$
where $W \to \Pb^1 \times X$ is a projective morphism with $W$ quasi-smooth over $k$, and where $W_0$ and $W_\infty$ are the fibers of $W \to \Pb^1$ over $0$ and $\infty$ respectively. 

Similarly, we define 
$$d\Omega_{k,\mathrm{naive}}^* := \Lb^* \otimes \Ab_{\Fc_k}^* / \Rc^\mathrm{fib},$$
where $\Rc^\mathrm{fib}(X \to Y)$ is the $\Lb$-module generated by elements of form
$$[W_0 \to X] - [W_\infty \to X] \in \Lb^* \otimes \Ab^*_k(X),$$
where $W \to \Pb^1 \times X$ is a projective morphism such that the composition $W \to \Pb^1 \times Y$ is quasi-smooth. Clearly, $d\Omega^{k,\mathrm{naive}}_*$ is the associated homology theory of $d\Omega_{k,\mathrm{naive}}^*$.
\end{cons}

Euler classes of vector bundles can be expressed in terms of vanishing loci in naive derived cobordism.

\begin{lem}\label{lem:eulerclasssectls}
Let $X \in \Cc_k$ and let $E$ be a vector bundle on $X$. Then, for every global section $s$ of $E$, 
$$e(E) = [V_X(s) \to X] \in d\Omega^*_{k,\naive}(X).$$
\end{lem}
\begin{proof}
Same as that of Lemma \ref{lem:eulerclasssect}.
\end{proof}

In particular, Euler classes of globally generated vector bundles are nilpotent. Let us denote by
$$F_\univ(x,y) \in \Lb[[x,y]]$$
the universal formal group law. 

\begin{cons}[\emph{Derived precobordism}, cf. \cite{lowrey--schurg} Definition 3.16 and Lemma 3.17]\label{cons:derprecob}
The theory $d\Omega^*_{k,\mathrm{pre}}$ is constructed from $\Omega^*_{k,\mathrm{naive}}$ in two steps. 
\begin{enumerate}
\item We define $\Rc^\mathrm{fgl}(X) \subset d \Omega^{k,\naive}_*(X)$ as consisting of elements of form
$$\big(e(\Ls_1 \otimes \Ls_2) - F_\univ(e(\Ls_1), e(\Ls_2))\big) \bullet 1_X,$$
where $X$ is smooth over $k$ and $\Ls_i$ are globally generated line bundles on $X$. Set
$$d \Omega'^{k, \pre}_* := d \Omega^{k,\naive}_* / \langle \Rc^\mathrm{fgl} \rangle_{d \Omega^{k,\naive}}.$$

Similarly, we define $\Rc^\mathrm{fgl}(X) \subset d \Omega_{k,\naive}^*(X)$ as consisting of elements of form 
$$e(\Ls_1 \otimes \Ls_2) - F_\univ(e(\Ls_1), e(\Ls_2)) $$
where $X$ is smooth over $k$ and $\Ls_i$ are globally generated line bundles on $X$, and set
$$d \Omega'^*_{k, \pre} := d \Omega^*_{k, \naive} / \langle \Rc^\mathrm{fgl} \rangle_{d \Omega_{k,\naive}}.$$
Since the orientations along smooth morphisms are isomorphisms (Proposition \ref{prop:veryspecialstrong}), $d\Omega'^{k,\pre}_*$ is the associated homology theory of $d\Omega'^*_{k,\pre}$.

\item Secondly, we set 
$$d\Omega^{k, \pre}_* := d\Omega'^{k, \pre}_* / \langle \Rc^\mathrm{fgl}_+ \rangle_{d \Omega'^{k, \pre}},$$ 
where the subset $\Rc^\mathrm{fgl}_{+}(X) \subset d\Omega'^{k, \pre}_*(X)$, for $X$ smooth over $k$, consists of elements of form
$$\Big( e(\Ls) - F_\univ\bigl(e(\Ls_1), i_\univ(e(\Ls_2)) \bigr)\Big) \bullet 1_X,$$
where $\Ls_1$ and $\Ls_2$ are globally generated line bundles on $X$, $\Ls \simeq \Ls_1 \otimes \Ls_2^\vee$, and $i_\univ$ is the formal inverse power series of $F_\univ$. 

Similarly, we set 
$$d\Omega_{k, \pre}^* := d\Omega'^*_{k, \pre} / \langle \Rc^\mathrm{fgl}_+ \rangle_{d \Omega'_{k, \pre}},$$ 
where $\Rc^\mathrm{fgl}_{+}(X) \subset d\Omega'^*_{k, \pre}(X)$, for $X$ smooth over $k$, consists of 
$$ e(\Ls) - F_\univ\bigl(e(\Ls_1), i_\univ(e(\Ls_2)) \bigr),$$
where $\Ls_1$ and $\Ls_2$ are globally generated line bundles on $X$, and where $\Ls \simeq \Ls_1 \otimes \Ls_2^\vee$. Clearly, $d\Omega^{k,\pre}_*$ is the associated homology theory of $d\Omega^*_{k,\pre}$.
\end{enumerate}
\end{cons}

\begin{lem}\label{lem:derprecobfgl}
The bivariant theory $d \Omega^*_{k, \pre}$ has good Euler classes (Definition \ref{def:fglthy}).
\end{lem}
\begin{proof}
Let us verify that $d \Omega^*_{k,\pre}$ satisfies the conditions of Definition \ref{def:fglthy}.

\begin{enumerate}
\item If $\Ls_1$ and $\Ls_2$ are globally generated, then their Euler classes are nilpotent, as is
$$e(\Ls_1 \otimes \Ls_2^\vee) = F_\univ\big(e(\Ls_1), i_\univ(e(\Ls_2))\big).$$
As every line bundle is equivalent to one of the form $\Ls_1 \otimes \Ls_2^\vee$ with $\Ls_i$ globally generated, if follows that Euler classes of line bundles are nilpotent.

\item This is Lemma \ref{lem:eulerclasssectls}.

\item This is proven similarly to Proposition \ref{prop:pbffgl}. \qedhere
\end{enumerate} 
\end{proof}

We are finally ready to construct $d \Omega^*_k$ and $d \Omega^k_*$.

\begin{cons}[\emph{Derived algebraic cobordism}, cf. \cite{lowrey--schurg} Definition 3.20]\label{def:deralgcob} 
We define $\Rc^\snc \subset d \Omega_*^{k,\pre}$ as the collection of homology elements containing
$$\zeta_{W,D} - 1_D \in d \Omega^{k,\pre}_*(D)$$
where $D \hookrightarrow W$ ranges over all $k$-snc divisors in all $W$ that are smooth over $k$, and $\zeta_{W,D}$ is as in Definition \ref{def:zetaclass}. Then the \emph{derived algebraic $k$-bordism} is defined as
$$d \Omega^k_*(X) := d \Omega^{k,\pre}_* / \langle \Rc^\snc \rangle_{d\Omega^{k,\pre}},$$
It is the homology theory associated to the bivariant theory
$$d\Omega^*_k := d \Omega^*_{k, \pre}/ \langle \Rc^\snc \rangle_{d \Omega_{k,\pre}}.$$
\end{cons}

For a quasi-projective classical $k$-scheme $X$, we denote by $\Omega_*(X)$ the algebraic bordism group of Levine--Morel \cite{levine-morel}. By the results of Levine--Pandharipande \cite{levine-pandharipande} it admits a presentation where the generators are isomorphisms classes of smooth $k$-varieties $V$ together with a projective map $V \to X$, and the relations are given by double point relations, the derived analogues of which we have used in he definition of universal precobordism. The following is Theorem 5.12 of \cite{lowrey--schurg}.

\begin{thm}\label{thm:lsmainresult}
If $k$ is a field of characteristic 0, then, for all quasi-projective derived $k$-schemes $X$, the map
$$\Omega_*(X_\cl) \to d \Omega^k_*(X),$$
given by the formula
$$[V \to X_\cl] \mapsto [V \to X],$$
is an isomorphism.
\end{thm}

\section{Comparison of $\Omega_k^*$ and $d\Omega_k^*$}

Here, we prove that the bivariant theories $\Omega^*_k$ and $d\Omega_k^*$ are equivalent. In particular, we obtain the following result.

\begin{cor}\label{cor:bivcobextendsLM}
Let $k$ be a field of characteristic 0. Then, for all $X \in \Cc_k$, the map
$$\Omega_*(X_\cl) \to \Omega^k_*(X),$$
given by the formula
$$[V \to X_\cl] \mapsto [V \to X],$$
is an isomorphism.
\end{cor}

In other words, in characteristic 0, the bivariant algebraic $k$-cobordism is a bivariant extension of Levine--Morel's algebraic bordism. We begin with the following observations.

\begin{lem}\label{lem:comparisonmap1}
There exists a unique orientation preserving Grothendieck transformation 
$$\eta: \PCob^*_k \to d\Omega^*_{k, \pre}.$$
\end{lem}
\begin{proof}
By Lemma \ref{lem:derprecobfgl}, $d\Omega^*_{k,\pre}$ has good Euler classes, so the claim follows from one of the universal properties of $\PCob^*_k$ given in Theorem \ref{thm:precobunivprop}.
\end{proof}

\begin{lem}\label{lem:comparisonmap2}
There exists a unique $\Lb$-linear orientation preserving Grothendieck transformation 
$$\eta': d\Omega^*_{k, \pre} \to \PCob^*_k,$$
where $\PCob^*_k$ has the $\Lb$-linear structure provided by its formal group law (Theorem \ref{thm:fgl}).
\end{lem}
\begin{proof}
Indeed, since $d\Omega^*_{k, \pre}$ is constructed as a quotient of $\Lb^* \otimes \Ab^*_{\Fc_k}$, there exists at most one such a Grothendieck transformation. The proof then reduces to checking that the unique transformation $\Lb^* \otimes \Ab^*_{\Fc_k} \to \PCob^*_k$ kills the relations imposed in Constructions \ref{cons:naivdercob} and \ref{cons:derprecob}, which is obvious.
\end{proof}

As the composition
$$\PCob^*_k \xto{\eta} d\Omega^*_{k, \pre} \xto{\eta'} \PCob^*_k$$
is the identity transformation, we have almost proven that $\eta$ and $\eta'$ are inverse isomorphisms. This is finished by the following observation.

\begin{lem}\label{lem:etaLlin}
The transformation $\eta: \PCob^*_k \to d\Omega^*_{k, \pre}$ is $\Lb$-linear.
\end{lem}
\begin{proof}
We wish to show that $\eta: \PCob^*_k(\Spec(k)) \to d\Omega^*_{k, \pre}(\Spec(k))$ is a map of $\Lb$-algebras. By the universal property of the Lazard ring, it suffices to show that $\eta$ sends the formal group law
$$F(x,y) = x + y + \sum_{i,j \geq 1} a_{ij} x^i y^j \in \PCob^*_k(\Spec(k))[[x,y]]$$
to the formal group law
$$F'(x,y) = x + y + \sum_{i,j \geq 1} a'_{ij} x^i y^j \in d\Omega^*_k(\Spec(k))[[x,y]].$$
As $\eta$ sends Euler classes to Euler classes, for all $n,m \geq 0$
\begin{align*}
&e(\Oc(1,0)) + e(\Oc(0,1)) + \sum_{i=1}^n \sum_{j=1}^m \eta(a_{ij})  e(\Oc(1,0))^i \bullet e(\Oc(0,1))^j \\
=& \eta(e(\Oc(1,1))) \\
=& e(\Oc(1,1)) \\
=& e(\Oc(1,0)) + e(\Oc(0,1)) + \sum_{i=1}^n \sum_{j=1}^m a'_{ij}  e(\Oc(1,0))^i \bullet e(\Oc(0,1))^j;
\end{align*}
however, since we do not know if $d\Omega^*_{k, \pre}$ satisfies the projective bundle formula, the claim does not follow immediately. Nonetheless, if $\eta(a_{ij})$ and $a'_{ij}$ differ for some $1 \leq i \leq n$ and $1 \leq j \leq m$, then, in the partial ordering where $(i,j) \leq (i',j')$ if $i \leq i'$ and $j \leq j'$, we can find a minimal pair $(i',j')$ where they differ, and then
\begin{align*}
&0\\
&= \pi_! \bigg( e(\Oc(1,0))^{n-i'} \bullet e(\Oc(0,1))^{m-j'} \\
& \bullet \sum_{i=1}^n \sum_{j=1}^m \big(\eta(a_{ij}) - a'_{ij}\big)  e(\Oc(1,0))^i \bullet e(\Oc(0,1))^j \bigg) \\
&= \eta(a_{i'j'}) - a'_{i'j'},
\end{align*}
where $\pi$ is the structure morphism $\Pb^n_k \times_{\Spec(k)} \Pb^n_k$, which is a contradiction since by assumption $\eta(a_{i'j'}) \not = a'_{i'j'}$.
\end{proof}

We have shown that the Grothendieck transformations $\eta$ and $\eta'$ are inverse isomorphisms. The following result is an immediate consequence of this.

\begin{prop}\label{prop:precobandcob}
There exists a unique orientation preserving Grothendieck transformation
$$\Omega_k^* \to d\Omega_k^*,$$
which is an isomorphism of bivariant theories. 
\end{prop}
\begin{proof}
Indeed, $\Omega_k^*$ and $d\Omega_k^*$ are obtained from $\PCob^*_k$ and $d\Omega^*_{k, \pre}$ by imposing equivalent relations.
\end{proof}

Corollary \ref{cor:bivcobextendsLM} follows immediately from the above and Theorem \ref{thm:lsmainresult}.

\chapter{Twisting bivariant theories}\label{ch:bivtwist}

In this Chapter, we study the process of altering the orientation of a bivariant theory called \emph{twisting}, which is inspired by the closely related twisting operation for oriented cohomology theories in topology \cite{quillen:1971} and oriented Borel--Moore homology theories in algebraic geometry \cite{levine-morel}. Most of these results appeared originally in \cite{annala-cob}.

Throughout the section, $S$ will be a finite-Krull-dimensional Noetherian derived scheme that admits an ample line bundle. We will work with the restricted functoriality $\Fc_S$, as we shall make use of Chern classes, which we have defined only for those Noetherian derived schemes that admit an ample line bundle.

\section{Twisting oriented bivariant theories}\label{sect:bivtwist}

Throughout this section, $\Bb$ will denote a stably oriented bivariant theory with functoriality $\Fc_S$ that has good Euler classes. By Theorem \ref{thm:precobunivprop}, there exists a unique orientation preserving Grothendieck transformation $\PCob^*_A \to \Bb$. In particular, $\Bb$ has a theory of Chern classes which satisfies the properties listed in Theorem \ref{thm:chernclass}.

\begin{defn}\label{def:toddclass}
Let $\tau := (\tau_0, \tau_1, \tau_2, ...)$ be an infinite sequence of elements of $\Bb(\Spec(A))$ with $\tau_0 = 1$. Then, for any line bundle $\Ls$ on $X \in \Cc_S$, we define its \emph{inverse ($\tau$-)Todd-class} by
$$\Td^{-1}_\tau(\Ls) := \sum_i \tau_i c_1(\Ls)^i \in \Bb^\bullet(X).$$
By splitting principle, there exists a unique extension of this to all vector bundles, satisfying, that, whenever
$$E' \to E \to E''$$
is a cofiber sequence of vector bundles on $X \in \Cc_S$, then
$$\Td^{-1}_\tau(E) = \Td^{-1}_\tau(E') \bullet \Td^{-1}_\tau(E'') \in \Bb(X).$$
Using this, $\Td^{-1}_\tau(\alpha)$ may be defined for all $K$-theory classes $\alpha \in K^0(X)$. The \emph{($\tau$-)Todd-class} $\Td_\tau(\alpha)$ is defined to be the inverse of $\Td^{-1}_\tau(\alpha)$.
\end{defn}

We can now define the twisting operation.

\begin{defn}\label{def:bivtwist}
Let $\tau$ be as above. Then we define the \emph{twisted} oriented bivariant theory $\Bb^{(\tau)}$ to be the bivariant theory $\Bb$ equipped with the orientation $\theta^{(\tau)}$ which, for a quasi-smooth morphism $f: X \to Y$, is defined by $\theta^{(\tau)}(f) := \Td_\tau(\Lb^\vee_{X/Y}) \bullet \theta(f).$ 
\end{defn}

\begin{lem}\label{lem:bivtwistor}
The twisted theory $\Bb^{(\tau)}$ is a stably oriented bivariant theory.
\end{lem}
\begin{proof}
Clearly $\theta^{(\tau)}(\Id) = 1$. Moreover, given quasi-smooth morphisms $f: X \to Y$ and $g: Y \to Z$, we compute that
\begin{align*}
\theta^{(\tau)}(f) \bullet \theta^{(\tau)}(g) &= \Td_\tau(\Lb^\vee_{X/Y}) \bullet \theta(f) \bullet \Td_\tau(\Lb^\vee_{Y/Z}) \bullet \theta(g) \\
&= \Td_\tau(\Lb^\vee_{X/Y}) \bullet \Td_\tau(f^*\Lb^\vee_{Y/Z}) \bullet \theta(f) \bullet \theta(g) & (\text{commutativity}) \\
&= \Td_\tau(\Lb^\vee_{X/Z}) \bullet \theta(g \circ f) \\
&= \theta^{(\tau)}(g \circ f);
\end{align*}
hence, $\theta^{(\tau)}$ is an orientation. It is stable under pullbacks as Todd classes of duals of relative cotangent complexes are, and as $\theta$ is a stable orientation.
\end{proof}

For clarity, we will distinguish Gysin pushforwards, Chern classes and Euler classes in the twisted theory $\Bb^{(\tau)}$ by marking them with the symbol $(\tau)$ in a convenient location. Moreover, we will use the notation
$$\lambda_\tau(x) := \sum_{i=0}^\infty \tau_i x^{i+1} \in \Bb^\bullet(S)[[x]].$$
We then have the following result.

\begin{lem}\label{lem:bivtwisteuler}
For all $X \in \Cc_S$ and all line bundles $\Ls$ on $X$, we have
$$e^{(\tau)}(\Ls) = \lambda_{\tau}(e(\Ls)) \in \Bb^\bullet(X).$$
\end{lem}
\begin{proof}
Indeed, denoting by $i$ the derived inclusion $V_X(0) \hook X$, where $0$ is the zero-section of $\Ls$, we have
\begin{align*}
e^{(\tau)}(\Ls) &= i^{(\tau)}_!\big(1_{V_X(0)}\big) \\
&= i_!\big( \Td_{\tau}(\Lb^\vee_{V_X(0) / X)}) \big) \\
&= i_!\big( \Td_{\tau}(\Ls\vert_{V_X(0)}[-1]) \big) \\
&= \Td^{-1}_{\tau}(\Ls) \bullet i_!(1_{V_X(0)}) \\
&= \Td^{-1}_{\tau}(\Ls) \bullet e(\Ls),
\end{align*}
from which the claim immediately follows.
\end{proof}

Let us denote by $\lambda_\tau^{-1}(x)$ the inverse of $\lambda_\tau(x)$ with respect to composition. Then, we define $\bar \tau := (\bar \tau_0, \bar \tau_1, \bar \tau_2, ...)$ to be the sequence satisfying
$$\lambda^{-1}_\tau(x) = \sum_{i=0}^\infty \bar\tau_i x^{i+1} \in \Bb^\bullet(S)[[x]].$$
The following result is useful.

\begin{lem}\label{lem:inversetwist}
Let $\Bb$ be a stably oriented bivariant theory that has good Euler classes. Then, for all $\alpha \in K^0(X)$, we have 
$$\Td^{-1}_{\bar \tau} (\alpha; \Bb^{(\tau)}) \bullet \Td^{-1}_{\tau} (\alpha) = 1_X \in \Bb^\bullet(X),$$
where $\Td^{-1}_{\bar \tau} (\alpha; \Bb^{(\tau)})$ denotes the inverse $\bar\tau$-Todd class of $E$, computed in the oriented bivariant theory $\Bb^{(\tau)}$. In particular, $\big(\Bb^{(\bar\tau)} \big)^{(\tau)} = \Bb = \big(\Bb^{(\tau)} \big)^{(\bar\tau)}$.
\end{lem}
\begin{proof}
By splitting principle and multiplicativity of Todd classes, it suffices to check that this formula holds for $\alpha = [\Ls]$, where $\Ls$ is a line bundle on $X$. This boils down to showing that
$$\frac{\lambda^{-1}_\tau}{x}\big( \lambda_\tau(x) \big) \frac{\lambda_\tau}{x}\big(x\big) = 1 \in \Bb^\bullet(S)((x)),$$
where for a power series $\psi(x)$ with zero constant coefficient, we have denoted by $(\psi/x)(x)$ the power series such that $x \cdot (\psi/x)(x) = \psi(x)$. This formula follows from the simple computation
\begin{align*}
\frac{\lambda^{-1}_\tau}{x}\big( \lambda_\tau(x) \big)  \frac{\lambda_\tau}{x}\big(x\big) &= \frac{\lambda^{-1}_\tau(\lambda_\tau(x))}{\lambda_\tau(x)} \frac{\lambda_\tau(x)}{x} \\
&= \frac{x}{\lambda_\tau(x)} \frac{\lambda_\tau(x)}{x} \\
&= 1,
\end{align*}
so we are done.
\end{proof}

By the next result, twisting preserves pleasant properties of Euler classes.

\begin{prop}\label{prop:bivtwistfglthy}
The bivariant theory $\Bb^{(\tau)}$ has good Euler classes. Moreover, the formal group law of $\Bb^{(\tau)}$ is given by
$$F_{\Bb^{(\tau)}}(x,y) = \lambda_{\tau}\big( F_\Bb\big(\lambda^{-1}_\tau(x), \lambda^{-1}_\tau(y) \big) \big).$$
\end{prop}
\begin{proof}
We verify that $\Bb^{(\tau)}$ satisfies the conditions of Definition \ref{def:fglthy}.
\begin{enumerate}
\item Nilpotence of the Euler classes of line bundles of $\Bb^{(\tau)}$ follows immediately from the corresponding fact for $\Bb$.
\item Suppose that $X \in \Cc_S$, $\Ls$ is a line bundle on $X$, and $s$ is a global section of $\Ls$. Then, denoting by $i_s$ the derived regular embedding $V_X(s) \hook X$, we have that  
\begin{align*}
e^{(\tau)}(\Ls) &= \sum_{i=0}^\infty b_i e(\Ls)^{i+1} \\
&= \sum_{i=0}^\infty b_i e(\Ls)^{i} \bullet i_{s!}(1_{V_X(s)}) \\
&= i_{s!} \bigg( \sum_{i=0}^\infty b_i e(\Ls)^{i} \bigg) \\
&= i_{s!} \big( \Td(\Ls^\vee_{V_X(s)/X}) \bullet \theta(i_s) \big) \\
&= i^{(\tau)}_{s!} (1_{V_X(s)}),
\end{align*}
as desired.

\item We compute that
\begin{align*}
e^{(\tau)}(\Ls_1 \otimes \Ls_2) &= \lambda_\tau \Big( F_\Bb\big(e(\Ls_1), e(\Ls_2)\big)\Big) \\
&=  \lambda_\tau \Big( F_\Bb\big(\lambda^{-1}_\tau\big(e^{(\tau)}(\Ls_1)\big), \lambda^{-1}_\tau\big(e^{(\tau)}(\Ls_2)\big)\big)\Big),
\end{align*}
as desired. \qedhere
\end{enumerate}
\end{proof}

\section{Grothendieck--Riemann--Roch theorem}\label{sect:grr}

Here, we prove our version of the Grothendieck--Riemann--Roch theorem. Essentially, it is a cohomological analogue of the Grothendieck--Riemann--Roch for singular varieties, proven by Baum--Fulton--MacPherson \cite{baum:1975}. We begin by defining the two bivariant theories that will play an important role in the theorem.

\begin{defn}\label{def:univaddfglthy}
The \emph{universal additive $S$-precobordism} is defined as
$$\PCob^*_{S,a} := \Zb_a \otimes_\Lb \PCob^*_S,$$
where $\Zb_a$ is the integers, considered as an $\Lb$-algebra via the additive formal group law $x + y$. Similarly, the \emph{universal additive $S$-cobordism}
is defined as
$$\Omega^*_{S,a} := \Zb_a \otimes_\Lb \Omega^*_S.$$
\end{defn}

By the results of Levine--Morel and Corollary \ref{cor:bivcobextendsLM}, if $A = k$ is a field of characteristic 0, then the associated homology theory of $\Omega^*_{k,a}$ recovers the Chow groups (see e.g. \cite{fulton:1998}). Hence, one should regard $\Omega^*_{S,a}$ as a ``bivariant intersection theory''; in particular the associated cohomology theory $\Omega^*_{S,a}$ provides a candidate for the elusive Chow cohomology theory. Its relationship with other (partial) candidates \cite{quillen:1973, fulton:1975, levine-weibel} remains poorly understood.

\begin{defn}\label{def:univmultfglthy}
The \emph{universal multiplicative $S$-precobordism} is defined as 
$$\PCob_{S,m} := \Zb_m \otimes_\Lb \PCob_S,$$
where $\Zb_m$ is the integers, considered as an $\Lb$-algebra via the multiplicative formal group law $x + y - xy$. Similarly, the \emph{universal multiplicative $S$-cobordism}
is defined as
$$\Omega_{S,m} := \Zb_m \otimes_\Lb \Omega_S.$$
\end{defn}

By Theorem \ref{thm:cf}, the associated cohomology theory of $\PCob_{S,m}$ is naturally equivalent to the $K$-theory of perfect complexes $K^0$. As the snc-relations (see Definition \ref{def:bivAcob}) are satisfied in bivariant $K$-theory, also $\Omega_{S,m}$ has $K^0$ as its associated cohomology theory.

Next, we will prove that the bivariant theories $\Qb \otimes \PCob^*_{S,a}$ and $\Qb \otimes \PCob_{S,m}$, and $\Qb \otimes \Omega^*_{S,a}$ and $\Qb \otimes \Omega_{S,m}$ are equivalent as unoriented bivariant theories. Moreover, we express the difference between the orientations using a relatively simple inverse Todd-class.

Let $\tau := (\tau_0, \tau_1, \tau_2, ...)$ be the sequence of rational numbers that satisfies 
$$\lambda_\tau(x) = 1 - e^{-x} \in \Qb[[x]].$$
One solves from the above that
$$\lambda^{-1}_\tau(x) = \lambda_{\bar \tau}(x) = - \ln(1 - x) \in \Qb[[x]].$$
Twisting a $\Qb$-linear oriented bivariant theory by these sequences has the following interesting property.

\begin{lem}\label{lem:addtomult}
Twisting by $\tau$ changes the additive formal group law $x + y$ to the multiplicative $x + y - xy$. Twisting by $\bar \tau$ has the opposite effect.
\end{lem}
\begin{proof}
It suffices to prove the first claim. Using Lemma \ref{prop:bivtwistfglthy}, we compute that if a $\Qb$-linear oriented bivariant theory $\Bb$ has formal group law $x + y$, then $\Bb^{(\tau)}$ has the formal group law
\begin{align*}
1 - \exp(\ln(1 - x) + \ln(1 - y)) &= 1 - (1-x)(1-y)\\
&= x + y -xy,
\end{align*}
as desired.
\end{proof}

Clearly $\PCob^*_{S,a}$ and $\PCob_{S,m}$ are the universal stably oriented bivariant theories that have good Euler classes with formal group laws $x+y$ and $x + y - xy$, respectively. On cohomology rings, the transformation induced by the universal property acquires a familiar form, after identifying the source with $K$-theory.

\begin{lem}\label{lem:chernchar}
Let $\ch$ be the unique orientation preserving Grothendieck transformation $\PCob_{S,m} \to \Qb \otimes \PCob^{(\tau)}_{S,a}$. Then, the induced homomorphism
$$K^0(X) \to \Qb \otimes \PCob^*_{S,a}(C)$$
of cohomology rings sends $[\Ls]$ to $\exp(e(\Ls))$, where $\Ls$ is a line bundle on $X$.
\end{lem}
\begin{proof}
By definition, 
\begin{align*}
[\Ls] &\mapsto 1 - c_1^{(\tau)}(\Ls^\vee) \\
&= 1 - \big( 1 - e^{-c_1(\Ls^\vee)} \big) \\
&= \exp(c_1(\Ls)),
\end{align*}
as desired.
\end{proof}

\begin{lem}\label{lem:cherncharisom}
The unique orientation preserving Grothendieck transformation $\ch: \PCob_{S,m} \to \Qb \otimes \PCob^{(\tau)}_{S,a}$ induces an isomorphism
$$\ch_\Qb: \Qb \otimes \PCob_{S,m} \to \Qb \otimes \PCob^{(\tau)}_{S,a}.$$
The result remains true if we replace $\PCob_{S,m}$ and $\PCob^*_{S,a}$ by $\Omega_{S,m}$ and $\Omega^*_{S,a}$, respectively.
\end{lem}
\begin{proof}
Denote by $\ch'$ the unique orientation preserving Grothendieck transformation $\PCob^*_{S,a} \to \Qb \otimes \PCob^{(\bar \tau)}_{S,m}$. Clearly the compositions
$$\Qb \otimes \PCob_{S,m} \xto{\ch_\Qb} \Qb \otimes \PCob^{(\tau)}_{S,a} \xto{\ch'_\Qb} \Qb \otimes \PCob_{S,m}$$
and
$$\Qb \otimes \PCob^*_{S,a} \xto{\ch'_\Qb} \Qb \otimes \PCob^{(\bar\tau)}_{S,m} \xto{\ch_\Qb} \Qb \otimes \PCob^*_{S,a}$$
are orientation preserving, and therefore identity transformations. Hence, $\ch_\Qb$ and $\ch'_\Qb$ are inverse isomorphisms, proving the first claim. 

In order to prove the second claim, we show that $\ch'_\Qb$ is compatible with the snc relations in the sense that, if $W$ is smooth over $S$ and $D \simeq n_1 D_1 + \cdots + n_r D_r$ is an $S$-snc divisor on $W$, then
$$\ch'_\Qb(\eta_{W,D} - 1_D) = u \bullet \big(\eta_{W,D} - 1_D \big) \in \PCob_*^{S,a}(D),$$
where $u$ is a unit in $\PCob_{S,m}(D)$. As
$$\ch'_\Qb(1_D) = \Td_{\bar \tau}(\Lb^\vee_{D / S}) \bullet 1_D,$$
it suffices to show that 
$$\ch'_\Qb(\eta_{W,D}) = \Td_{\bar\tau}(\Lb^\vee_{D / S}) \bullet \eta_{W,D}.$$
which is the content of Lemma \ref{lem:zetatwist}.
\end{proof}

Suppose that $\Bb$ is a stably oriented bivariant theory on $\Fc_S$ that has good Euler classes, and let $\tau = (1, \tau_1, \tau_2,...)$ be a sequence of elements of $\Bb^\bullet(\Spec(A))$. Then, if $W$ is quasi-smooth and quasi-projective derived $S$-scheme, and if $D \simeq n_1 D_1 + \cdot + n_r D_r$ is a virtual Cartier divisor on $W$, we denote by 
$$\zeta_{W, D, D_1,...,D_r}, \zeta^{(\tau)}_{W, D, D_1,...,D_r} \in \Bb_\bullet(D)$$
the images of $\zeta_{W, D, D_1,...,D_r} \in \PCob^S_*(D)$ under the unique orientation preserving Grothendieck transformations $\PCob^*_S \to \Bb$ and $\PCob^*_S \to \Bb^{(\tau)}$, respectively.

\begin{lem}\label{lem:zetatwist}
Let everything be as above. Then,
$$\zeta^{(\tau)}_{W, D, D_1,...,D_r} = \Td_{\bar\tau}(\Lb^\vee_{D/S}) \bullet \zeta_{W, D, D_1,...,D_r}$$
in $\Bb_\bullet(D)$.
\end{lem}
\begin{proof}
Suppose that, for each $I \subset \{1,2,...,r\}$, we have a formal power series $f_I(x_1,...,x_r) \in \Bb^\bullet(S)[[x_1,...,x_r]]$, such that
\begin{enumerate}
\item $f_\emptyset = 0$;
\item $\sum_{I \subset [r]} \xbf^I f_I(x_1,...,x_r) =  [n_1]_{F_\Bb} \cdot x_1 +_{F_\Bb} \cdots +_{F_\Bb} [n_r]_{F_\Bb} \cdot x_R$,
\end{enumerate}
where $\xbf^I = \prod_{i \in I} x_i$. Then, applying the projection formula, we observe that
$$\zeta_{W, D, D_1,...,D_r} = \sum_{I \subset [r]} \iota^I_*\Big(f_I \big( e(\Oc(D_1)), ...,e(\Oc(D_r)) \big) \bullet 1_{D_I} \Big) \in \Bb_\bullet(D),$$
where $\iota^I$ is the closed immersion $\bigcap_{i \in I} D_i \hook D$.

Let us denote the formal group law of $\Bb$ and $\Bb^{(\tau)}$ by $F$ and $F_{(\tau)}$, respectively, and let $F^{n_1,...,n_r}$ and $F^{n_1,...,n_r}_{(\tau)}$ be the formal power series corresponding to the formal linear combination $[n]_{F} \cdot x_1 +_{F} \cdots +_{F} [n_r]_F \cdot x_r$ in $\Bb$ and $\Bb^{(\tau)}$, respectively. Moreover, let
$$F^{n_1,...,n_r}(x_1,...,x_r) = \sum_{I \subset [r]} \xbf^I F^{n_1,...,n_r}_I(x_1,...,x_r)$$
be the unique expression, where $F^{n_1,...,n_r}_I(x_1,...,x_r)$ contains only those variables, whose index lies in $I$. Then,
\begin{align*}
&´F^{n_1,...,n_r}_{(\tau)}(x_1,...,x_r) \\
&= \lambda_\tau\Big(F^{n_1,...,n_r}\big(\lambda_\tau^{-1}(x_1),...,\lambda_\tau^{-1}(x_1)\big)\Big) \\
&= \lambda_\tau\Big(F^{n_1,...,n_r}\big(\lambda_\tau^{-1}(x_1),...,\lambda_\tau^{-1}(x_1)\big)\Big) \\
&= \frac{\lambda_\tau}{x}\Big(F^{n_1,...,n_r}\big(\lambda_\tau^{-1}(x_1),...,\lambda_\tau^{-1}(x_1)\big)\Big) F^{n_1,...,n_r}\big(\lambda_\tau^{-1}(x_1),...,\lambda_\tau^{-1}(x_1)\big)\\
&= \sum_{I \subset [r]} \lambda_\tau^{-1}(\xbf)^I \frac{\lambda_\tau}{x}\Big(F^{n_1,...,n_r}\big(\lambda_\tau^{-1}(x_1),...,\lambda_\tau^{-1}(x_1)\big)\Big)\\
&\cdot F_I^{n_1,...,n_r}\big(\lambda_\tau^{-1}(x_1),...,\lambda_\tau^{-1}(x_1)\big) \\
&= \sum_{I \subset [r]}\xbf^I \frac{\lambda_\tau^{-1}}{x} (\xbf)^I \frac{\lambda_\tau}{x}\Big(F^{n_1,...,n_r}\big(\lambda_\tau^{-1}(x_1),...,\lambda_\tau^{-1}(x_1)\big)\Big) \\
&\cdot F_I^{n_1,...,n_r}\big(\lambda_\tau^{-1}(x_1),...,\lambda_\tau^{-1}(x_1)\big).
\end{align*}
Let us denote by $\Ls_i$ the restrictions of $\Oc_W(D_i)$ to $D$, by $\Ls$ the restriction of $\Oc_W(D)$ to $D$, and, for $I \subset [r]$, by $\Ls_I$ the sum $\bigoplus_{i \in I} \Ls_i$. By applying the observation we made at the beginning of the proof, we obtain the desired formula
\begin{align*}
& \zeta^{(\tau)}_{W, D, D_1,...,D_r} \\
&= \sum_{I \subset [r]} \iota^I_*\Big( \Td^{-1}_{\bar\tau}(\Ls_I; \Bb^{(\tau)}) \bullet \Td^{-1}_\tau(\Ls) \\
&\bullet F^{n_1,...,n_r}_I\big(e(\Ls_1),...,e(\Ls_r)\big) \bullet 1^{(\tau)}_{D_I} \Big) \\
&= \sum_{I \subset [r]} \iota^I_*\Big( \Td^{-1}_{\bar\tau}(\Ls_I; \Bb^{(\tau)}) \bullet \Td^{-1}_\tau(\Ls) \\
&\bullet F^{n_1,...,n_r}_I\big(e(\Ls_1),...,e(\Ls_r) \big) \bullet \Td_\tau(\Lb_{D_I/S}^\vee) \bullet 1_{D_I} \Big) \\
&= \sum_{I \subset [r]} \iota^I_*\Big( \Td^{-1}_{\bar\tau}(\Ls_I; \Bb^{(\tau)}) \bullet \Td^{-1}_\tau(\Ls) \bullet \Td_\tau(\Lb_{D_I/S}^\vee) \\
& \bullet F^{n_1,...,n_r}_I\big(e(\Ls_1),...,e(\Ls_r) \big) \bullet 1_{D_I} \Big) \\
&= \sum_{I \subset [r]} \iota^I_*\Big( \Td^{-1}_{\bar\tau}(\Ls_I; \Bb^{(\tau)}) \bullet \Td^{-1}_\tau(\Ls_I) \bullet \Td_\tau(\Lb_{D/S}^\vee) \\
& \bullet F^{n_1,...,n_r}_I\big(e(\Ls_1),...,e(\Ls_r) \big) \bullet 1_{D_I} \Big) \\
&= \Td_\tau(\Lb_{D/S}^\vee) \bullet \sum_{I \subset [r]} \iota^I_*\Big( \Td^{-1}_{\bar\tau}(\Ls_I; \Bb^{(\tau)}) \bullet \Td^{-1}_\tau(\Ls_I) \\
& \bullet  F^{n_1,...,n_r}_I\big(e(\Ls_1),...,e(\Ls_r) \big) \bullet 1_{D_I} \Big) \\
&= \Td_\tau(\Lb_{D/S}^\vee) \bullet \sum_{I \subset [r]} \iota^I_*\Big( F^{n_1,...,n_r}_I\big(e(\Ls_1),...,e(\Ls_r) \big) \bullet 1_{D_I} \Big) \\
&= \Td_\tau(\Lb_{D/S}^\vee) \bullet \zeta_{W,D,D_1,...,D_r},
\end{align*}
where, in the second-to-last step, we have used Lemma \ref{lem:inversetwist}.
\end{proof}

Combining the above results we obtain the following theorem.

\begin{thm}\label{thm:rr}
The maps
$$\ch: \Qb \otimes K^0(X) \to \Qb \otimes \PCob^*_{S,a}(X)$$
defined by the formula
$$[\Ls] \mapsto \exp(c_1(\Ls))$$
on line bundles, and extended additively to all elements of $\Qb \otimes K^0(X)$ using the splitting principle, are isomorphisms of rings that commute with pullbacks along morphisms in $\Cc_S$. Moreover, if $f: X \to Y$ is projective and quasi-smooth, and if $\alpha \in \Qb \otimes K^0(X)$, then
$$\ch(f_!(\alpha)) = f_!\big(\ch(\alpha) \bullet \Td_\tau(\Lb^\vee_{X/Y})\big),$$
where $\Td_\tau([\Ls]) = \frac{1 - e^{-c_1(\Ls)}}{c_1(\Ls)}$ and $\Td_\tau$ is extended to all elements in $K^0$ in a multiplicative fashion using the splitting principle. \qed
\end{thm}

\section{Algebraic cobordism with rational coefficients}\label{sect:Qcob}

Here, we show that the bivariant $S$-(pre)cobordism can be obtained from the universal additive theory by tensoring it with the rational Lazard ring $\Lb_\Qb := \Qb \otimes \Lb$ and twisting.

We begin by recalling some classical results on formal group laws (see e.g. \cite{hazewinkel:1978}). Consider the polynomial ring $B = \Qb[b_1, b_2, ...]$ with infinitely many generators. Denoting 
$$\lambda_\bbf(x) := x + b_1 x^2 + b_2 x^3 + \cdots \in B[[x]],$$ 
we define a formal group law on $B$ by the formula
$$F'(x,y) := \lambda_\bbf\big(\lambda^{-1}_\bbf(x) + \lambda^{-1}_\bbf(y)\big).$$
In classical terminology, $\lambda_\bbf(x)$ is the \emph{(inverse) logarithm} of the formal group law $F'$. The formal group law $F'$ is classified by a ring homomorphism $\Lb_\Qb \to B$. Let us denote by $\lambda_\Lb(x) \in \Lb_\Qb[[x]]$ the inverse logarithm of the universal formal group law. The coefficients of $\lambda_\Lb$ induce a morphism $B \to \Lb_\Qb$ that sends $\lambda_\bbf$ to $\lambda_\Lb$. As the compositions $\Lb_Q \to B \to \Lb_Q$ and $B \to \Lb_\Qb \to B$ are identities, we have identified $\Lb_\Qb$ with $B$.

Denote by $\bbf$ the sequence $(1,b_1,b_2,...)$ of elements of $B$. 

\begin{lem}\label{lem:Qcob}
The unique orientation preserving Grothendieck transformation
$$\Qb \otimes \PCob^*_S \to B \otimes \PCob^{(\bbf)}_{S,a}$$
is an isomorphism. The result remains true if we replace $\PCob^*_S$ with $\Omega^*_S$ and $\PCob^*_{S,a}$ with $\Omega^*_{S,a}$.
\end{lem}
\begin{proof}
By construction and Proposition \ref{prop:bivtwistfglthy}, the formal group law of $\PCob^{(\bbf)}_{S,a}$ is given by $F'$. Moreover, as $\Qb \otimes \PCob^{(\bar\bbf)}_S$ has the additive formal group law, there exists a unique orientation preserving Grothendieck transformation
$$B \otimes \PCob^*_{S,a} \to \Qb \otimes \PCob^{(\bar\bbf)}_S$$
that sends $\lambda_\bbf$ to the image of $\lambda_\Lb$ in $Q \otimes \PCob^*_S(S)$. As the composition
$$\Qb \otimes \PCob^*_S \to B \otimes \PCob^{(\bbf)}_{S,a} \to \Qb \otimes \PCob^*_S$$
is orientation preserving, and as
$$B \otimes \PCob^*_{S,a} \to \Qb \otimes \PCob^{(\bar\bbf)}_S \to B \otimes \PCob^*_{S,a}$$
is orientation preserving and preserves $B$, they are the identity transformations. This proves the first claim. By Lemma \ref{lem:zetatwist} the identification we have obtained is compatible with the snc relations, and therefore the second claim follows. 
\end{proof}

The following result follows from the above by combining it with Theorem \ref{thm:cf} and Lemma \ref{lem:cherncharisom}

\begin{thm}\label{thm:QcobvsKthy}
For each quasi-projective derived $S$-scheme, there exists an isomorphism of rings
$$\Qb \otimes \PCob^*_S(X) \to \Lb_\Qb \otimes K^0(X).$$
Moreover, these maps commute with pullbacks along morphisms in $\Cc_S$ and Gysin pushforwards along quasi-smooth and projective morphisms $f: X \to Y$ such that $[\Lb_{X/Y}] = 0 \in K^0(X)$. \qed
\end{thm}

In particular, we can compute the precobordism with rational coefficients of local rings.

\begin{cor}\label{cor:Qcoboflocr}
Let $A$ be a local Noetherian derived ring of finite Krull dimension. Then $\Qb \otimes \PCob^*(\Spec(A)) = \Lb^*_\Qb$. \qed
\end{cor}


\chapter{Conclusion}\label{ch:concl}

Algebraic cobordism was first introduced as an extraordinary motivic cohomology theory in \cite{voevodsky:1996} and later in a geometric form in \cite{levine-morel}. In this thesis, we have constructed a geometric theory of bivariant algebraic cobordism, thus providing a vast generalization of the theory considered in \cite{levine-morel}. Our theory, and its associated cohomology theory, cannot be constructed employing motivic homotopy theory, as they are not $\Ab^1$-invariant (see e.g. \cite{deglise:2018} for bivariant theories in motivic homotopy theory). 

By proving the validity of projective bundle formula for the bivariant cobordism, we were able to construct cobordism Chern classes of vector bundles. This in turn, analogously to the well-known Conner--Floyd theorem \cite{conner:1966, conner:1969}, enabled us to recover algebraic $K$-theory rings---another non-$\Ab^1$-invariant theory of fundamental importance in algebraic geometry---from the cobordism cohomology rings, thus providing strong evidence that at least the algebraic cobordism cohomology theory we have constructed is the correct one. Moreover, this result allowed us to express the algebraic cobordism cohomology rings, taken with rational coefficients, in terms of algebraic $K$-theory.

\section{Open problems}

Even though some progress has been made in the articles \cite{annala-cob, annala-yokura, annala-chern, annala-pre-and-cob, annala-base-ind-cob, annala-spivak, annala-iwasa}, the study of the derived-geometric approach to algebraic cobordism is still in its early stages. Here, we list several interesting problems that remain open.

\begin{enumerate}
\item \emph{Localization.}  One of the most important computational tools in the study of Levine--Morel's algebraic bordism is the localization exact sequence. We expect that, for a regular Noetherian scheme $S$ of finite Krull dimension that admits an ample line bundle, and for a closed immersion $i: Z \hook X$ of derived $S$-schemes with open complement $j: U \hook X$, the sequence
$$\Omega^S_*(Z) \xto{i_*} \Omega^S_*(X) \xto{j^!} \Omega^S_*(U) \to 0$$
is exact. If $S = \Spec(k)$, where $k$ is a field of characteristic 0, then this is a theorem of Levine and Morel. The surjectivity of $j^!$ was proven in \cite{annala-spivak} with $\Zb[e^{-1}]$-coefficients in the case $S = \Spec(A)$, where $A$ is either a field or an excellent Henselian discrete valuation ring, and where $e$ is the (residual) characteristic exponent of $A$.

\item \emph{Comparison with $MGL$.} In \cite{levine:2009}, Levine constructs natural isomorphisms
$$\Omega_n(X) \xto \cong MGL'_{2n,n}(X)$$
for all quasi-projective varieties $X$ over a field of characteristic 0, where $MGL'_{*,*}$ is the motivic Borel--Moore homology theory associated to the algebraic cobordism spectrum $MGL$. In other words, the geometric algebraic bordism of Levine--Morel recovers a slice of ``higher algebraic bordism''. The same should be true over a more general regular basis scheme $S$. Note that this is a sensible expectation: the groups $\Omega^S_*(X)$ are expected to be $\Ab^1$-invariant, and, in fact, the $\Ab^1$-invariance was proven in \cite{annala-spivak} with $\Zb[e^{-1}]$-coefficients in the case $S = \Spec(A)$, where $A$ is either a field or an excellent Henselian discrete valuation ring, and where $e$ is the (residual) characteristic exponent of $A$. In the proof of \cite{levine:2009}, the localization exact sequence plays a fundamental role. However, it is conceivable that a geometric model of $MGL$, similar to that studied in \cite{EHKSY2} but which allows studying also the non-connective part of $MGL$, might lead to a direct proof that avoids the employment of the localization exact sequence.

\item \emph{Algebraic cobordism of a local ring.} It was proven by Levine and Morel that the algebraic cobordism of $\Spec(k)$, where $k$ is a field of characteristic 0, is naturally isomorphic to the Lazard ring $\Lb$. We have proven the analogous result for a general Noetherian local (derived) ring $A$ of finite Krull dimension, after taking rational coefficients. In \cite{annala-spivak}, it was proven that, for a field $k$ of characteristic $p > 0$, $\Omega^*_k(\Spec(k))[p^{-1}]$ and $\Lb^*[p^{-1}]$ coincide in degrees $\geq -2$. This follows from the fact, also proven in \cite{annala-spivak}, that $\Omega^*_k(\Spec(k))[p^{-1}]$ is generated by classes of regular projective $k$-schemes, combined with well-known resolution of singularities and birational factorization results for surfaces.

\item \emph{Comparison with Chow groups.} Levine and Morel prove that, for an algebraic scheme over a characteristic 0 field $k$, the natural map $\eta_\CH: \Zb_a \otimes_\Lb \Omega_*(X) \to \CH_*(X)$ defined by
$$[V \xto{f} X] \mapsto f_*(1_V),$$
is an isomorphism, where $\CH_*(X)$ is the Chow group of $X$ \cite{fulton:1998}. It would be interesting to establish a similar result over fields $k$ of characteristic $p > 0$ as well. With $\Zb[p^{-1}]$-coefficients, $\eta_\CH$ is a surjection, because every closed subvariety $Z$ of $X$ admits a generically finite morphism (\emph{resolution} by a $p$-alteration) of $p$-power degree from a regular $k$-variety \cite{dejong:1996, temkin:2017}. The main difficulty in proving that $\eta_\CH$ is an isomorphism with $p$-inverted coefficients is comparing the classes of different resolutions in $\Zb_a[p^{-1}] \otimes_\Lb \Omega^k_*(X)$, which in characteristic 0 can be achieved by employing the weak factorization theorem \cite{abramovich:2002}.

\item \emph{Comparison with candidates of Chow cohomology.} This is an open problem even over a field $k$ of characteristic 0, but of course the question makes sense over all fields $k$, and even over more general base schemes $S$. Namely, what is the relationship between the additive cobordism ring $\Omega^*_{k,a}(X)$ and the left-Kan-extended Chow rings \cite{fulton:1975}? The same question can be asked about the rings constructed by applying Bloch's formula to singular schemes \cite{quillen:1973}. If $X$ is of pure dimension $d$, then what is the relationship of $\Omega^d_{k,a}(X)$ with Levine--Weibel's group of zero-cycles on $X$ \cite{levine-weibel}?

\item \emph{Bivariant (or homological) Conner--Floyd theorem.} It would be interesting to extend the cohomological Conner--Floyd theorem to a bivariant Conner--Floyd theorem. In other words, for any finite-Krull-dimensional Noetherian derived scheme $S$ admitting an ample line bundle, we would expect that the natural map
$$\Zb_m \otimes_\Lb \Omega^*_{S}(X \to Y) \to K^0(X \to Y)$$
given by the formula
$$[V \xto{f} X] \mapsto [f_* \Oc_V]$$
is an isomorphism. Note that, using the twisting results of Chapter \ref{ch:bivtwist}, such a result would allow us to express the groups $\Qb \otimes \Omega^*_{S}(X \to Y)$ in terms of bivariant $K$-theory. This would be interesting even in the special case where $S$ is regular and $Y \simeq S$: for example, we could conclude that the the localization sequence with rational coefficients is exact\footnote{As $j$ is an open immersion, twisting the orientation does not affect Gysin pullbacks along it.}. Moreover, if $S = \Spec(k)$, where $k$ is a field, then it would follow that $\Qb_a \otimes \Omega^k_*(X) = \Qb \otimes \CH_*(X)$ \cite{baum:1975}.

\item \emph{Higher algebraic cobordism.} In this thesis, we have considered only what should be the truncation of a spectrum-valued invariant of maps derived schemes, the higher (bivariant) algebraic cobordism. Finding the correct model for higher algebraic cobordism is an interesting open problem: not only does it provide a more refined invariant, but there are also many properties of the spectrum-valued invariant that do not hold for the truncation, the most obvious being Nisnevich descent. In a joint project with Ryomei Iwasa, we study candidates of higher algebraic cobordism, which are motivated by the fact that universal precobordism is the universal bivariant theory satisfying the weak projective bundle formula, see \cite{annala-iwasa}. It is still too early to tell how successful this approach turns out to be.

The higher algebraic cobordism cohomology theory $\Omega$ is expected to satisfy at least the following properties:
\begin{enumerate}
\item $\Omega$ should be a Nisnevich sheaf of spectra on a suitable $\infty$-category of derived schemes;
\item $\Omega$ should satisfy the projective bundle formula;
\item $\Omega$ should satisfy a higher analogue of Conner--Floyd theorem, i.e., a non-$\Ab^1$-invariant analogue of \cite{panin:2008}; in particular, $\Omega$ should not be $\Ab^1$-invariant;
\item by applying a suitable $\Ab^1$-localization functor to $\Omega$, one should recover the algebraic cobordism theory represented by $MGL$ in the stable motivic homotopy category. 
\end{enumerate}

\item \emph{Excision in derived blowups}. Given a derived blowup square
$$
\begin{tikzcd}
\Ec \arrow[r]{}{i'} \arrow[d]{}{p} & \bl_Z(X) \arrow[d]{}{\pi^*} \\
Z \arrow[r]{}{i} & X,
\end{tikzcd}
$$
is the sequence
$$\PCob^*(X) \xto{i^* + \pi^*} \PCob^*(Z) \oplus \PCob^*(\bl_Z(X)) \xto{i'^* - p^*} \PCob^*(\Ec)$$
(or the analogous sequence for $\Omega^*$) exact?

There are two motivations for this question:
\begin{enumerate}
\item this is the analogue of the zero-level ``shadow'' of a spectrum-level result that holds for all additive invariants of derived categories (e.g. algebraic $K$-theory) \cite{khan:dcdh};

\item imposing an analogous spectrum-level property is useful when trying to construct higher algebraic cobordism; however, something like the exactness of the above sequence  would be needed for comparing the zeroth homotopy group of higher algebraic cobordism with $\PCob^*$ (or $\Omega^*$).
\end{enumerate}
\end{enumerate}

\begin{epigraph}
\emph{I have all the symptoms of fright....It really seems like I'm cut loose and very vulnerable....Still, I have a feeling of \emph{strength}....I'm feeling it internally now, a sort of surging up of force...something really big and strong. And yet at first it was almost a physical feeling of just being out \emph{alone}, and sort of cut off from a support I have been carrying around...(pause)...I have the feeling that now I am going to begin to \emph{do} more things.} ---~A client to Carl Rogers \cite{rogers:1980}.
\end{epigraph}

\begin{singlespace}
\raggedright
\bibliographystyle{abbrvnat}
\bibliography{../../References/references}
\end{singlespace}

\appendix

\backmatter

\end{document}